\documentclass[onecolumn,12pt]{elsarticle}

\usepackage{amssymb, bm, color}
\usepackage{amsthm, amsmath}

\usepackage{graphicx}
\usepackage{caption}
\usepackage{subcaption}
\usepackage{multirow}
\usepackage{mathtools,xparse}
\usepackage{algorithm}
\usepackage{algorithmic}

\makeatletter
\newenvironment{smallermatrix}[1][c]
{\null\,\vcenter\bgroup
  \Let@\restore@math@cr\default@tag
  \baselineskip0pt \lineskip0.4pt \lineskiplimit0pt
  \ialign\bgroup\if#1l\else\hfil\fi$\m@th\scriptstyle##$\if#1r\else\hfil\fi&&\thickspace\hfil
  $\m@th\scriptstyle##$\hfil\crcr
}{%
  \crcr\egroup\egroup\,%
}
\makeatother

\NewDocumentCommand{\ts}{O{c} e{^?_}}{
  \begin{smallermatrix}[#1]
  \mathstrut\IfValueT{#2}{#2} \\
  \mathstrut\IfValueT{#3}{#3} \\
  \mathstrut\IfValueT{#4}{#4}
  \end{smallermatrix}%
}

\newcommand{\revise}[1]{{\color{black}{#1}}}

\journal{arXiv}

\begin{document}

\begin{frontmatter}



\title{Optimal transport for mesh adaptivity and shock capturing of compressible flows}


\author[inst1]{Ngoc Cuong Nguyen}

\author[inst1]{R. Loek Van Heyningen}

\author[inst1]{Jordi Vila-P\'erez}

\author[inst1]{Jaime Peraire}

\affiliation[inst1]{organization={Center for Computational Engineering, Department of Aeronautics and Astronautics, Massachusetts Institute of Technology},
            addressline={77 Massachusetts
Avenue}, 
            city={Cambridge},
            state={MA},
            postcode={02139}, 
            country={USA}}

\begin{abstract}
We present an optimal transport approach for mesh adaptivity and shock capturing of compressible flows. Shock capturing is based on a viscosity regularization of the governing equations by introducing an artificial viscosity field as  solution of the Helmholtz equation. Mesh adaptation is based on the optimal transport theory by formulating a mesh mapping as solution of Monge-Amp\`ere equation. The marriage of optimal transport and viscosity regularization for compressible flows leads to a coupled system of the compressible Euler/Navier-Stokes equations, the Helmholtz equation, and the Monge-Amp\`ere equation. We propose an iterative  procedure to solve the coupled system in a sequential fashion using homotopy continuation to minimize the amount of artificial viscosity while enforcing positivity-preserving and smoothness constraints on the numerical solution. We explore various mesh monitor functions for computing r-adaptive meshes in order to reduce the amount of artificial dissipation and improve the accuracy of the numerical solution. The hybridizable discontinuous Galerkin method is used for the spatial discretization of the governing equations to obtain high-order accurate solutions. Extensive numerical results are presented to demonstrate the optimal transport approach on  transonic, supersonic,  hypersonic  flows in two dimensions. The approach is found to yield accurate, sharp yet smooth solutions within a few mesh adaptation iterations.  
\end{abstract}









\end{frontmatter}


\section{Introduction}
\label{sec:intro}

Compressible flows at high Mach number lead to shock waves which pose one of the most challenging problems for numerical methods. For high-order numerical methods, insufficient resolution or an inadequate treatment of shocks can result in Gibbs oscillations, which grow rapidly and contribute to numerical instabilities. Effective treatment of shock waves requires both  shock capturing  and mesh adaptation.

Shock capturing methods lie within one of the following two categories: limiters and artificial viscosity. Limiters, in the form of flux limiters \cite{Burbeau2001,Cockburn1989,Krivodonova2007}, slope limiters \cite{Cockburn1998a,MR2056921,Lv2015,Sonntag2017}, and WENO-type schemes \cite{Luo2007,Qiu2005,Zhu2008,Zhu2013} pose implementation difficulties for implicit time integration schemes and high-order methods on complex geometries. As for artificial viscosity methods, Laplacian-based \cite{Barter2010,Hartmann2013,Lv2016,Moro2016,Nguyen2011a,persson06:_shock_capturing,Persson2013} and physics-based \cite{Abbassi2014,Chaudhuri2017,Cook2004,Cook2005,Fernandez2018,Fiorina2007,Kawai2008,Kawai2010,Mani2009,Olson2013,persson06:_shock_capturing} approaches have been proposed. When the amount of viscosity is properly added in a neighborhood of shocks, the solution can converge uniformly except in the region around shocks, where it is smoothed and spread out over some length scale. Artificial viscosity has been widely used in finite volume methods~\cite{Jameson1995}, streamline upwind Petrov-Galerkin (SUPG) methods~\cite{HughesMalletMisukamiII86},  spectral methods~\cite{madtad93,Tadmor89}, as well as DG methods~\cite{Barter2010,Ching2019,HartmannHoustonCompressible02,persson06:_shock_capturing,Bai2022a,Vila-Perez2021}. Both Laplacian-based \cite{Hartmann2013,Lv2016,Nguyen2011a,Moro2016,persson06:_shock_capturing,Persson2013} and physics-based \cite{Abbassi2014,Bhagatwala2009,Mani2009,Cook2005,Cook2007,Fiorina2007,Kawai2008,Kawai2010,Olson2013,Premasuthan2013} artificial viscosity methods have been used for shock capturing.
 


Recent advances lead to shock-fitting methods that do not require limiter and artificial viscosity to stabilize shocks. The recent work \cite{Zahr2018,Zahr2020,Shi2022} introduces a high-order implicit shock tracking (HOIST) method for resolving discontinuous solutions of conservation laws with high-order numerical discretizations that support inter-element solution discontinuities, such as discontinuous Galerkin or finite volume methods \cite{Zahr2020}. The method aims to align mesh elements with shock waves by deforming the computational mesh in order to obtain accurate high-order solutions. It requires solution of a PDE-constrained optimization problem for both the computational mesh and the numerical solution using sequential quadratic programming
solver. Recently, a moving discontinuous Galerkin finite element method with interface condition enforcement (MDG-ICE) \cite{Corrigan2019,Kercher2021,Kercher2021a} is formulated for shock flows by enforcing the interface condition separately from the conservation laws. In the MDG-ICE method, the grid coordinates are treated as additional unknowns to detect interfaces and satisfy the interface condition, thereby directly fitting shocks and preserving high-order accurate solutions. The Levenberg-Marquardt method is used to solve the coupled system of the conservation laws and the interface condition to obtain the numerical solution and the shock-fitted mesh. 

In a recent work \cite{Nguyen2023c}, we introduce an adaptive viscosity regularization scheme for the numerical solution of nonlinear conservation laws with shock waves. The scheme solves a sequence of viscosity-regularized problems by using homotopy continuation to minimize the amount of viscosity subject to relevant physics and smoothness constraints on the numerical solution. The scheme is coupled to mesh adaptation algorithms that identify the shock location and generate shock-aligned meshes in order to further reduce the amount of artificial dissipation. In particular, shocks curves are constructed by determining shock-containing elements and finding a collection of points at which the artificial viscosity reaches its maximum value along streamline directions. A shock-aligned grid is generated by replicating shock curves along streamline directions. While the mesh alignment   procedure is simple, it is not practical for complex geometries and shock patterns.  

In this paper, we present an optimal transport approach for mesh adaptivity and shock capturing of compressible flows. Shock capturing is based on the viscosity regularization method introduced in \cite{Nguyen2023c}. Mesh adaptation is based on the optimal transport theory by formulating a mesh mapping as solution of Monge-Amp\`ere equation \cite{nguyen2023hybridizable}. The marriage of optimal transport and viscosity regularization for compressible flows leads to a coupled system of the compressible Euler/Navier-Stokes equations, the Helmholtz equation, and the Monge-Amp\`ere equation. We propose an iterative solution procedure to solve the coupled system in a sequential manner. We explore various mesh monitor functions for computing r-adaptive meshes in order to reduce the amount of artificial dissipation and improve the accuracy of the numerical solution. The hybridizable discontinuous Galerkin (HDG) method  is used for the spatial discretization of the governing equations owing to its efficiency and high-order accuracy \cite{Nguyen2012,Fernandez2018a,Moro2011a,Peraire2010,Vila-Perez2021,Woopen2014c,Fidkowski2016,Fernandez2017a,williams2018entropy}. 

Extensive numerical results are presented to demonstrate the proposed approach on a wide variety  of transonic, supersonic,  hypersonic  flows in two dimensions. The approach is found to yield accurate, sharp yet smooth solutions within a few mesh adaptation iterations. It is capable of moving mesh points to resolve complex shock patterns without creating new mesh points or modifying the  connectivity of the initial mesh. The generated r-adaptive mesh can significantly improve the accuracy of the numerical solution relative to the initial mesh. Accurate prediction of drag forces and heat transfer rates for viscous shock flows requires meshes to resolve both shocks and boundary layers. We show that the approach can generate r-adaptive meshes that resolve not only shocks but also boundary layers for viscous shock flows. The approach predicts  heat transfer coefficient accurately by adapting the initial mesh to resolve bow shock and boundary layer when it is applied to  viscous hypersonic flows past a circular cylinder.

The paper is organized as follows. We describe the adaptive viscosity regularization method in Section~\ref{sec:viscosity} and the optimal transport approach in Section 3. In Section~\ref{sec:results}, we present numerical results to assess the performance of the proposed approach on a wide variety of transonic, supersonic, and hypersonic flows. Finally, in Section~\ref{sec:conclusions}, we conclude the paper with some remarks and future work.

\section{Adaptive Viscosity Regularization Method} \label{sec:viscosity}

\subsection{Governing equations}

We consider the steady-state conservation laws of $m$  state variables, defined on a physical domain $\Omega \in \mathbb{R}^d$ and subject to appropriate  boundary conditions, as follows
\begin{equation}
\label{eq1}
 \nabla \cdot \bm F(\bm u, \nabla \bm u) = 0  \quad \mbox{in }\Omega,    
\end{equation}
where $\bm u(\bm x) \in \mathbb{R}^m$ is the solution of the system of conservation laws at $\bm x \in \Omega$ and the  physical fluxes $\bm F = (\bm f_1(\bm u, \nabla \bm u), \ldots, \bm f_d(\bm u, \nabla \bm u)) \in \mathbb{R}^{m \times d}$ include $d$ vector-valued functions of the solution. This paper focuses on the compressible Euler and Navier-Stokes equations. 

For the compressible Euler equations, the state vector and physical fluxes are given by
\begin{equation}
    \bm u = \begin{pmatrix}
        \rho \\
        \rho v_i \\
        \rho E \\
    \end{pmatrix}, \qquad  
    \bm F(\bm u) = \begin{pmatrix}
        \rho v_j \\ 
        \rho v_i v_j + \delta_{ij} p \\ 
        \rho v_j H
    \end{pmatrix}
    \label{eq:euler}
\end{equation}
with density $\rho$, velocity $\bm v$, total energy $E$, total specific enthalpy $H = E + p/\rho$ and pressure $p$ given by the ideal gas law $p = (\gamma - 1)  \rho ( E-\frac{1}{2}v_i \ v_i)$. Let $\Gamma_{\rm wall} \subset \partial \Omega$ be the wall boundary. The boundary condition at the wall boundary $\Gamma_{\rm wall}$ is $\bm v \cdot \bm n = 0$, where $\bm v$ is the velocity field and $\bm n$ is the unit normal vector outward the boundary. For supersonic and hypersonic flows, supersonic inflow and outflow conditions are imposed on the inflow and outflow boundaries, respectively. For transonic flows, a freastream boundary condition is imposed at the far field boundary by using the freestream state $\bm u_\infty$. The freestream Mach number $M_\infty$ enters through the non-dimensional freestream pressure $p_{\infty} = 1 / (\gamma M^2_{\infty})$. Here $\gamma$ denotes the specific heat ratio.

For the compressible Navier-Stokes equations, the fluxes are given by
\begin{equation}
\label{flux}
\bm{F}(\bm{u},\nabla \bm u) = \left( \begin{array}{c}
\rho v_j \\
\rho v_i v_j + \delta_{ij} p \\
 \rho v_j H
\end{array}
\right) - \left( \begin{array}{c}
0 \\
\tau_{ij}  \\
v_i \tau_{ij} + f_j
\end{array}
\right) .
\end{equation}
For a Newtonian, calorically perfect gas in thermodynamic equilibrium, the non-dimensional viscous stress tensor and heat flux are given by
\begin{equation}
\label{closures}
\tau_{ij} = \frac{1}{Re} \bigg[ \Big( \frac{\partial v_i}{\partial x_j}+\frac{\partial v_j}{\partial x_i} \Big) -\frac{2}{3}\frac{\partial v_k}{\partial x_k}\delta_{ij} \bigg] ,  \qquad f_j = - \frac{\gamma}{Re \ Pr} \ \frac{\partial T}{\partial x_j} ,
\end{equation}
respectively. Here $Re$ denotes the Reynolds number, and $Pr$ the Prandtl number. For high Mach number flows, Sutherland's law is used to obtain the dynamic viscosity, thereby rendering the Reynolds number dependent on the temperature. The boundary conditions at the wall are zero velocity and either isothermal or adiabatic temperature. Other boundary conditions are similar to those of the compressible Euler equations.


\subsection{Viscosity regularization of compressible flows}



Shock waves have always been a considerable source of difficulties toward a rigorous numerical solution of compressible flows. In order to treat shock waves, we follow the recent work \cite{Barter2010,Ching2019} by considering the  viscosity regularization of the conservation laws (\ref{eq1}) as follows
\begin{subequations}
\label{eq3}
\begin{alignat}{2}
 \nabla \cdot \bm F(\bm u, \nabla \bm u) - \lambda_1 \nabla \cdot \bm G(\bm u,   \nabla \bm u,  \eta) = 0  \quad \mbox{in }\Omega, \\
\eta - \lambda_2^2 \nabla  \cdot \left( 
\ell^2 \nabla  \eta  \right)- s(\bm u, \nabla \bm u) = 0 \quad \mbox{in }\Omega ,
\end{alignat}
\end{subequations}
where $\eta(\bm x)$ is the solution of the Helmholtz  equation (\ref{eq3}b) with homogeneous  Neumann boundary conditions 
\begin{equation}
\eta = 0 \quad \mbox{on } \Gamma_{\rm wall}, \qquad 
\ell^2 \nabla  \eta \cdot \bm n = 0 \quad \mbox{on }\partial \Omega \backslash \Gamma_{\rm wall} \ .
\end{equation}
\revise{Here $\lambda_1$ is the first regularization parameter that controls the amplitude of artificial viscosity, and $\lambda_2$ is the second regularization parameter that controls the thickness of artificial viscosity. Furthermore, $\ell$ is an appropriate length scale which is chosen as the smallest mesh size $h_{\rm min}$. For notational convenience, we denote $\bm \lambda = (\lambda_1, \lambda_2)$.}


The artificial fluxes $\bm G$ provide a viscosity regularization  to smooth out discontinuities in the shock region. There are a number of different options for the artificial fluxes $\bm G$. In this paper, we use the Laplacian fluxes of the form
\begin{equation}
\bm G(\bm u,   \nabla \bm u,  \eta) = \mu(\eta) \nabla \bm u  ,  
\end{equation} 
where 
\begin{equation}
\mu(\eta) = (\bar{\eta}-\bar{\eta}_{\rm T})\left(\frac{\arctan(100(\bar{\eta} -\bar{\eta}_{\rm T}))}{\pi} + \frac{1}{2} \right) - \frac{\arctan(100)}{\pi} + \frac{1}{2}  
\end{equation} 
is a smooth approximation of a ramp function. Here $\bar{\eta} = \eta/\|\eta\|_\infty$ is the normalized function with $\|\eta\|_{\infty} = \max_{\bm x \in \Omega}  |\eta(\bm x)|$ being the $L_\infty$ norm. Note that $\bar{\eta}_{\rm T}$ is the artificial viscosity threshold that makes $\mu(\eta)$ vanish to zero when $\bar{\eta} \le \bar{\eta}_{\rm T}$. In other words, artificial viscosity is only added to the shock region where $\bar{\eta}$ exceeds $\bar{\eta}_{\rm T}$. Therefore, the threshold $\bar{\eta}_{\rm T}$ will help remove excessive artificial viscosity. Since $\|\bar{\eta}\|_\infty = 1$, $\bar{\eta}_{\rm T} = 0.2$ is a sensible choice. Note that the artificial viscosity field is equal to $\lambda_1 \mu(\bm x)$, where  $\mu(\bm x)$ is bounded by $\mu(\bm x) \in [0, 1 - \bar{\eta}_{\rm T}]$ for any $\bm x \in \Omega$. We can also consider a more general form $\bm G = \mu(\eta) \nabla \bm u^*$ \cite{Barter2010,Nguyen2011a}, where $\bm u^*$ is a modified state vector.  Another option is physics-based artificial viscosity  by taking $\bm G$ to be the viscous stress tensor and the heat flux of the Navier-Sokes equation and adding the artificial viscosity to the physical viscosities and thermal conductivity \cite{CuongNguyen2022,Nguyen2023a}. 

 
The source term $s$ in (\ref{eq3}b) is required to determine $\eta$ and defined as follows
\begin{equation}
\label{avsource}
s(\bm u, \nabla \bm u) =  {g}(S(\bm u, \nabla \bm u)) 
\end{equation}
where  $g(S)$ is a smooth approximation of the following step function 
\begin{equation}
\label{eq8g}
\tilde{g}(S) = \left\{
\begin{array}{cl}
   0  & \mbox{if } S < 0, \\
   S  & \mbox{if } 0 \le S \le s_{\rm max}, \\
   s_{\rm max} & \mbox{if } S > s_{\rm max} . 
\end{array}
\right.
\end{equation}
The quantity $S(\bm u, \nabla \bm u)$ is a measure of the shock strength which is given by
\begin{equation}
S(\bm u, \nabla \bm u) = -\nabla \cdot \bm v \ ,
\end{equation}
where $\bm v$ is the non-dimensional velocity field that is determined from the state vector $\bm u$.  The use of the velocity divergence as shock strength for defining an artificial viscosity field follows from \cite{Fernandez2018,Moro2016,Nguyen2011a}. The parameter $s_{\max}$ is used to put an upper bound on the source term when the divergence of the velocity becomes too negatively large. Herein we choose $s_{\max} = 0.5 \|S\|_\infty$, where $\|S\|_{\infty} = \max_{\bm x \in \Omega}  |S(\bm x)|$ is the $L_\infty$ norm. \revise{Since $S$ depends on the solution, so its norm may not be known prior. In practice, we employ a homotopy continuation scheme to  iteratively solve the problem (\ref{eq3}).  Hence, $s_{\max}$ is computed by using the numerical solution at the previous iteration of the homotopy continuation.} 



It remains to determine $\lambda_1$ and $\lambda_2$ in order to close the system (\ref{eq3}). We propose to solve the following minimization problem
\begin{subequations}
\label{eq5}
\begin{alignat}{2}
\min_{\lambda_1 \in \mathbb{R}^+, \lambda_2 \ge 1, \bm u, \eta} & \quad \lambda_1 \lambda_2  \\
\mbox{s.t.} & \quad \mathcal{L}(\bm u, \eta, \bm \lambda) = 0 \\
 & \quad \bm u \in \mathcal{C} .
\end{alignat}
\end{subequations}
Here $\mathcal{L}$ represents the spatial discretization of the coupled system (\ref{eq3}) by a numerical method and $\mathcal{C}$ represents a set of constraints on the numerical solution. The objective function is to minimize the amount of artificial viscosity which is proportional to $\lambda_1 \lambda_2$. The constraints are specified to rule out unwanted solutions of the discrete system (\ref{eq5}b) and play an important role in yielding a high-quality numerical solution. Hence, the optimization problem (\ref{eq5}) is to minimize the amount of artificial viscosity while ensuring the smoothness of the numerical solution. 

\subsection{Solution constraints}

We  introduce the constraints to ensure the quality of the numerical solution. The physical constraints are that pressure and density must be positive. In order to establish a smoothness constraint on the numerical solution, we express an approximate scalar variable $\xi$ of degree $k$ within each element in terms of an orthogonal basis and its truncated expansion of degree $k-1$ as
\begin{equation}
\xi = \sum_{i=1}^{N(k)} \xi_i \psi_i, \qquad \xi^* = \sum_{i=1}^{N(k-1)} \xi_i \psi_i    
\end{equation}
where $N(k)$ is the total number of terms in the $k$-degree expansion and $\psi_i$ are the basis functions \cite{persson06:_shock_capturing}. Here $\xi$ is chosen to be either density, pressure, or local Mach number.  We introduce the following quantity
\begin{equation}
\label{eq10}
\sigma(\bm \lambda) = \max_{K \in \mathcal{T}_h^{\rm shock}} \sigma_K(\bm \lambda), \qquad  \sigma_K(\bm \lambda) \equiv   \frac{\int_K  |\xi/\xi^* - 1| d \bm x}{\int_K d \bm x} ,
\end{equation}
where $\mathcal{T}_h^{\rm shock}$ is the set of elements defining the shock region
\begin{equation}
\mathcal{T}_h^{\rm shock} = \{K \in \mathcal{T}_h \ : \ \int_K  \bar{\eta}  d \bm x  \ge \bar{\eta}_{\rm T} |K|  \} 
\end{equation}
and $\mathcal{T}_h$ is a collection of  high-order elements on the physical domain $\Omega$ 
\begin{equation}
\label{homesh}
\mathcal{T}_h = \{K_n \in \Omega \ : \cup_{n=1}^{N_e} \bar{K}_n = \bar{\Omega}, \bm x|_{K_n} \in [\mathcal{P}_k(K_{\rm ref})]^d, 1 \le n \le N_{\rm e}\} .
\end{equation}
Here $N_e$ is the number of elements, $K_{\rm ref}$ is the master element, and $\mathcal{P}_k(K_{\rm ref})$ is the space of polynomials of degree $k$ on $K_{\rm ref}$.  The constraint set $\mathcal{C}$ in (\ref{eq5}) consists of the following contraints    
\begin{equation}
 \rho(\bm x) > 0, \quad  p(\bm x) > 0,   \quad   \sigma(\bm \lambda)  \le C_0 \, \sigma(\bm \lambda_0) ,
\end{equation}
where $\bm \lambda_0$ is an initial value and the constant $C_0$ is set to 5. The first two constraints enforce the positivity of density and pressure, while the last constraint guarantees the smoothness of the numerical solution. The smoothness constraint  imposes a degree of regularity on the numerical solution and plays a vital role in yielding sharp and smooth solutions. 

 
 \subsection{Homotopy continuation of the regularization parameters}

The pair of regularization  parameters $\bm \lambda = (\lambda_1, \lambda_2)$ controls the magnitude and thickness of the artificial viscosity in order to obtain accurate solutions. On the one hand, if $\bm \lambda$ is too small then the numerical solution can develop oscillations across the shock waves. On the other hand, if $\bm \lambda$ is too large the solution becomes less accurate in the shock region, which in turn  affects the accuracy of the solution in the remaining region. Therefore, we propose a homotopy continuation method to determine $\bm \lambda$. The key idea is to solve the regularized system with a large value of $\bm \lambda$ first and then gradually decrease $\bm \lambda$ until any of the physics or smoothness constraints on the numerical solution are violated. At this point, we take the value of $\bm \lambda$ from the previous iteration where the numerical solution still satisfies all of the physics and smoothness constraints. This procedure is summarized in the following algorithm:


\revise{
\begin{itemize}
 \item Given an initial value $\bm \lambda_0 = (\lambda_{0,1}, \lambda_{0,2})$ and $\eta_0$ such that $\|\eta_0\|_{\infty} = 1$, solve the regularized system  (\ref{eq3}a) with $\lambda_1 = \lambda_{0,1}, \eta = \eta_0$ to obtain the initial solution $\bm u_0$. 
 \item   Set $\lambda_{n,1} = \zeta^{n-1} \lambda_{{n-1},1}$ and $\lambda_{n,2} = 1 + \zeta^{n-1}(\lambda_{{n-1},2} -1)$ for some constant $\zeta \in (0,1)$; solve the Helmholtz equation (\ref{eq3}b) with $\lambda_2 = \lambda_{n,2}$ and the source term from $\bm u_{n-1}$ to obtain $\eta_{n}$; and solve the regularized system (\ref{eq3}a) with $\lambda_1 = \lambda_{n,1}, \eta = \eta_{n}$ to obtain the solution $\bm u_{n}$ for $n = 1, 2, \ldots$ until $\bm u_{n}$ violates any of the constraints.
 \item Finally, we accept $\bm u_{n-1}$ as the numerical solution of the compressible Euler/Navier-Stokes equations.
\end{itemize}

The initial function $\eta_0$ can be set to 1 on most of the physical domain $\Omega$ except near the wall boundary where it vanishes smoothly to zero at the wall. The initial value $\lambda_{0,1}$ is conservatively large to make the initial solution $\bm u_0$ very smooth. The initial value $\lambda_{0,2}$ depends on the type of meshes used to compute the numerical solution. For regular meshes that have the elements of the same size in the shock region, $\lambda_{0,2} = 1.5$ is a sensible choice. For adaptive meshes that are refined toward the shock region, we choose $\lambda_{0,2} = 5$ since $\ell = h_{\min}$ is extremely small for shock-adaptive meshes. In any case, $\lambda_{n,2}$ will decrease from $\lambda_{0,2}$ toward 1 during the homotopy iteration. Hence, the choice of $\lambda_{0,2}$ can be flexible.
}


\subsection{Finite element approximations} 
\label{sec:FEmethods}

The homotopy continuation solves the Helmholtz equation (\ref{eq3}b) separately from the regularized system (\ref{eq3}a). Hence, different  numerical methods can be used to solve (\ref{eq3}a) and (\ref{eq3}b) separately. In this paper, we employ the hybridizable discontinuous Galerkin (HDG) method to solve the former and the continuous Galerkin (CG) method to solve the latter. We use the CG method since it allows us to obtain a continuous artificial viscosity field. The HDG method \cite{Nguyen2012,Fernandez2018a,Moro2011a,Peraire2010,Vila-Perez2021,Woopen2014c,Fidkowski2016,Fernandez2017a,williams2018entropy} is suitable for solving the regularized conservation laws because of its efficiency and high-order accuracy.

\section{Mesh adaptation via optimal transport} \label{sec:ot}

\subsection{Optimal transport theory}

The optimal transport (OT) problem is described as follows. Suppose we are given two probability densities: $\varrho(\bm x)$ supported on $\Omega \in \mathbb{R}^d$ and $\varrho'(\bm x')$ supported on $\Omega' \in \mathbb{R}^d$. The source density $\varrho(\bm x)$ may be discontinuous and  even vanish. The target density $\varrho'(\bm x')$ must be strictly positive and Lipschitz continuous. The OT problem is to find a map $\bm \phi : \Omega \to \Omega'$ such that it minimizes the following functional
\begin{equation}
\inf_{\bm \phi \in \mathcal{M}} \int_{\Omega} \|\bm x - \bm \phi(\bm x)\|^2 \varrho(\bm x) d \bm x ,     
\end{equation}
where 
\begin{equation}
\label{otset}
\mathcal{M} = \{\bm \phi : \Omega \to \Omega', \ \varrho'(\bm \phi(\bm x)) \det (\nabla \bm \phi(\bm x)) = \varrho(\bm x), \ \forall \bm x \in \Omega  \} ,
\end{equation}
is the set of mappings which map the source density $\varrho(\bm x)$ onto the target density $\varrho'(\bm x')$.  Here $\det$ denotes determinants for $d \times d$ matrices. Whenever the infimum is achieved by some map $\bm \phi$, we say that $\bm \phi$ is an optimal map.

In \cite{Brenier1991}, Brenier gave the proof of the existence and uniqueness of the solution of the OT problem. Furthermore, the optimal map $\bm \phi$ can be written as the gradient  of a unique (up to a constant) convex potential $u$, so that $\bm \phi(\bm x) = \nabla u(\bm x)$, $\Delta u(\bm x) > 0$. Substituting $\bm \phi(\bm x) = \nabla u(\bm x)$ into (\ref{otset}) results in the Monge–Amp\`ere equation
\begin{equation}
\label{mae}
 \varrho'(\nabla u(\bm x)) \det (D^2 u(\bm x)) = \varrho(\bm x) \quad \mbox{in } \Omega, 
\end{equation}
along with the restriction that $u$ is convex.  This equation lacks standard boundary conditions. However, it is geometrically constrained by the fact that the gradient map takes $\partial \Omega$ to $\partial \Omega'$: 
\begin{equation}
\label{maebc}
 \nabla u(\bm x) \in \partial \Omega', \quad \forall \bm x \in \Omega .
\end{equation}
This constraint is referred to as the second boundary value problem for the Monge–Amp\`ere equation. If the boundary $\partial \Omega'$ can be expressed by 
$$\partial \Omega' = \{\bm x' \in \Omega' : c(\bm x') = 0\},$$ 
then the  boundary constraint (\ref{maebc}) becomes the following Neumann boundary condition
\begin{equation}
\label{maenm}
 c(\nabla u(\bm x)) = 0, \quad \forall \bm x \in \partial \Omega .
\end{equation}
The scalar potential $u$ is required to satisfy $\int_{\Omega} u(\bm x) d \bm x = 0$ for uniqueness. For problems where densities are periodic, it is natural and convenient to use periodic boundary conditions instead. 
\subsection{Equidistribution principle}

 Mesh adaptation is based on the equidistribution principle that equidistributes the target density function $\varrho'$ so that the source density $\varrho$ is uniform on $\Omega$ \cite{Delzanno2008,Chacon2011}. The equidistribution principle leads to a constant source density $\varrho(\bm x) = \theta$, where $\theta = \int_{\Omega'} \varrho'(\bm x') d \bm x' / \int_{\Omega} d \bm x$. Using the optimal transport theory, the optimal map is sought by solving the  Monge–Amp\`ere equation:
\begin{equation}
\label{maem}
\begin{split}
\varrho'(\nabla u(\bm x))  \det (D^2 u(\bm x)) & = \theta, \quad \mbox{in } \Omega, \\
 c(\nabla u(\bm x)) & = 0, \quad  \mbox{on } \partial \Omega ,
 \end{split}
\end{equation}
with the constraint $\int_{\Omega} u(\bm x) d \bm x = 0$. In the context of mesh adaptation, the target  boundary $\partial \Omega'$ coincides with  $\partial \Omega$. Hence, the root of the equation $c(\bm x) = 0$ defines $\partial \Omega$. It means that $c(\nabla u(\bm x)) = c(\bm x) = 0, \forall \bm x \in \partial \Omega$.

\subsection{Mesh density function}

In the context of mesh adaptation, $\varrho'(\bm x')$ is the mesh density function and $\mathcal{T}_h$ is the initial mesh. The optimal map $\bm \phi(\bm x) = \nabla u(\bm x)$ drives the coordinates of the initial mesh to concentrate around a region where the mesh density function is high. Therefore, we need to make $\varrho'(\bm x')$ large in the shock region and small in the smooth region. It is also necessary for $\varrho'(\bm x')$ to be sufficiently smooth, so that the numerical approximation of the Monge–Amp\`ere equation (\ref{maem}) is convergent. To this end, we compute $\varrho'(\bm x')$ as solution of the Helmholtz equation 
\begin{equation}
\label{hmmd}
\varrho'(\bm x') - \nabla \cdot \left( \ell^2 \nabla \varrho'(\bm x') \right) = b(\bm x') \quad \mbox{in }\Omega, 
\end{equation}
with homogeneous Neumann boundary condition. Here $b$ is a resolution indicator function that is large in the shock region and small elsewhere.

We explore two different options for the indicator function. The first option is to define it as a function of the velocity divergence as
\begin{equation}
b(\bm x') = \sqrt{1 + \beta s(\bm u, \nabla \bm u)}
\end{equation}
where $s$ is given by (\ref{avsource}) and $\beta$ is a specified constant. The second option is to define it as a function of the density gradient as
\begin{equation}
\label{mdf2}
b(\bm x') = \sqrt{1 + \beta g(|\nabla \rho(\bm x')|) } 
\end{equation}
where $g(\cdot)$ is given by (\ref{eq8g}). Other indicator functions are possible, such as those based on some combination of physics-based sensors that can distinguish between shocks, large temperature gradients, and other sharp features. 

The Helmholtz equation (\ref{hmmd}) is numerically solved by using the CG method in which the same polynomial spaces are used to represent both the numerical solution and the geometry. In this case, the value of the mesh density function at any given point $\bm x' \in K_i \subset \mathcal{T}_h$ is calculated as
\begin{equation}
\label{mdp1}
\varrho'(\bm x')|_{K_i} = \sum_{j=1}^{N_p} \rho_{ij} \varphi_j(\bm \xi(\bm x')) 
\end{equation}
where $\bm \xi(\bm x')$ is found by solving the following system
\begin{equation}
\label{mdp2}
\sum_{j=1}^{N_p} \bm x_{ij} \varphi_j(\bm \xi) = \bm x' .
\end{equation}
Here $N_p$ is the number of polynomials per element, $\bm x_{ij}$ are the mesh nodes on element $K_n$, $\rho_{ij}$ are the degrees of freedom of the function $\rho'$ on $K_n$, and $\varphi_j(\bm \xi), 1 \le j \le N_p,$ are polynomials of degree $k$ defined on the master element $K_{\rm ref}$. We note that the system (\ref{mdp2}) is linear for $k=1$ and nonlinear for $k > 1$. 

The mesh density function is the numerical solution of the Helmholtz equation (\ref{hmmd}) whose source term depends on the flow state $\bm u$. In practice, we compute the approximate solution of the flow state by using the adaptive viscosity regularization method to solve the problem (\ref{eq5}) on the initial mesh $\mathcal{T}_h$ or on the previous adaptive mesh during the mesh adaptation procedure described in subsection 3.5.

\subsection{Numerical solution of the Monge–Amp\`ere equation}

In a recent paper \cite{nguyen2023hybridizable}, we introduce HDG methods for numerically solving the Monge–Amp\`ere equation in which the mesh density function is an analytical function. In order to solve the Monge–Amp\`ere equation in which the mesh density function is approximated by local spaces of polynomials in (\ref{mdp1})-(\ref{mdp2}), we propose to extend the HDG methods introduced in \cite{nguyen2023hybridizable}.

In two dimensions, the Monge-Amp\`ere equation (\ref{maem}) can be rewritten as a first-order system of equations
\begin{equation}
\begin{array}{rcll}
\bm{H} - \nabla \bm q & = & 0, \quad & \mbox{in } \Omega, \\
\bm{q} - \nabla u & = & 0, \quad & \mbox{in } \Omega, \\
f(\bm H, \bm q) - \nabla \cdot \bm q & = & 0 , \quad & \mbox{in }  \Omega, \\
c(\bm q) & = & 0, \quad & \mbox{on } \partial \Omega , \\
\int_{\Omega} u(\bm x) d \bm x &= & 0,
\end{array}
\label{maem3}
\end{equation}
where $f(\bm H,\bm q) = \sqrt{H_{11}^2 + H_{22}^2 + H_{12}^2 + H_{21}^2 + 2 \theta/\varrho'(\bm q)}$. The HDG discretization of the system (\ref{maem3}) is to find $(\bm{H}_h, \bm{q}_h,u_h,\widehat{u}_h) \in \bm{W}_{h}^p \times \bm{V}_{h}^p  \times U_{h}^p \times M_h^p$ such that
\begin{equation}
\label{eq82}
\begin{array}{rcl}
 \left(\bm{H}_h, \bm{G}\right)_{\mathcal{T}_h} +  \left(\bm q_h,  \nabla \cdot \bm{G}\right)_{\mathcal{T}_h} -  \left\langle \widehat{\bm q}_h, \bm{G} \cdot \bm{n} \right\rangle_{\partial \mathcal{T}_h} & = & 0,  \\
 \left(\bm{q}_h, \bm{v}\right)_{\mathcal{T}_h} +  \left(u_h,  \nabla \cdot \bm{v}\right)_{\mathcal{T}_h} -  \left\langle \widehat{u}_h, \bm{v} \cdot \bm{n} \right\rangle_{\partial \mathcal{T}_h} & = & 0,  \\
\left(\bm{q}_h, \nabla w\right)_{\mathcal{T}_h}  -  \left\langle \widehat{\bm{q}}_h \cdot \bm{n}, w \right\rangle_{\partial \mathcal{T}_h} + (f(\bm H_h, \bm q_h),w)_{\mathcal{T}_h} & = & 0, \\
\left\langle \widehat{\bm{q}}_h  \cdot \bm{n} , \mu \right\rangle_{\partial \mathcal{T}_h \backslash \partial \Omega} + \left\langle c(\bm q_h) + \tau (\widehat{u}_h - u_h), \mu \right\rangle_{\partial \Omega} & =  & 0, \\
(u_h,1)_{\mathcal{T}_h} & = & 0,
\end{array}
\end{equation}
for all $(\bm G, \bm{v}, w, \mu) \in  \bm{W}_h^p \times 
 \bm{V}_h^p \times U_h^p \times M_h^p$, where
\begin{equation}
\widehat{\bm{q}}_h = {\bm{q}_h} -
\tau (u_h - \widehat{u}_h) \bm{n}, \quad \mbox{on } \mathcal{E}_h.
\label{fluxdef2}
\end{equation}
We are going to use the fixed point method to solve this nonlinear system of equations.

To deal with the nonlinear boundary condition $c(\bm q) = 0$, we linearize it around the previous solution $\bm q^{\ell -1}$ to obtain 
\begin{equation}
c(\bm q^{l-1}) + \partial c_{\bm q}(\bm q^{l-1}) \cdot \left( \bm q^{l} - \bm q^{l-1} \right) = 0 ,
\end{equation}
where $\partial c_{\bm q}$ denotes the partial derivative of $c$ with respect to $\bm q$. Starting from an initial guess $(\bm H^0_h, \bm q_h^0, u_h^0)$  we  find $(\bm{q}_h^l,u_h^l,\widehat{u}_h^l) \in \bm{V}_{h}^k  \times U_{h}^k \times M_h^k$ such that
\begin{equation}
\label{nfpHDG}
\begin{array}{rcl}
  \left(\bm{q}_h^l, \bm{v}\right)_{\mathcal{T}_h} +  \left(u_h^l,  \nabla \cdot \bm{v}\right)_{\mathcal{T}_h} -  \left\langle \widehat{u}_h^l, \bm{v} \cdot \bm{n} \right\rangle_{\partial \mathcal{T}_h} & = & 0,  \\
\left(\bm{q}_h^l, \nabla w\right)_{\mathcal{T}_h}  -  \left\langle \widehat{\bm{q}}_h^l \cdot \bm{n}, w \right\rangle_{\partial \mathcal{T}_h}  & = & - (f(\bm H_h^{l-1}, \bm q_h^{l-1}),w)_{\mathcal{T}_h}, \\
\left\langle \widehat{\bm{q}}_h^l  \cdot \bm{n} , \mu \right\rangle_{\partial \mathcal{T}_h \backslash \partial \Omega}  + \left\langle \partial c_{\bm q}(\bm q_h^{l-1}) \cdot \bm q_h^l + \tau ( \widehat{u}_h^l -  u_h^l) , \mu \right\rangle_{\partial \Omega} & =  &  -\left\langle  a(\bm q^{l-1}) , \mu \right\rangle_{\partial \Omega}, \\
(u_h^l,1)_{\mathcal{T}_h} & = & 0, 
\end{array}
\end{equation}
for all $(\bm{v}, w, \mu) \in  \bm{V}_h^k \times U_h^k \times M_h^k$, and then compute $\bm{H}_h^l\in \bm{W}_{h}^k$ such that
\begin{equation}
\label{neqH}
 \left(\bm{H}_h^l, \bm{G}\right)_{\mathcal{T}_h} =  - \left(\bm q_h^l,  \nabla \cdot \bm{G}\right)_{\mathcal{T}_h} +  \left\langle \widehat{\bm q}_h^l, \bm{G} \cdot \bm{n} \right\rangle_{\partial \mathcal{T}_h}, \quad \forall \bm G \in \bm{W}_h^k .
\end{equation}
Note here that $a(\bm q_h^{l-1}) = c(\bm q_h^{l-1}) - \partial c_{\bm q}(\bm q_h^{l-1}) \cdot  \bm q_h^{l-1}$, and that the numerical flux $\widehat{\bm q}_h^l$ is defined by (\ref{fluxdef2}). We refer to \cite{nguyen2023hybridizable} for the definition of the finite element spaces  associated with the fixed-point HDG formulation (\ref{nfpHDG})-(\ref{neqH}) and the detailed implementation.

At each iteration of the fixed-point HDG method, the weak formulation (\ref{nfpHDG}) yields a matrix system which can be solved efficiently by locally eliminating the degrees of freedom of $(\bm{q}_h^l,u_h^l)$ to obtain a global linear system in terms of the degrees of freedom of $\widehat{u}_h^l$. While it is straightforward to form the matrix, computing the right-hand side vector is more complicated because we need to evaluate $f(\bm H_h^{l-1}, \bm q_h^{l-1})$. Henceforth, we must compute $\varrho'(\bm q_h^{l-1})$ by replacing $\bm x'$ with $\bm q_h^{l-1}$ in (\ref{mdp1}) and solve the resulting system (\ref{mdp2}) by using Newton's method for all quadrature points.


\subsection{Mesh adaptation procedure}

We start mesh adaptation with an initial mesh $\mathcal{T}_h$ and compute the initial solution $\bm u_h$. Next, we compute a mesh density function based on $\bm u_h$ and solve the Monge-Amp\`ere equation to obtain an adaptive mesh $\mathcal{T}_h^*$. Finally, we interpolate $\bm u_h$ onto $\mathcal{T}_h^*$ and use it as an initial guess to solve for the final solution $\bm u_h^*$ on the adaptive mesh. The mesh adaptation procedure is described in Algorithm 1. The adaptation procedure can be repeated by using the adaptive mesh as an initial mesh in the next iteration until $\|\bm u_h^* - \bm u_h\|_{\Omega}$ is less than a specified tolerance. It should be pointed out that we do not perform the homotopy continuation at every mesh adaptation iterations. We perform the homotopy continuation to compute the numerical solution on the final adaptive mesh only. This will considerably reduce the number of times we solve the compressible Euler/Navier-Stokes equations.

\begin{algorithm}
\begin{algorithmic}[1]
\REQUIRE{The initial mesh $\mathcal{T}_h$.}
\ENSURE{The r-adaptive mesh $\mathcal{T}_h^{*}$  and the numerical solution $\bm u_h^*$ on $\mathcal{T}_h^{*}$.}
\STATE{Solve (\ref{eq5}) for $\bm u_h$ on $\mathcal{T}_h$ using the adaptive viscosity regularization method.}
\STATE{Compute the mesh density $\varrho'_h(\bm x')$ based on $\bm u_h$ by solving  \eqref{hmmd}.}
\STATE{Solve the Monge-Amp\`ere equation \eqref{maem3} on $\mathcal{T}_h$ using the fixed-point HDG method.}
\STATE{Average $\bm q_h$ at duplicate nodes to obtain the adaptive mesh $\mathcal{T}_h^{*}$.}
\STATE{Interpolate $\bm u_h$ onto $\mathcal{T}_h^{*}$ and use it as the initial guess.}
\STATE{Solve (\ref{eq5}) for $\bm u_h^*$ on $\mathcal{T}_h^{*}$ using the adaptive viscosity regularization method.}
\end{algorithmic}
\caption{Mesh adaptation procedure.}
\end{algorithm}

We demonstrate the action of Algorithm 1 by applying it to an inviscid supersonic flow in a channel with a 4\% thick circular bump \cite{Nguyen2011d}. The length and height of the channel are 3 and 1, respectively. The inlet Mach number is $M_\infty = 1.4$. Supersonic inlet/outlet conditions are prescribed at the left/right boundaries, while inviscid wall boundary condition is used on the top and bottom sides. Isoparametric elements with the polynomials of degree $k = 4$ are used to represent both the numerical solution and geometry. Representative inputs and outputs are shown in Figure 1.

\begin{figure}[ht]
\centering
	\centering
	\begin{subfigure}[b]{0.49\textwidth}
		\centering		
  \includegraphics[width=\textwidth]{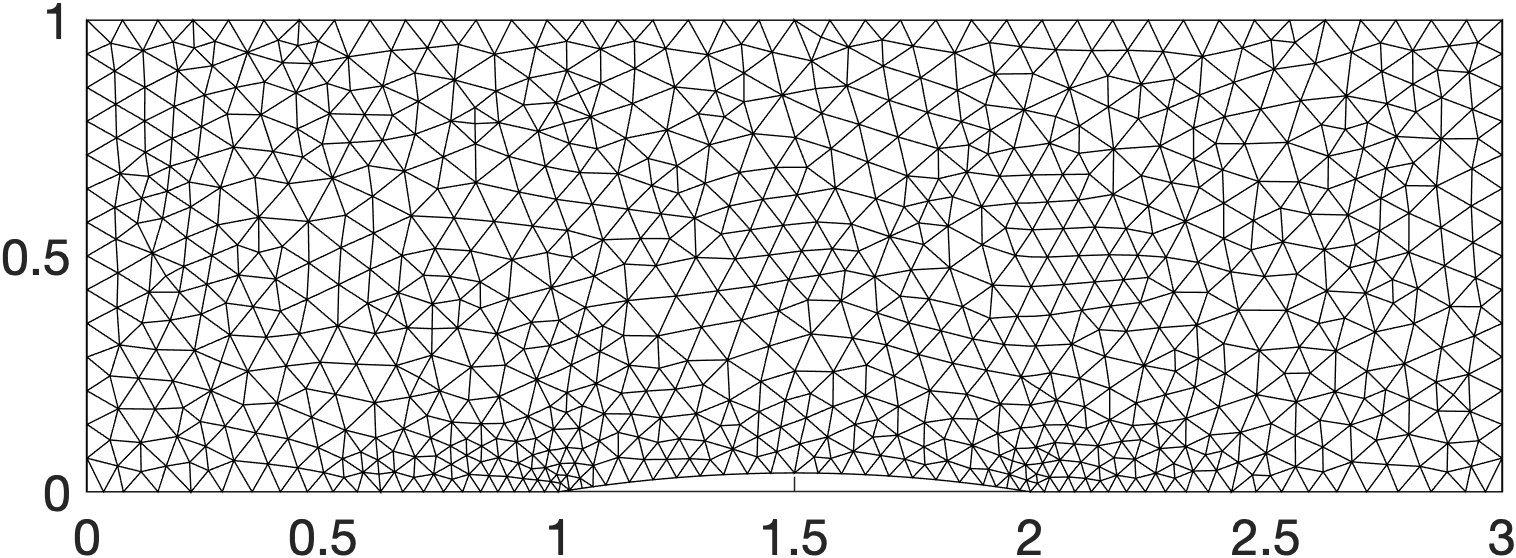}
  \caption{Initial mesh $\mathcal{T}_h$}    
	\end{subfigure}
	\hfill
	\begin{subfigure}[b]{0.49\textwidth}
		\centering
		\includegraphics[width=\textwidth]{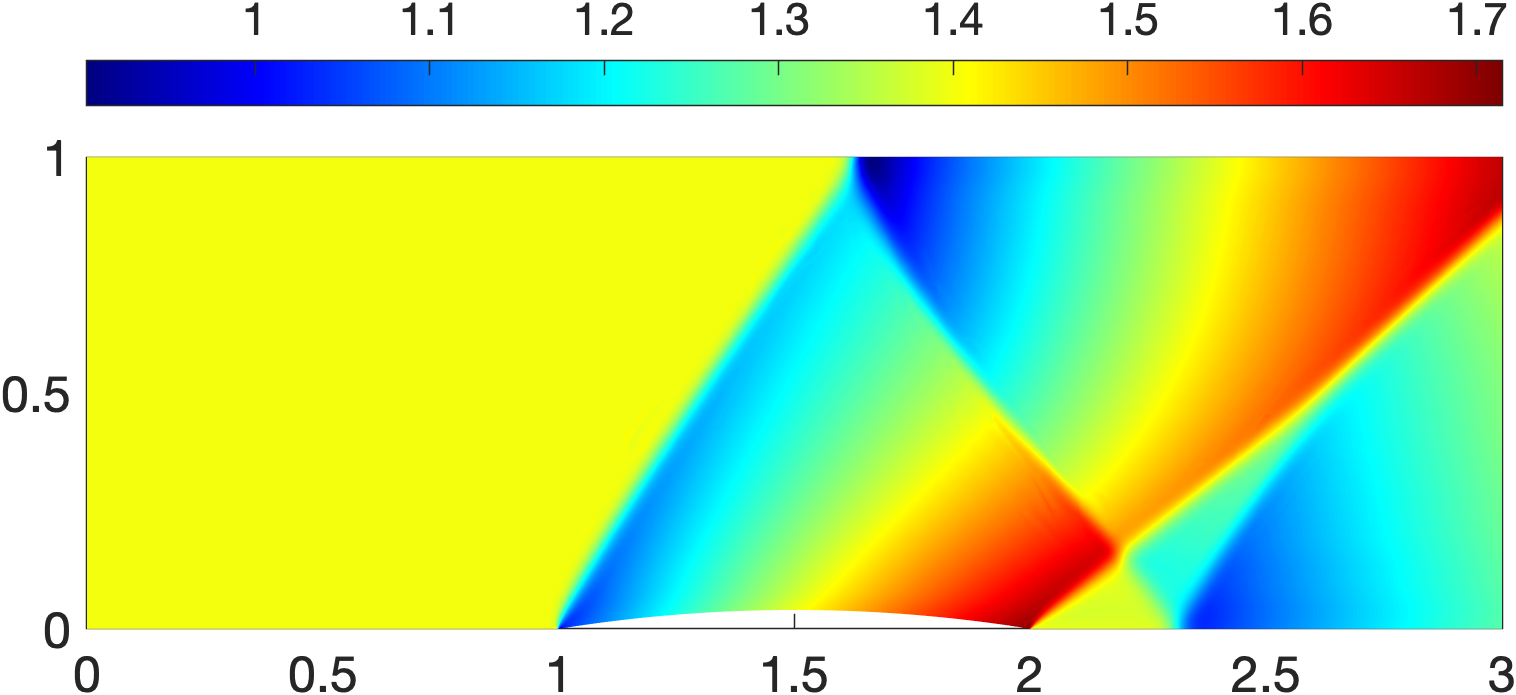}
  \caption{Step 1: Solution $\bm u_h$ on $\mathcal{T}_h$}
	\end{subfigure} 
        \\[1ex]
	\begin{subfigure}[b]{0.49\textwidth}
		\centering
		\includegraphics[width=\textwidth]{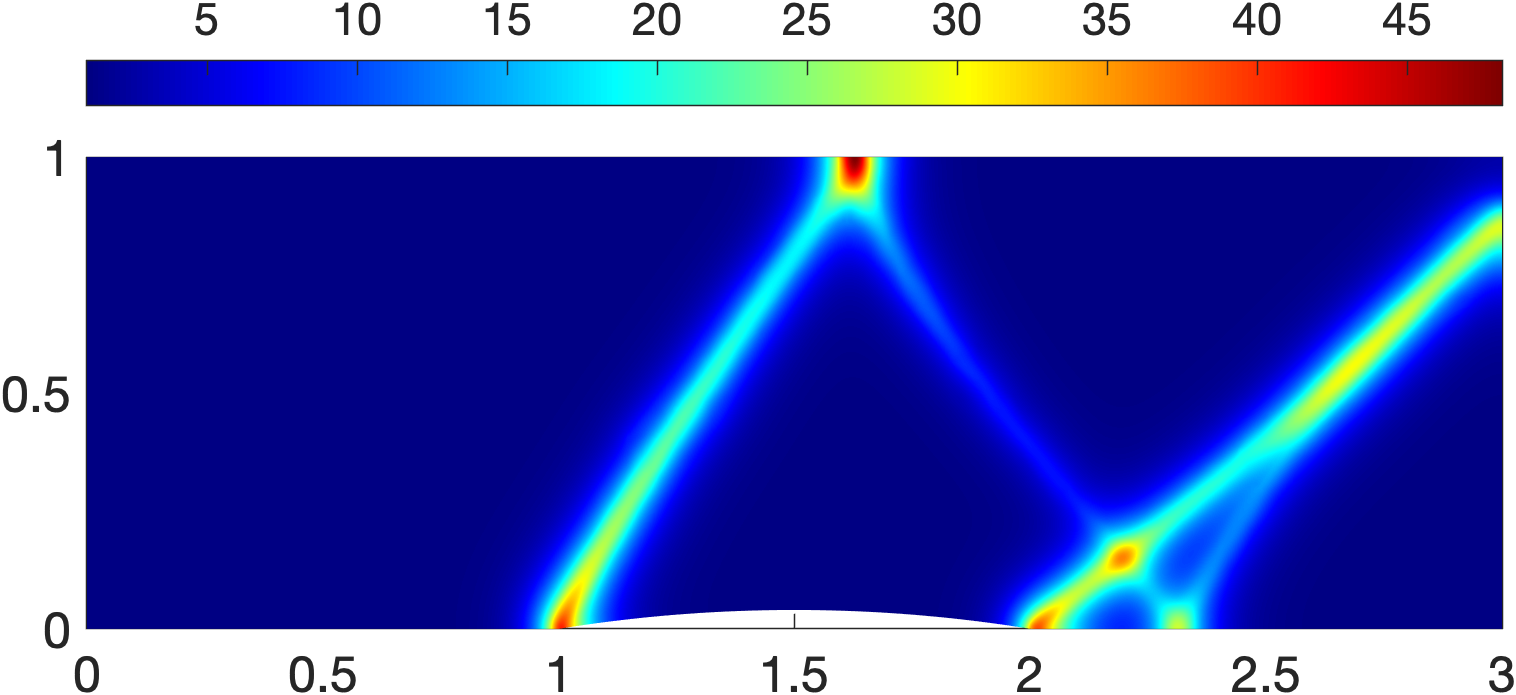}
  \caption{Step 2: Mesh density $\varrho'_h$ on $\mathcal{T}_h$}
	\end{subfigure} 
 	\begin{subfigure}[b]{0.49\textwidth}
		\centering		
  \includegraphics[width=\textwidth]{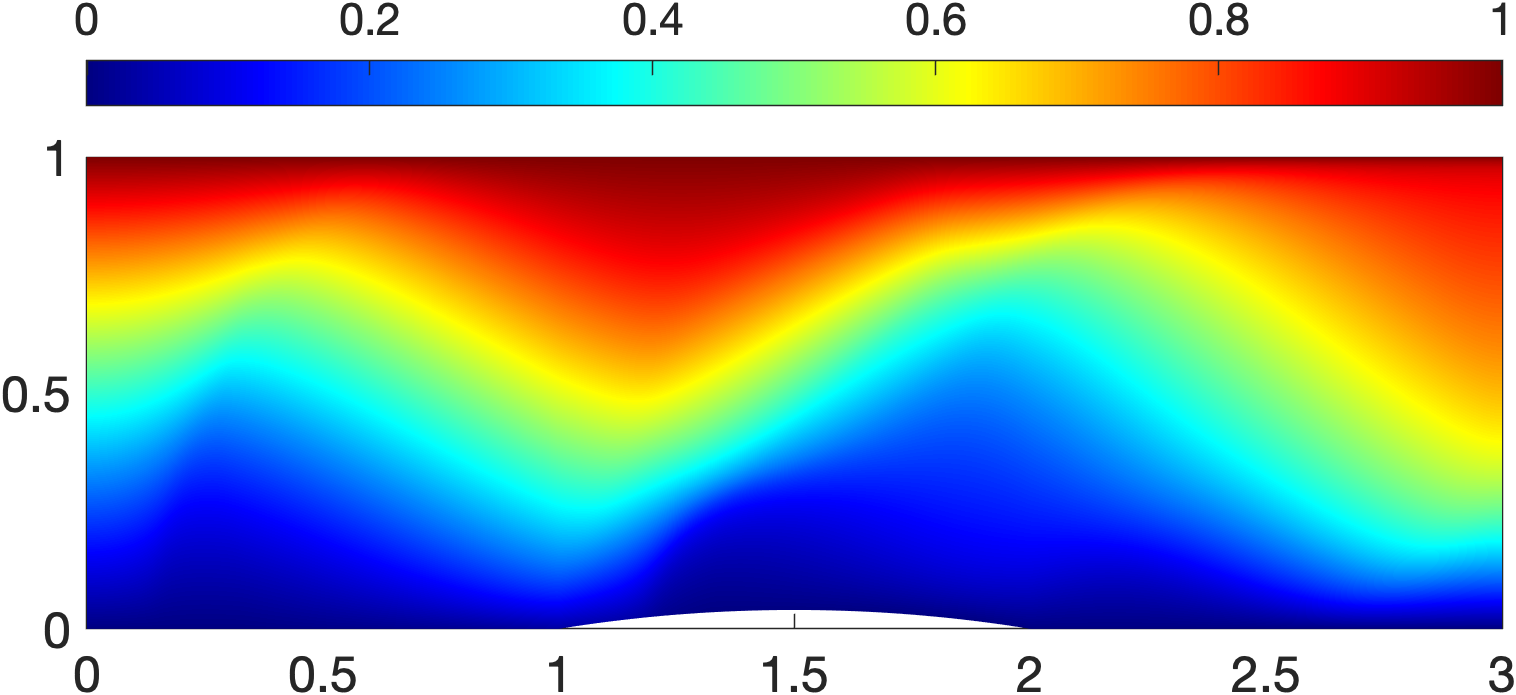}
  \caption{Step 3: Monge-Amp\'ere solution on $\mathcal{T}_h$}    
	\end{subfigure}
	\\[1ex]
	\begin{subfigure}[b]{0.49\textwidth}
		\centering		\includegraphics[width=\textwidth]{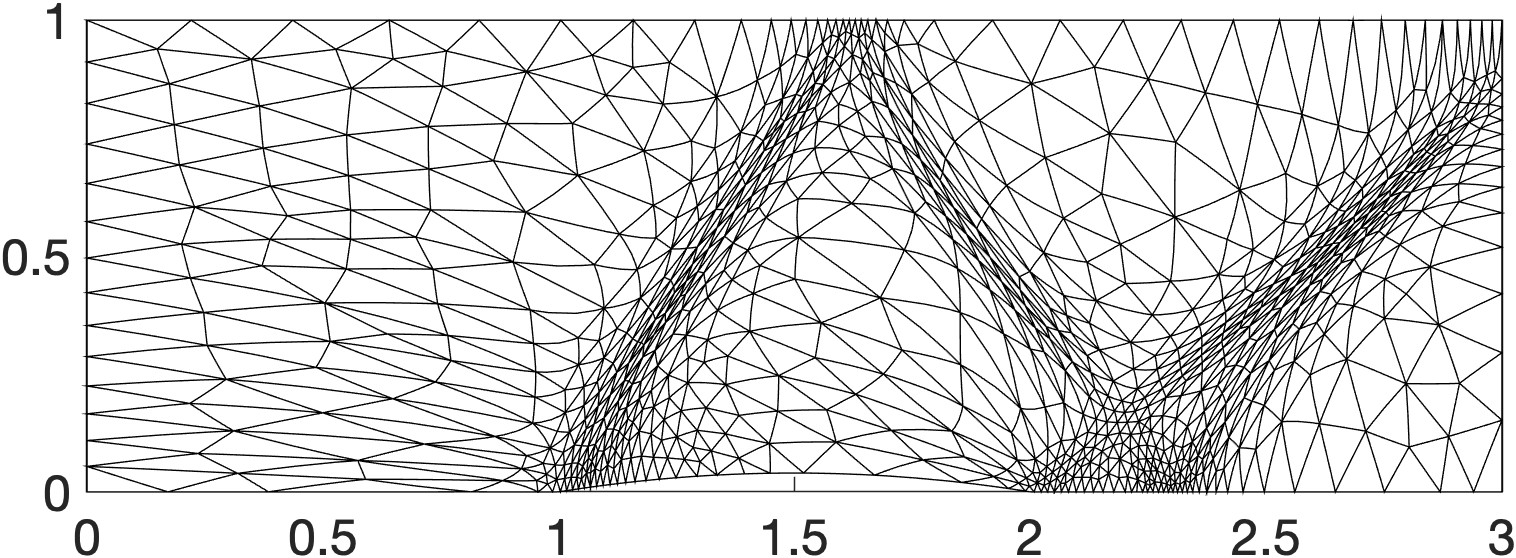}
        \caption{Step 4: Adaptive mesh $\mathcal{T}_h^*$}    
	\end{subfigure} 
        \hfill
	\begin{subfigure}[b]{0.49\textwidth}
		\centering
		\includegraphics[width=\textwidth]{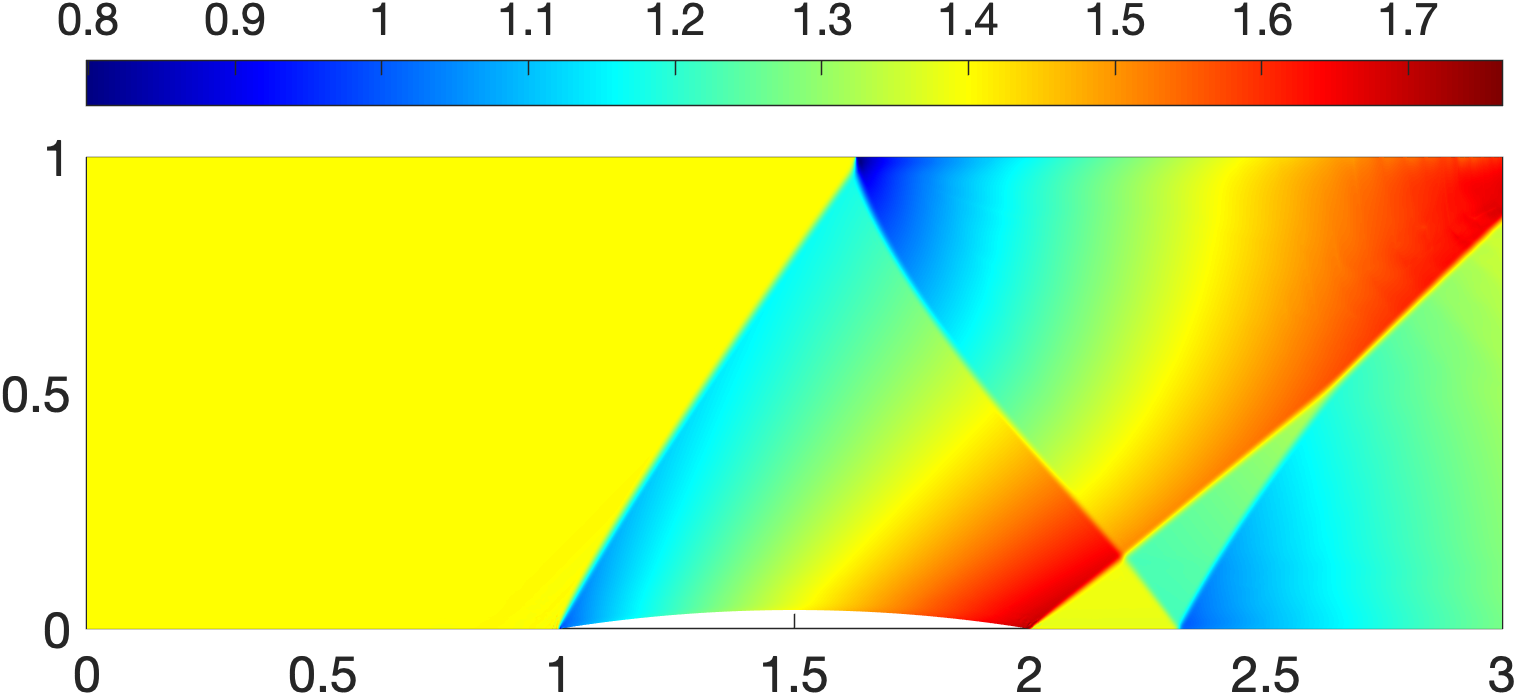}
        \caption{Step 5-6: Solution on $\bm u_h^*$ on $\mathcal{T}_h^*$}    
	\end{subfigure} 
\caption {Illustration of Algorithm 1 applied to the inviscid supersonic flow in a channel with a 4\% thick circular bump.}
\label{fig1}
\end{figure}

\section{Numerical Results} \label{sec:results}

In this section, we present numerical results for a number of well-known test cases to demonstrate the proposed approach. Unless otherwise specified, polynomial degree $k = 4$ is used to represent both the numerical solution and the geometry. Although the polynomial degree $k=4$ is relatively high for shock flows, our approach can compute the numerical solution without using the  solutions computed with lower polynomial degrees. 

\subsection{Inviscid transonic flow past NACA 0012 airfoil}

\begin{figure}[h]
	\centering
	\begin{subfigure}[b]{0.49\textwidth}
		\centering		\includegraphics[width=\textwidth]{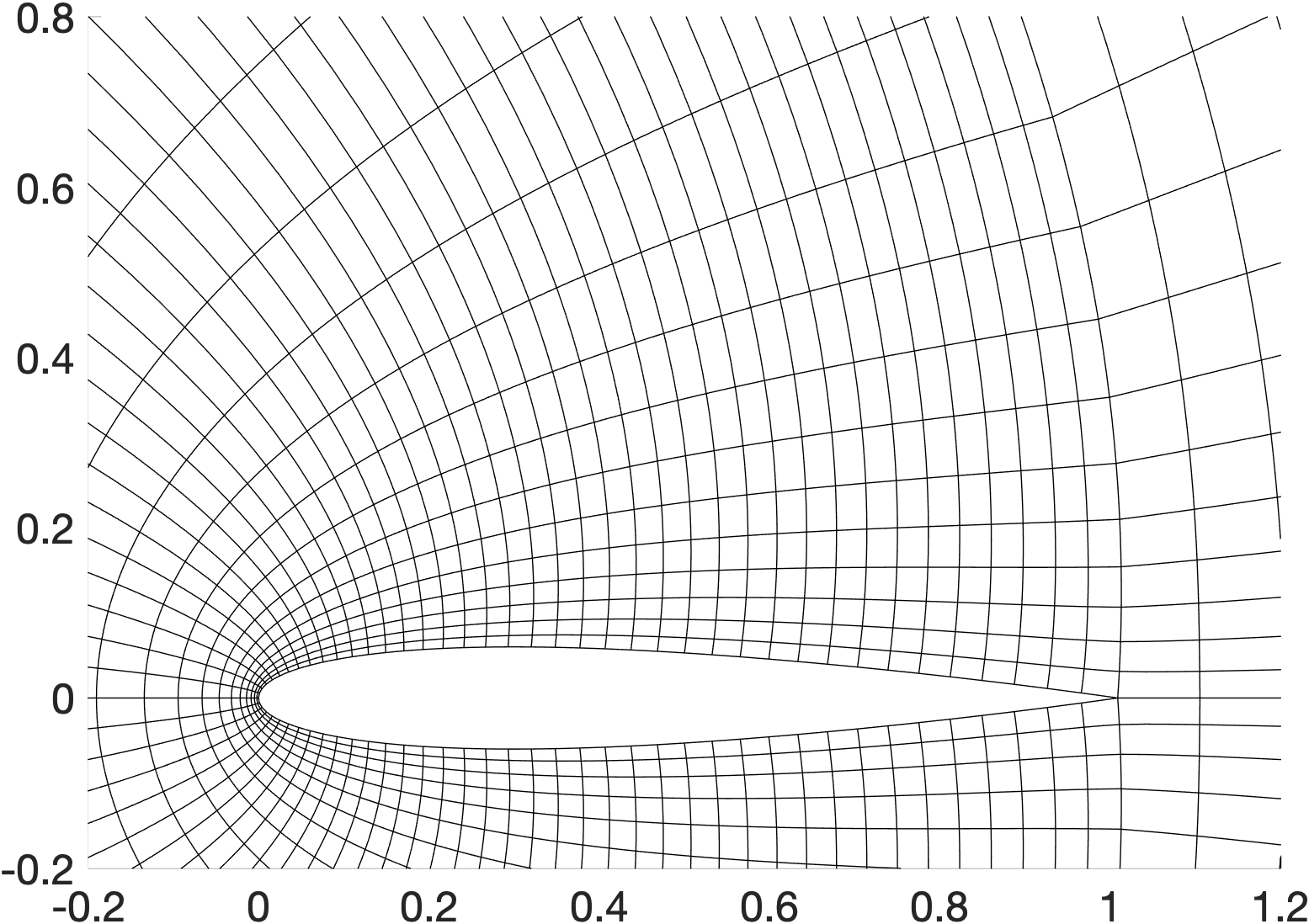}
  \caption{Initial mesh}
	\end{subfigure}
	\hfill
	\begin{subfigure}[b]{0.49\textwidth}
		\centering		\includegraphics[width=\textwidth]{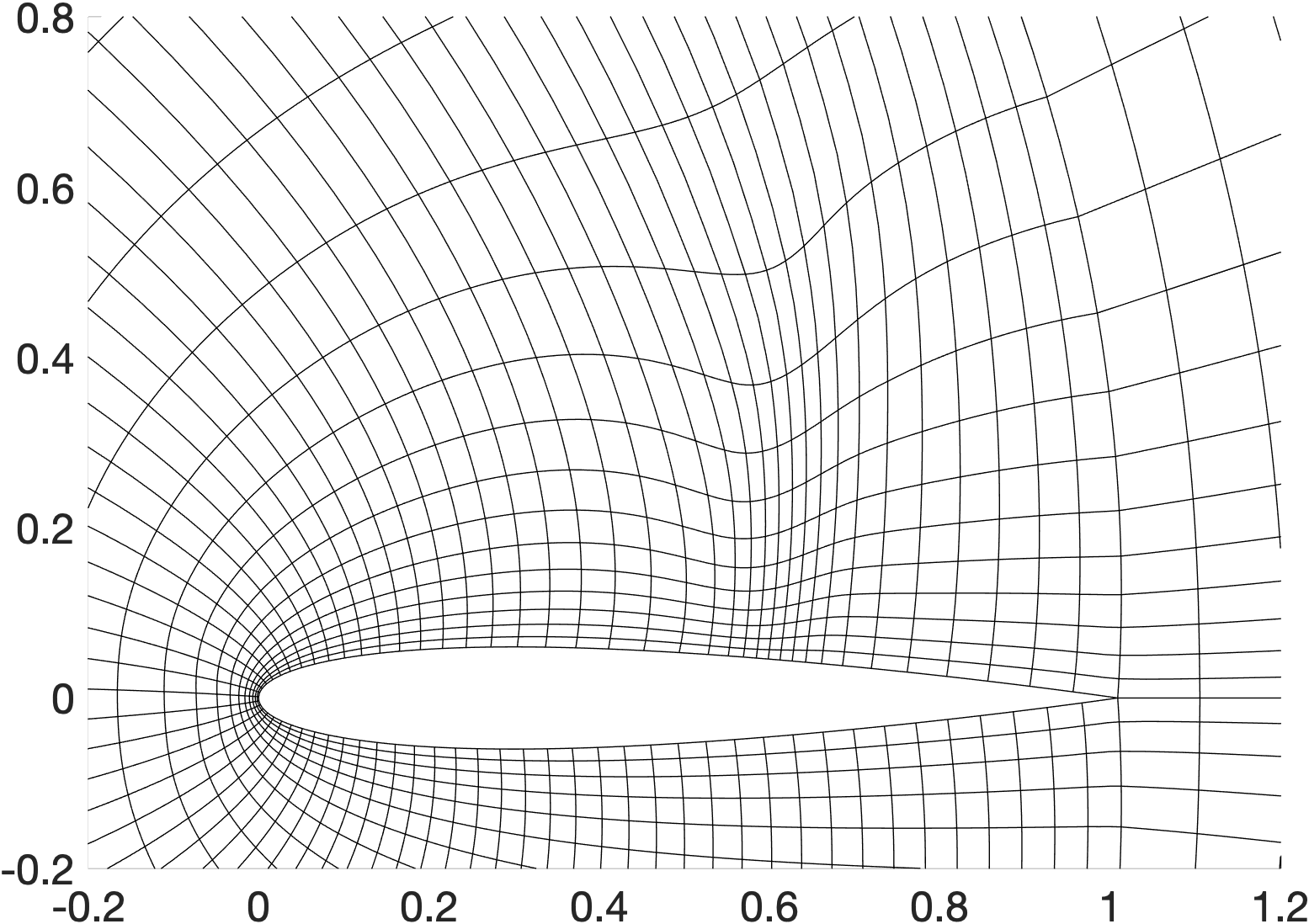}
  \caption{First adaptive mesh}
	\end{subfigure} 
 	\begin{subfigure}[b]{0.49\textwidth}
		\centering		\includegraphics[width=\textwidth]{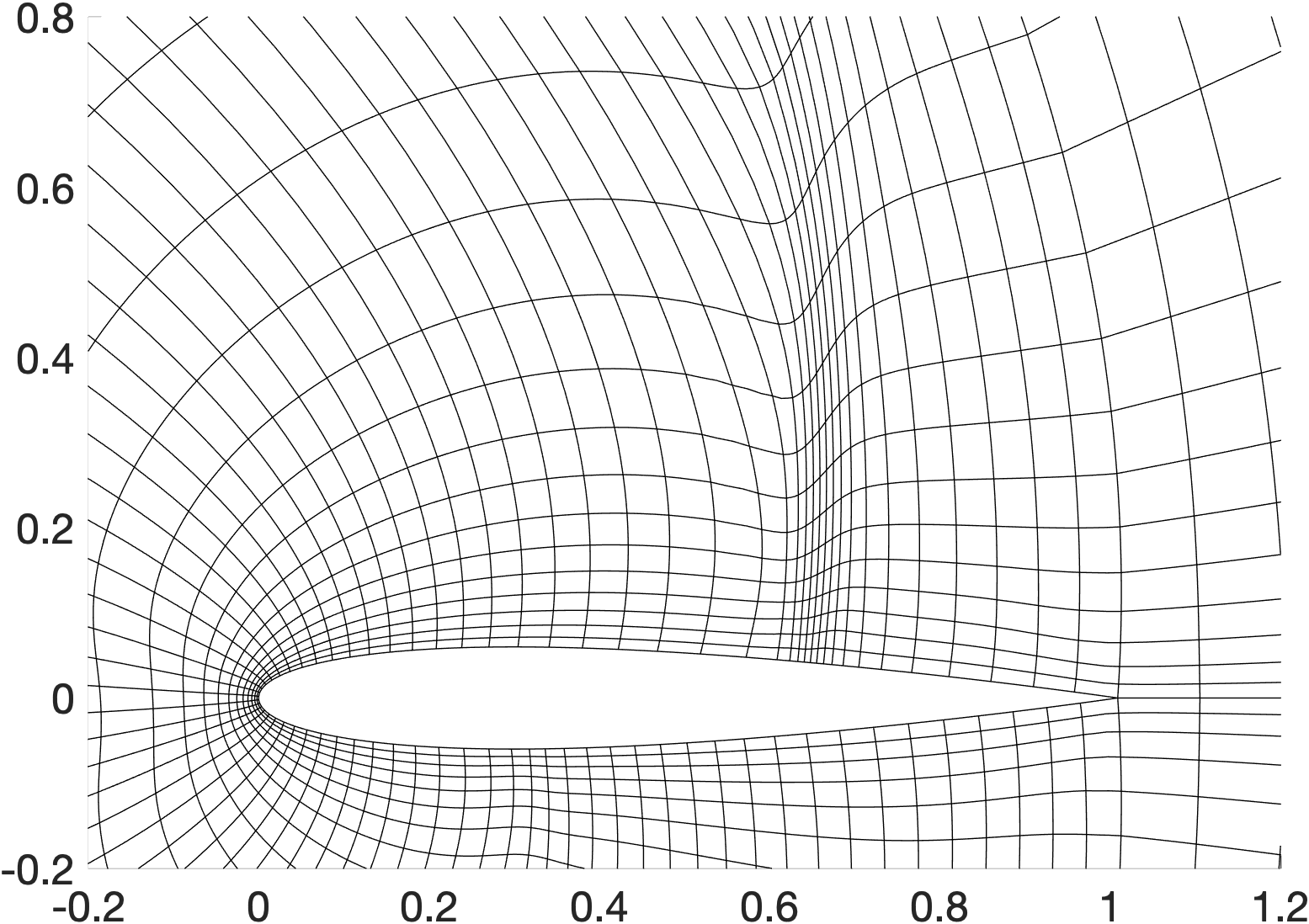}
		\caption{Second adaptive mesh}
	\end{subfigure}
	\hfill
	\begin{subfigure}[b]{0.49\textwidth}
		\centering		\includegraphics[width=\textwidth]{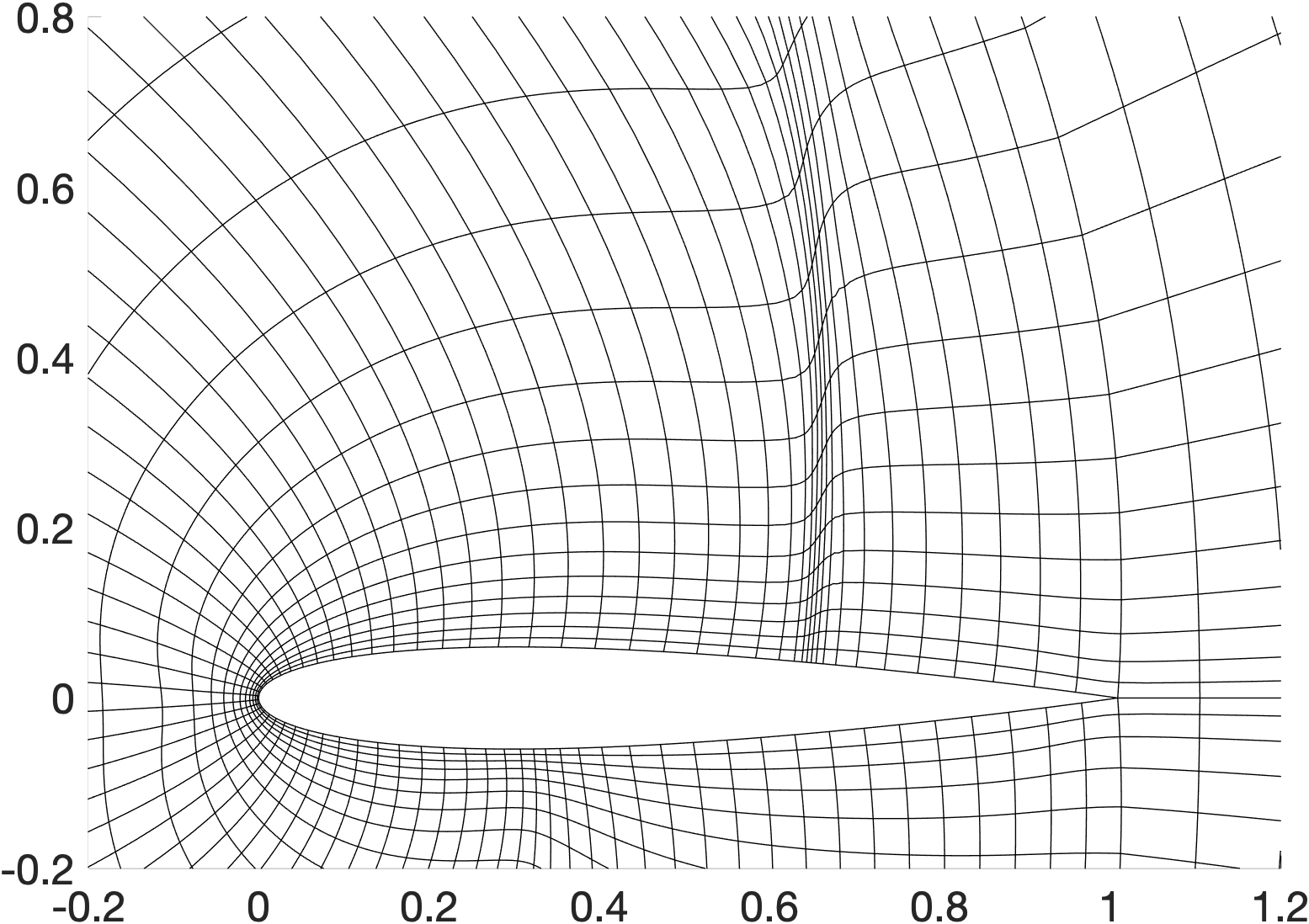}
		\caption{Third adaptive mesh}
	\end{subfigure} 
	\caption{Close-up view near the NACA0012 airfoil of the sequence of adaptive meshes for inviscid transonic flow at angle of attack $\alpha = 1.5^{\rm o}$ and freestream Mach number $M_\infty = 0.8$.}
 \label{fignaca2}
\end{figure}

\begin{figure}[htb]
	\centering
	\begin{subfigure}[b]{0.49\textwidth}
		\centering		\includegraphics[width=\textwidth]{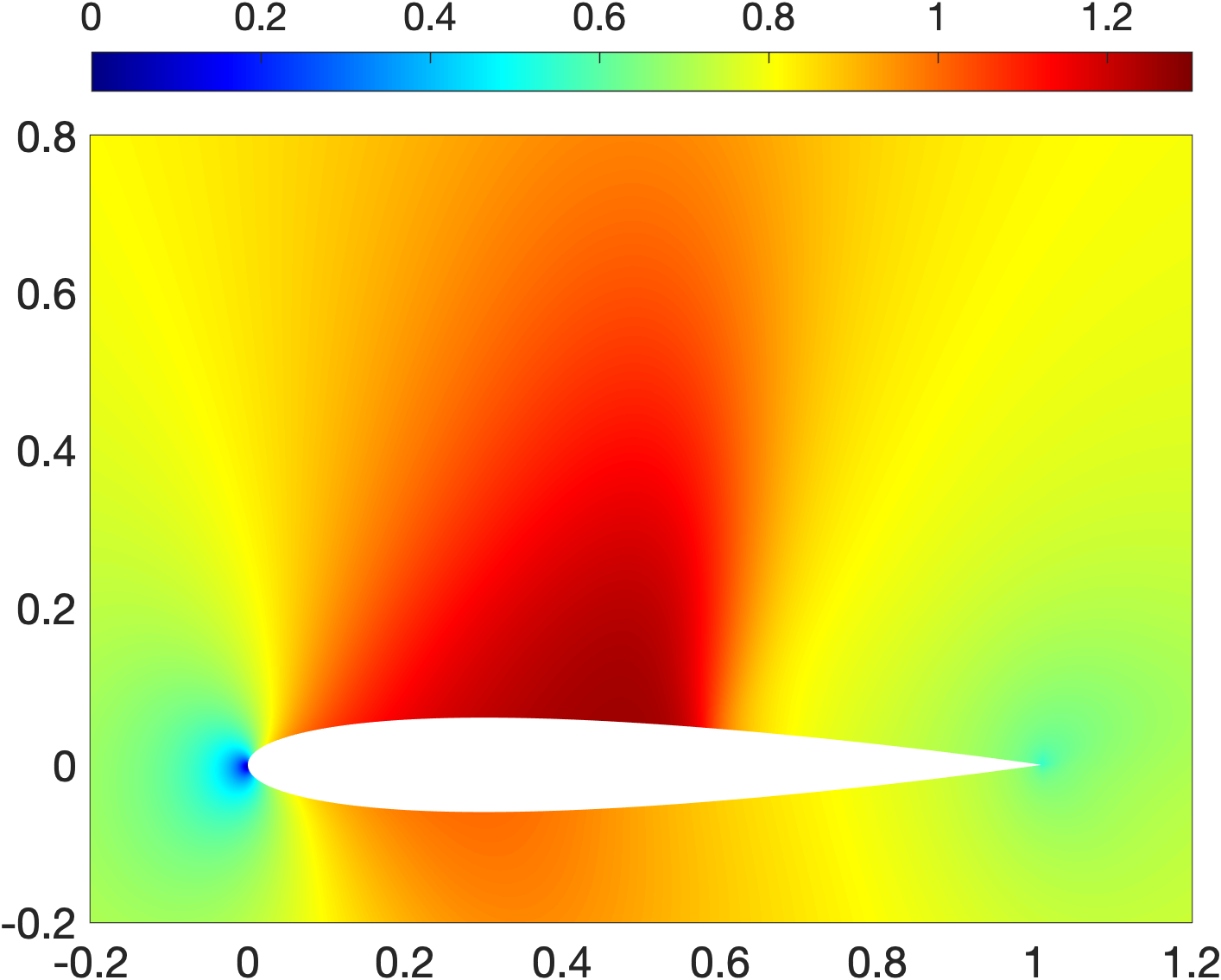}
  \caption{Initial mesh}
	\end{subfigure}
	\hfill
	\begin{subfigure}[b]{0.49\textwidth}
		\centering		\includegraphics[width=\textwidth]{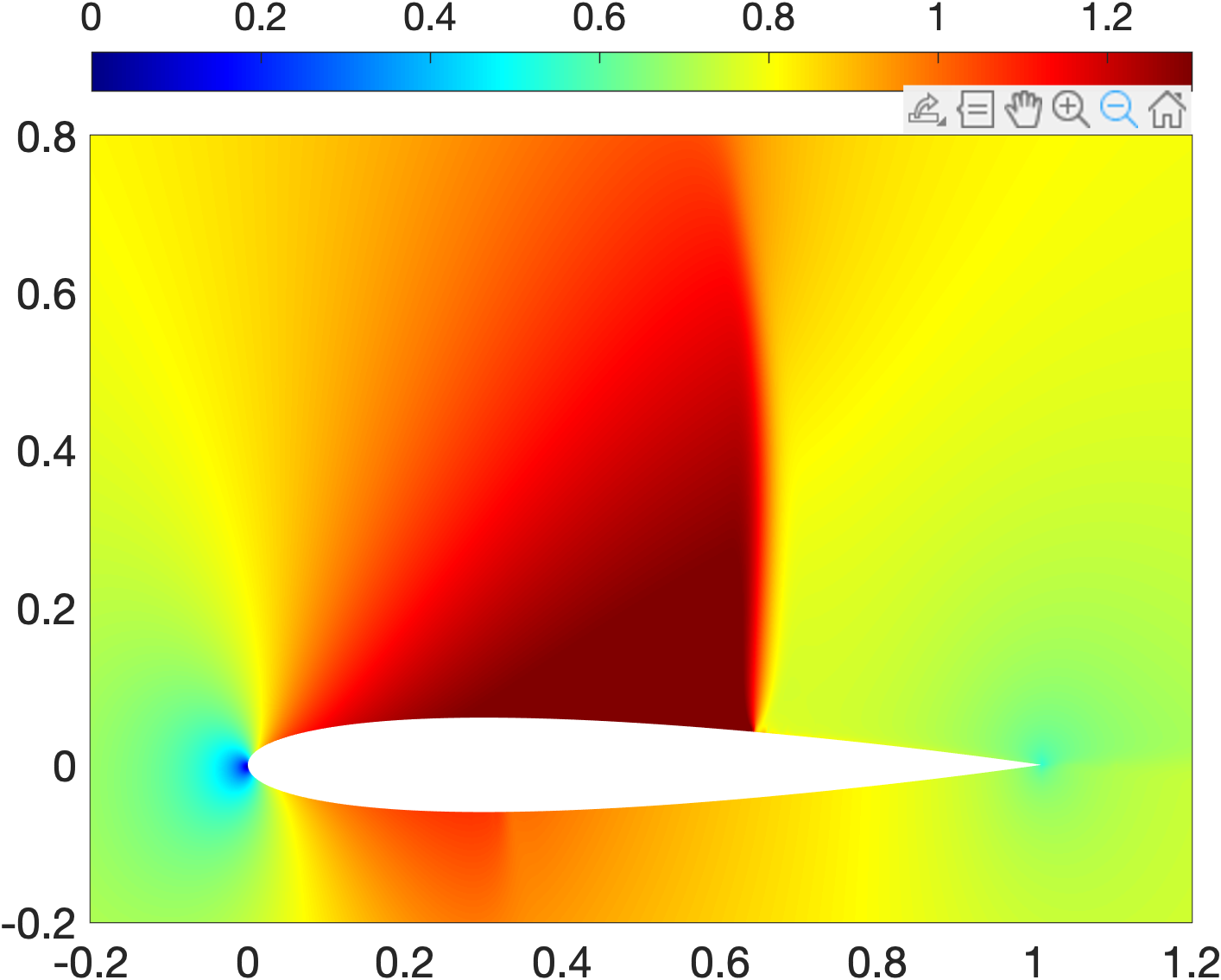}
  \caption{First adaptive mesh}
	\end{subfigure} 
 	\begin{subfigure}[b]{0.49\textwidth}
		\centering		\includegraphics[width=\textwidth]{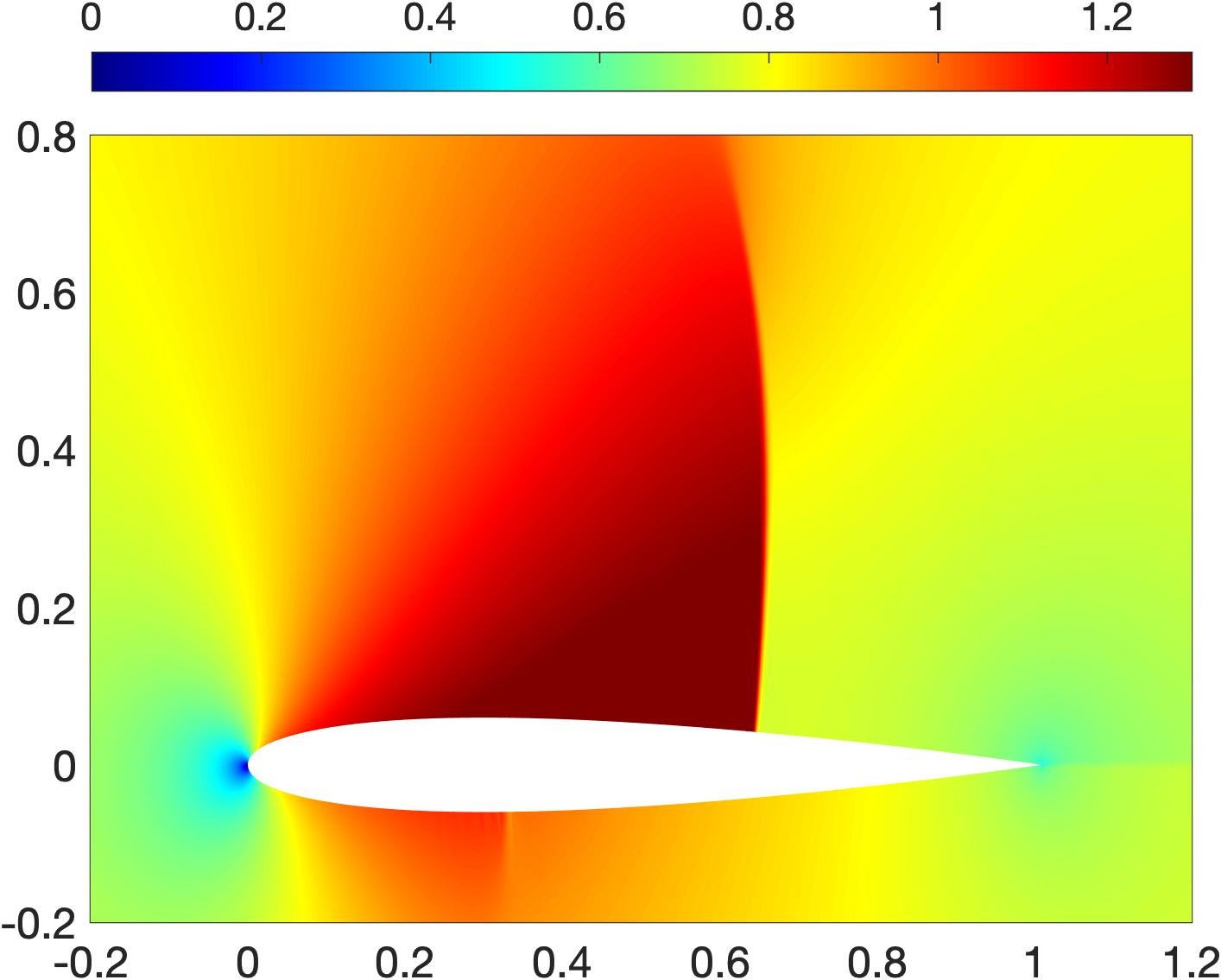}
		\caption{Second adaptive mesh}
	\end{subfigure}
	\hfill
	\begin{subfigure}[b]{0.49\textwidth}
		\centering		\includegraphics[width=\textwidth]{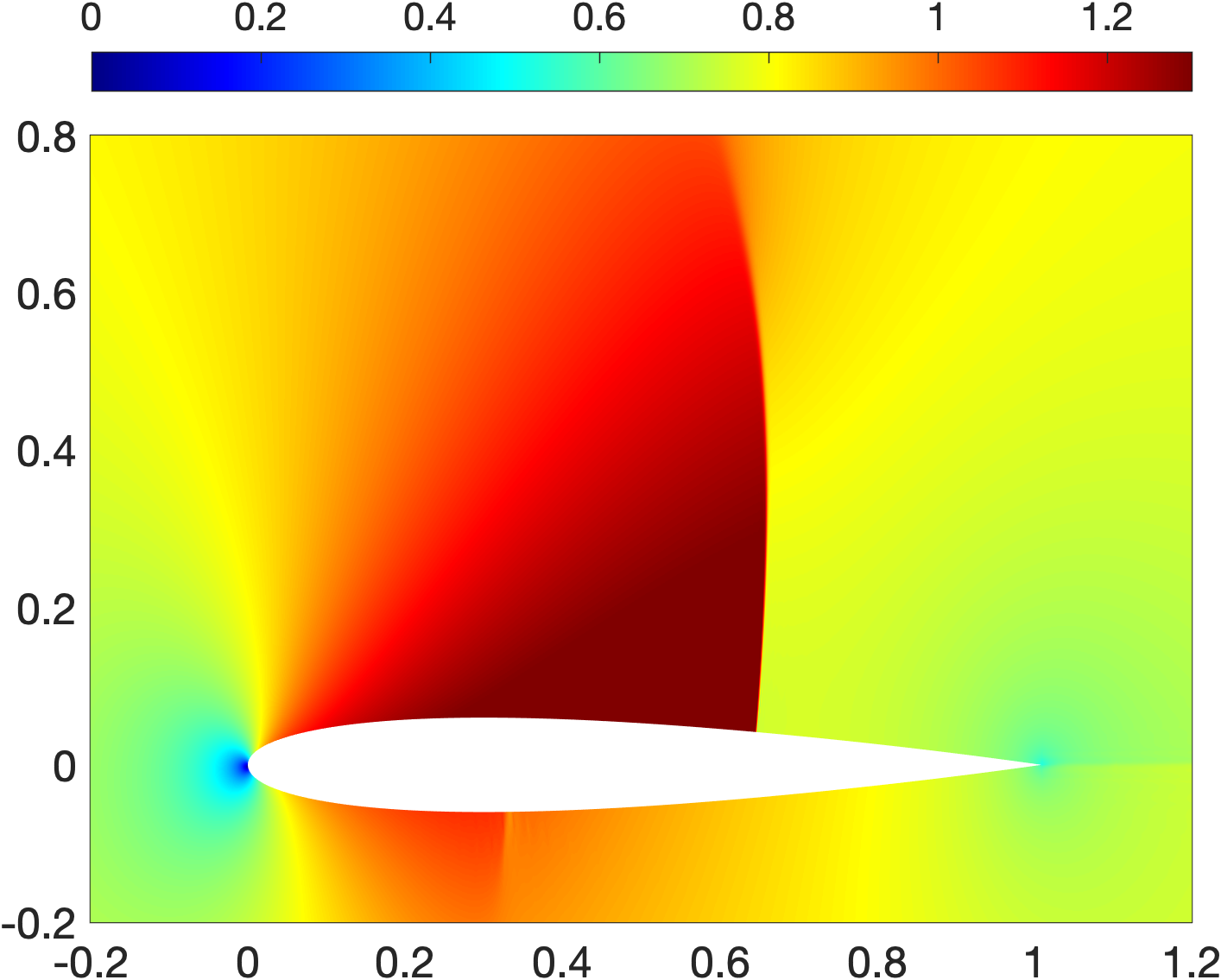}
		\caption{Third adaptive mesh}
	\end{subfigure} 
	\caption{Mach number computed on the initial  and adaptive meshes for the inviscid transonic flow past NACA 0012 airfoil.}
 \label{fignaca3}
\end{figure}

The first test case is an inviscid transonic flow past a NACA 0012 airfoil at angle of attack $\alpha = 1.5^{\rm o}$ and freestream Mach number $M_\infty = 0.8$ \cite{Nguyen2011a}. Slip velocity boundary condition is imposed on the airfoil, while far-field boundary condition is imposed on the rest of the boundary. A  shock is formed on the upper surface, while another  weaker shock is formed under the lower surface. Figure~\ref{fignaca2} depicts the initial mesh and three consecutive adaptive meshes near the airfoil surface. Figure \ref{fignaca3} shows the Mach number computed on the initial mesh and the r-adaptive meshes. 



\begin{figure}[htbp]
	\centering
 \includegraphics[width=0.65\textwidth]{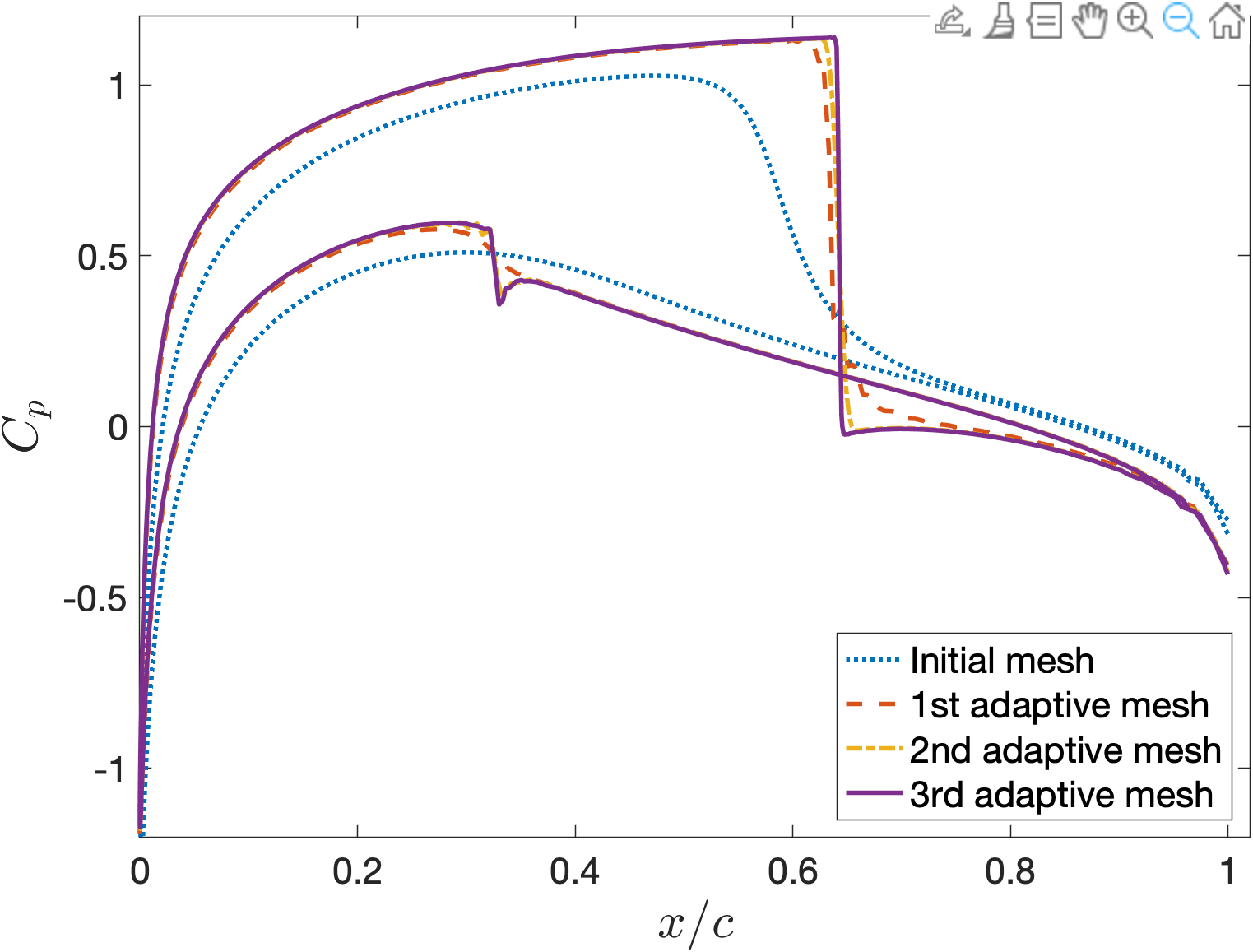}
 \caption{Profiles of the pressure coefficient computed on the initial and adaptive meshes for the inviscid transonic flow past NACA 0012 airfoil.}
	\label{fignaca4}
\end{figure}

It is interesting to see how the elements of the initial mesh are moved to create new meshes that align well with the shocks. The results also show how the numerical solution is improved and how the shocks are better resolved over each iteration of the mesh adaptation procedure. We observe that the shocks are well resolved on the final adaptive mesh and that the solution on the final mesh is accurate. This can be clearly seen from the profiles of the computed pressure coefficient shown in Figure \ref{fignaca4}. We see that the pressure coefficient profiles converge rapidly and that the profile computed on the second adaptive mesh is very similar to that computed on the third adaptive mesh. We emphasize that the profile on the third adaptive mesh is very sharp at the shocks, yet there is no oscillation and overshoot.

\subsection{Inviscid supersonic flow over a double ramp}

This test case is used in \cite{carnes2019code} as a building block towards more complicated double wedge and cone flows. 
The geometry is a double-ramp with a $25^{\circ}$  incline for the first ramp and $37^{\circ}$  incline for the second. 
Note that the second angle is shallower than typical hypersonic double wedge or cone flows \cite{olejniczak1997numerical}. 
We consider supersonic flow at a free-stream Mach number of 3.6, for which the resulting flow-field is relatively simple. Two shocks are expected to emanate from the corners and intersect to form a third shock. 

\begin{figure}[htbp]
\centering
        \begin{subfigure}[b]{0.3\textwidth}
		\centering		\includegraphics[width=0.8\textwidth]{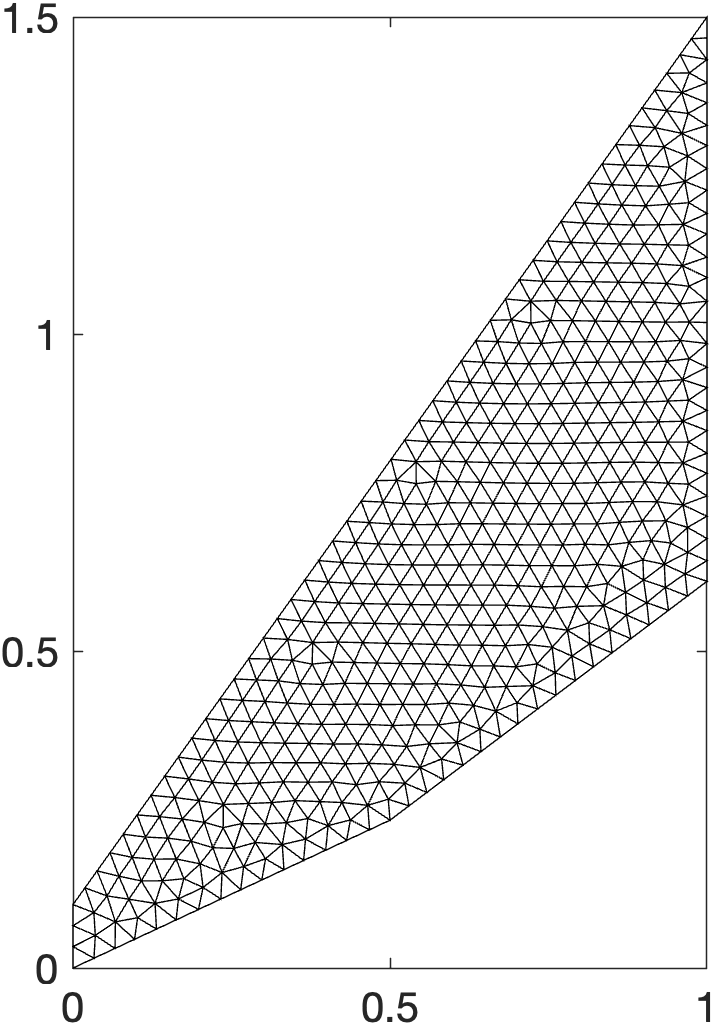}
  \caption{Initial mesh}
	\end{subfigure}
	\begin{subfigure}[b]{0.34\textwidth}
		\centering		\includegraphics[width=0.8\textwidth]{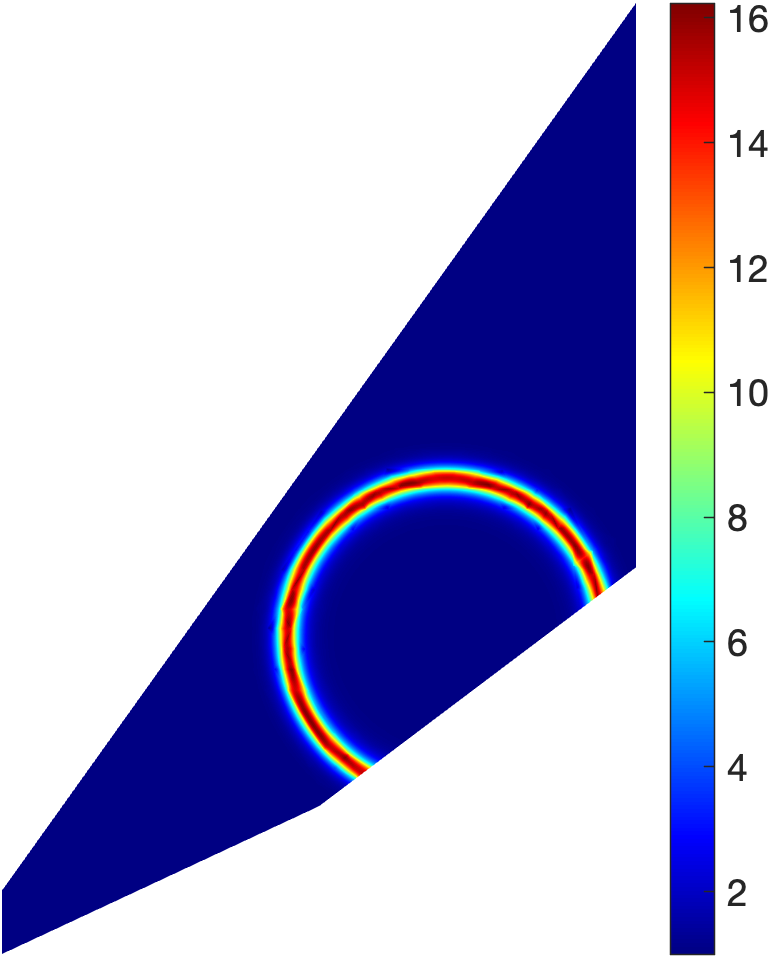}
  \caption{Artificial sensor $\varrho'$}
	\end{subfigure}
        \begin{subfigure}[b]{0.3\textwidth}
		\centering		\includegraphics[width=0.8\textwidth]{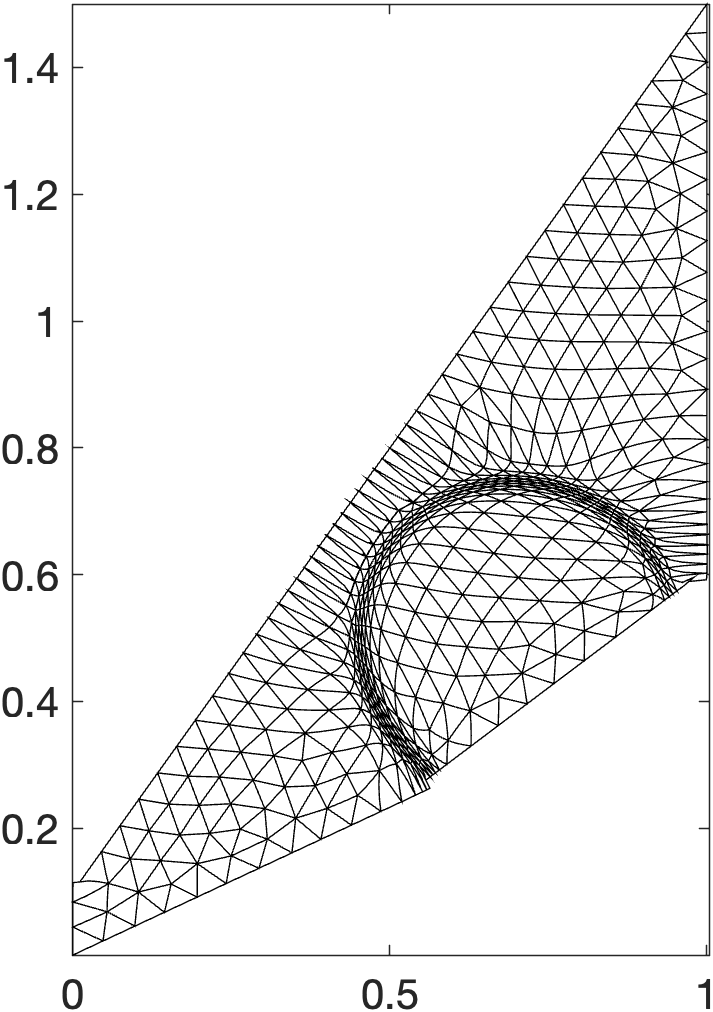}
  \caption{Adapted mesh}
	\end{subfigure}
\caption{A demonstration of the corner separation that that can occur if boundary nodes are not allowed to move between boundaries. }
\label{fig:cornersep}	
\end{figure}

\begin{figure}[htbp!]
\centering
        \begin{subfigure}[b]{0.37\textwidth}
		\centering		\includegraphics[width=0.8\textwidth]{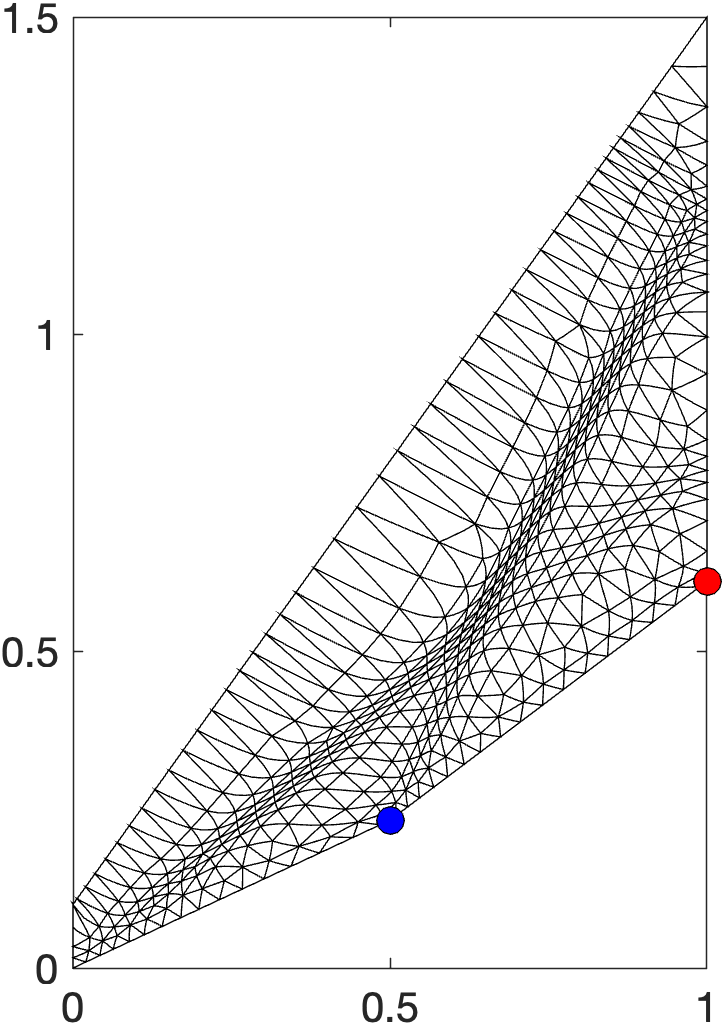}
  \caption{Monge-Amp\`ere adaptive mesh on double ramp geometry.}
	\end{subfigure}
	\begin{subfigure}[b]{0.62\textwidth}
		\centering		\includegraphics[width=0.8\textwidth]{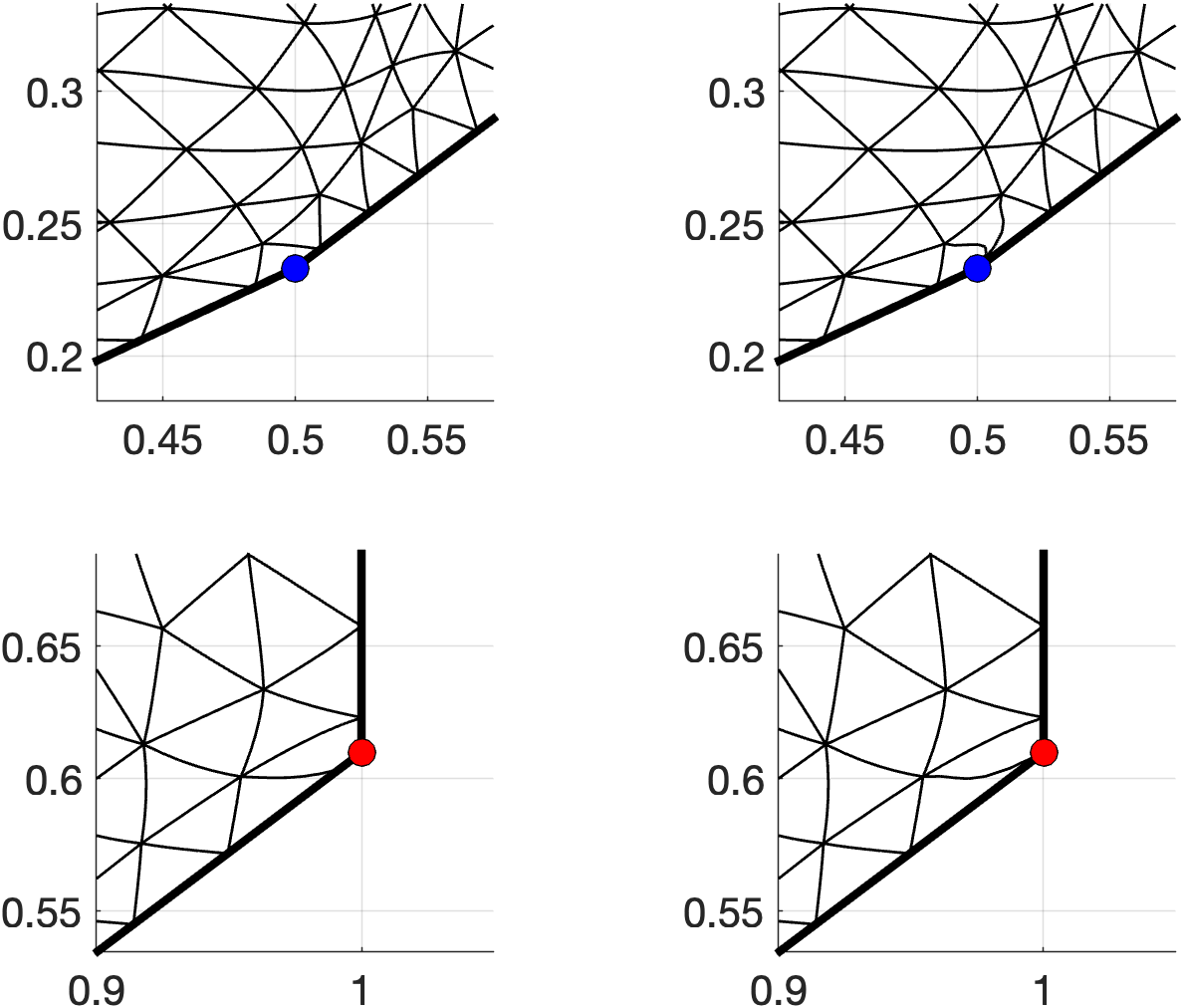}
  \caption{Before (left) and after (right) corner fix is applied}
	\end{subfigure}
\caption{Adapted mesh (left) conforms to boundary but requires a manual reassignment of nodes to the geometry corners (right). }
\label{fig:corners}	
\end{figure}

The purpose of this example is to examine the ability of the Monge-Amp\`ere solver to refine the mesh on polygonal domains with flow over corners of the domain. 
The boundary consists of six line segments $\{\bm \Gamma_i\}_{i=1}^6$ defined as $c_i(\bm x) := \bm A_i \bm x + \bm b_i = 0$. 
Enforcing that each $\bm q_h$ on $\Omega_i$ must satisfy $c_i(\bm q_h) = 0$ led to meshes that would detach at the corners, hampering convergence. 
This is demonstrated in Figure \ref{fig:cornersep} with an artificial target density. 
Whether this phenonema is a result of the HDG discretization or the formulation of Monge-Ampe\`re on this domain remains to be determined. 
This issue is  addressed by changing the Neumann boundary condition to obey a global description of the geometry; instead of enforcing that $\bm q_h$ at $\bm \Gamma_i$ must satisfy $c_i(\bm q_h) = 0$, it is allowed to transition onto adjacent faces if the value of $\bm q_h$ leaves the bounds of $\bm \Gamma_i$. 
In this way, boundary nodes are allowed to slide along the boundary and move from one face to another. 
Other $r$-adaptive methods have found it advantageous to fix nodes at boundaries rather than let them transition from one boundary to another. 
Since the domain mapping is determined as the gradient of a scalar potential, we cannot explicitly fix the location of certain nodes. 
Instead, after the adaptive mesh is formed, the element that crosses a corner is identified and its closest vertex is moved to that same corner, in order to not change the definition of the geometry. 
This procedure is illustrated in Figure \ref{fig:corners}.

\begin{figure}[htbp]
\centering
        \begin{subfigure}[b]{0.32\textwidth}
		\centering		\includegraphics[width=0.8\textwidth]{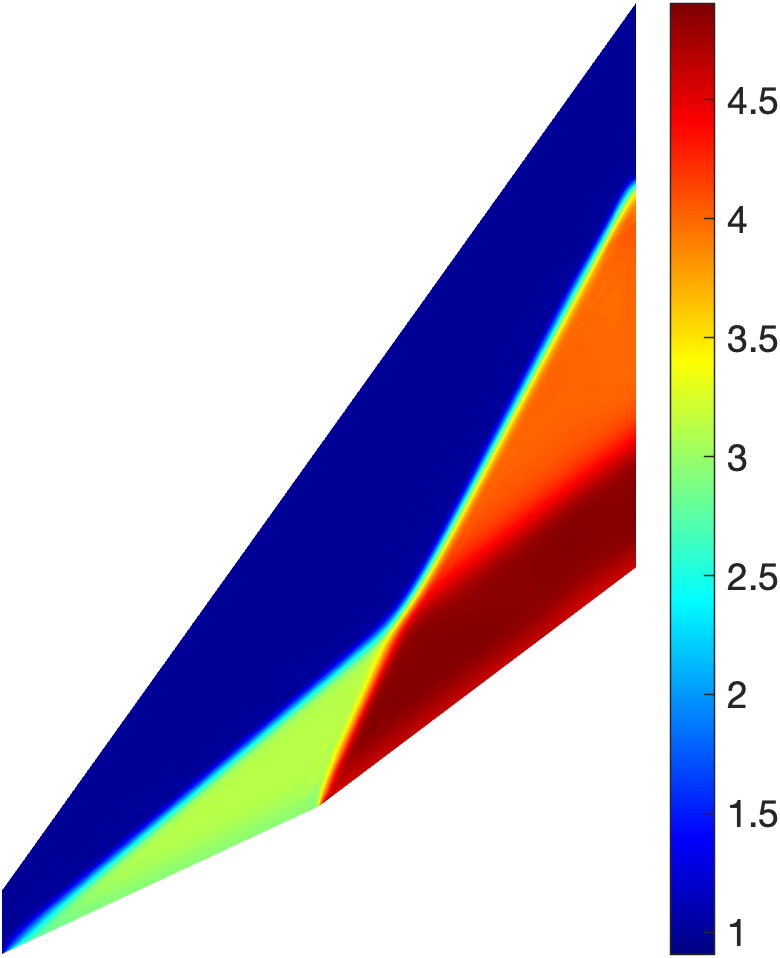}
  \caption{$\rho/\rho_{\infty}$}
	\end{subfigure}
	\begin{subfigure}[b]{0.32\textwidth}
		\centering		\includegraphics[width=0.8\textwidth]{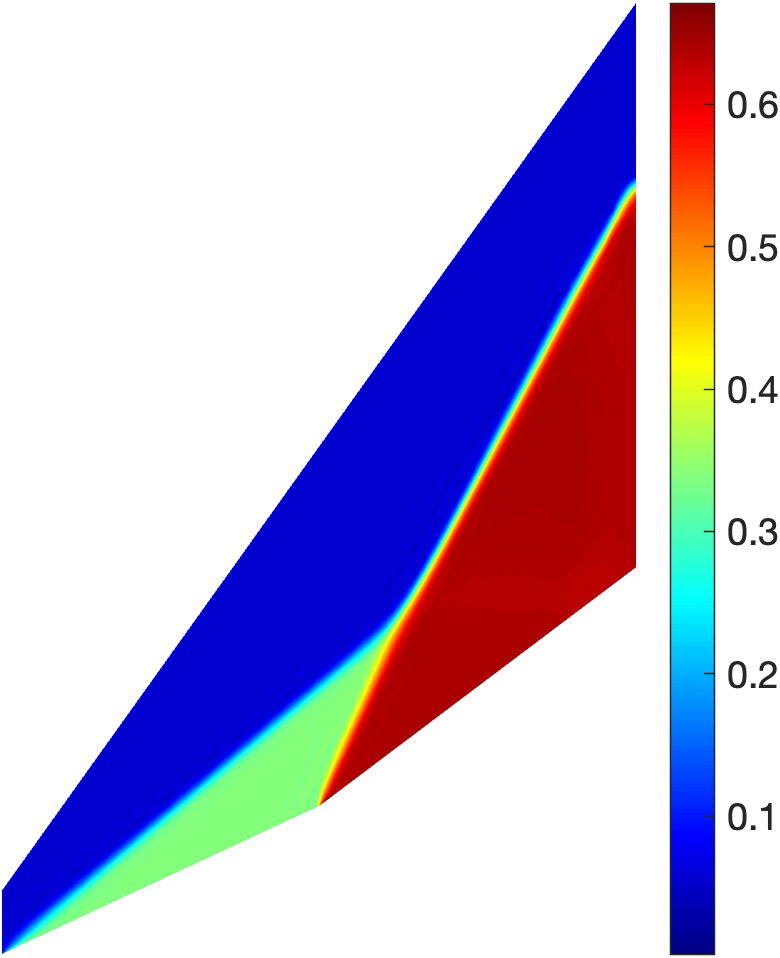}
  \caption{$p/p_{\infty}$}
	\end{subfigure}
        \begin{subfigure}[b]{0.32\textwidth}
		\centering		\includegraphics[width=0.8\textwidth]{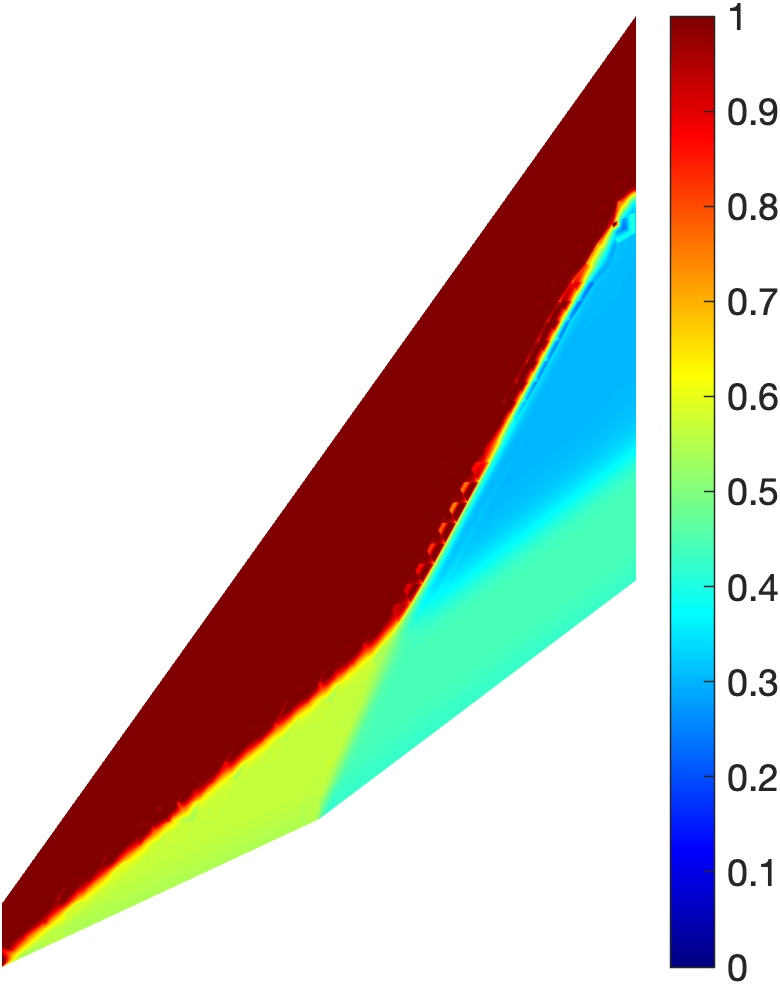}
  \caption{$M/M_{\infty}$}
	\end{subfigure}
\caption{Numerical solution computed on the initial mesh for supersonic inviscid flow over a double ramp.}
\label{fig:wedge}
\end{figure}

\begin{figure}[htbp]
\centering
        \begin{subfigure}[b]{0.32\textwidth}
		\centering		\includegraphics[width=0.8\textwidth]{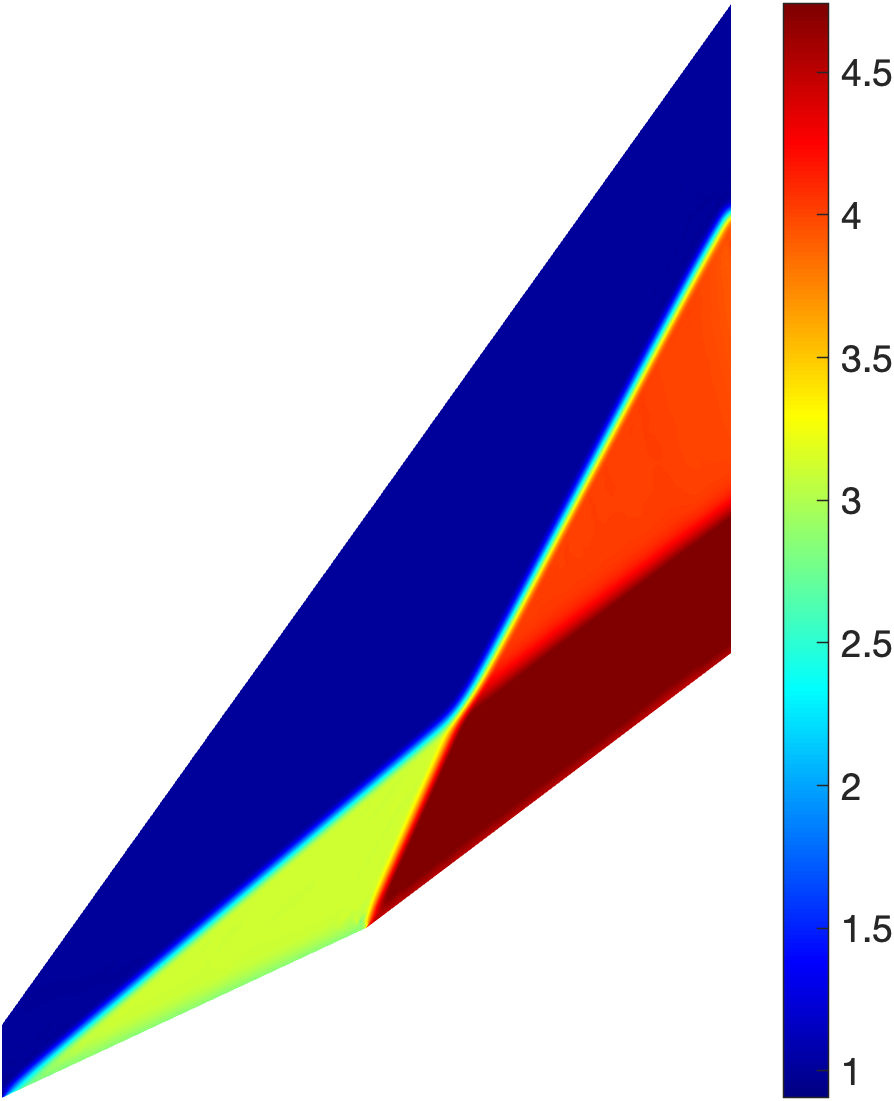}
  \caption{$\rho/\rho_{\infty}$}
	\end{subfigure}
	\begin{subfigure}[b]{0.32\textwidth}
		\centering		\includegraphics[width=0.8\textwidth]{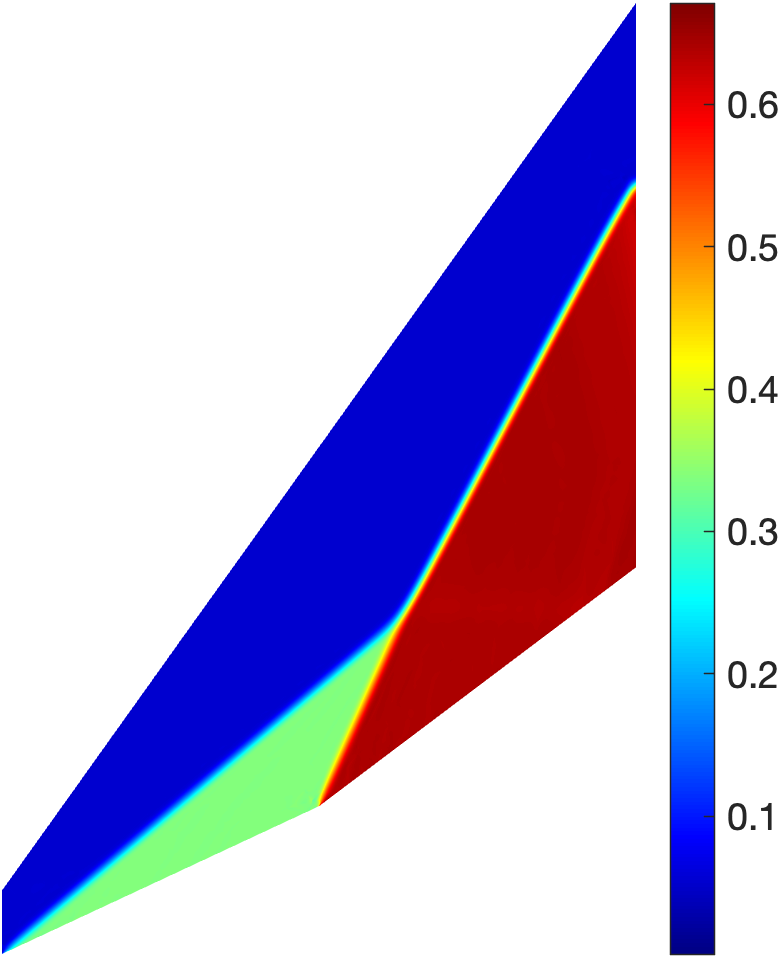}
  \caption{$p/p_{\infty}$}
	\end{subfigure}
        \begin{subfigure}[b]{0.32\textwidth}
		\centering		\includegraphics[width=0.8\textwidth]{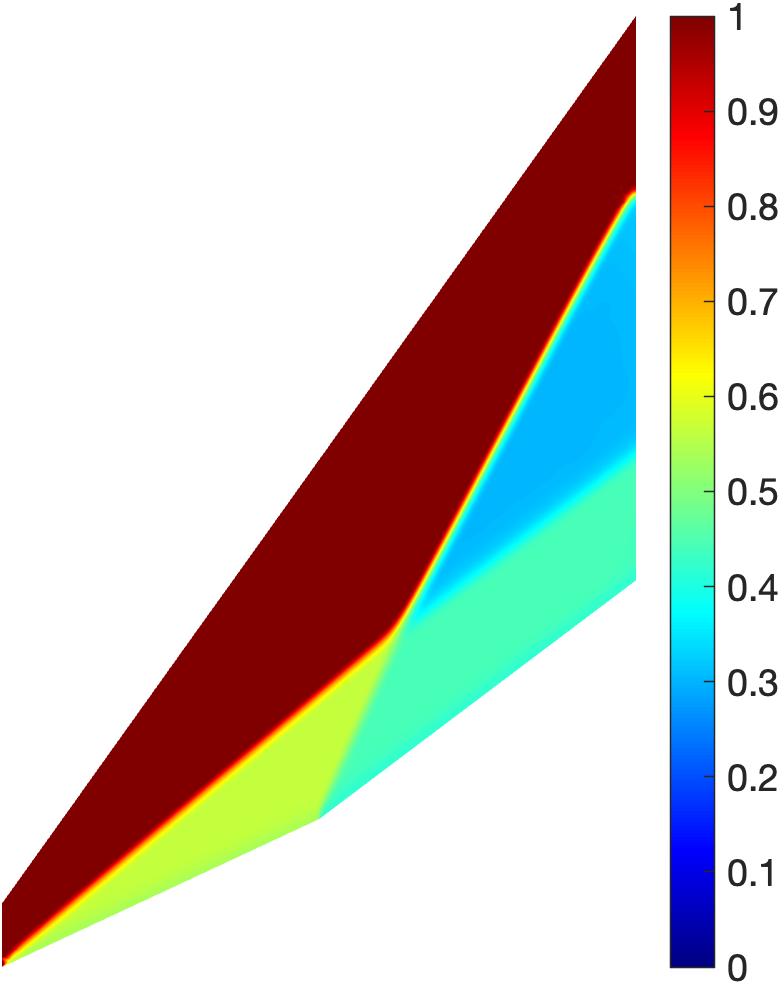}
  \caption{$M/M_{\infty}$}
	\end{subfigure}
\caption{Numerical solution computed on the adapted mesh for supersonic inviscid flow over a double ramp.}
\label{fig:wedgeadapt}	
\end{figure}

The starting grid consists of 909 elements and polynomial order $k=3$. The results on the initial mesh are shown in Figure \ref{fig:wedge}. 
We use the sensor based off the gradient of the physical density \eqref{mdf2} with $\beta = 1$ in order to get some refinement along the contact discontinuity, which would be missed with the sensor based off the divergence of the velocity. 
While the starting mesh is fine enough to capture the density and pressure well, visible oscillations are present in the Mach number field. These oscillations are not visible with mesh adaptation and the primary shocks and contact discontinuities are sharper than on the starting mesh.  See Figure \ref{fig:wedgeadapt}.

\subsection{Inviscid hypersonic flow past unit circular cylinder}

\begin{figure}[htbp]
\centering
	\begin{subfigure}[b]{0.3\textwidth}
		\centering		\includegraphics[width=0.6\textwidth]{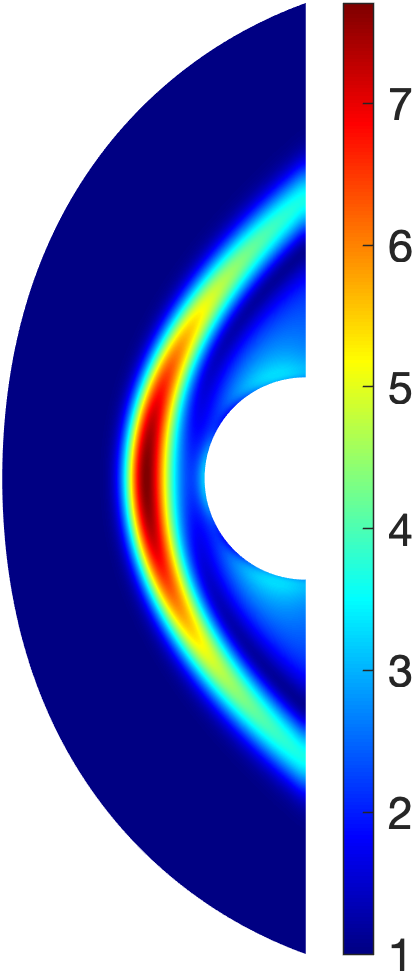}
  \caption{1st adaptive mesh}
	\end{subfigure}
	\hfill
	\begin{subfigure}[b]{0.3\textwidth}
		\centering		\includegraphics[width=0.66\textwidth]{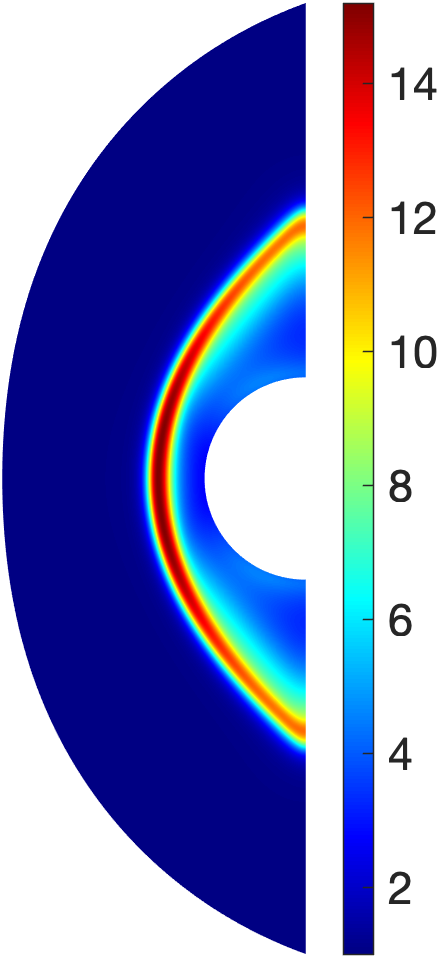}
  \caption{2nd adaptive mesh}
	\end{subfigure}
        \hfill
        \begin{subfigure}[b]{0.3\textwidth}
		\centering		\includegraphics[width=0.65\textwidth]{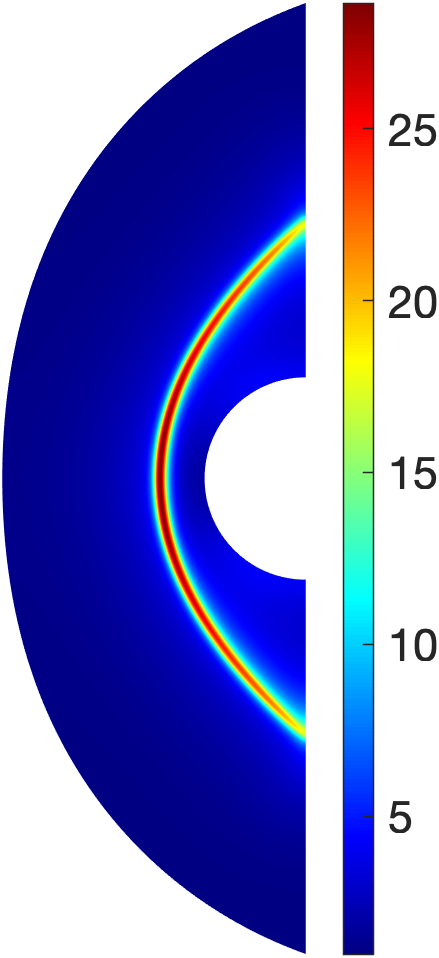}\caption{3rd adaptive mesh}	
	\end{subfigure} 
\caption{Mesh density functions used in the numerical solution of the Monge-Amp\`ere equation to generate the three adaptive meshes shown in Figure \ref{cyl70}.}
\label{cyl7a}	
\end{figure}

\begin{figure}[h!]
\centering
	\begin{subfigure}[b]{0.24\textwidth}
		\centering		\includegraphics[width=0.7\textwidth]{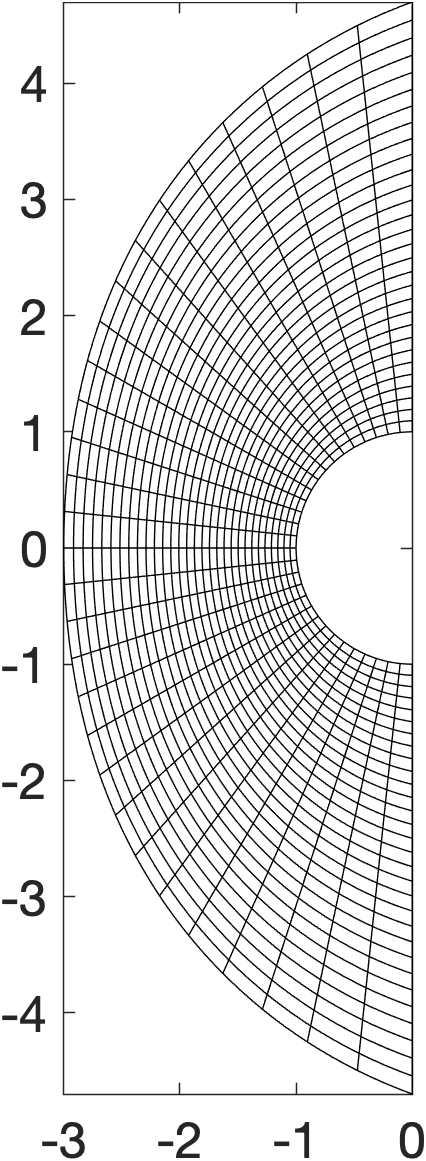}
  \caption{Initial mesh}
	\end{subfigure}
	\hfill
	\begin{subfigure}[b]{0.24\textwidth}
		\centering		\includegraphics[width=0.7\textwidth]{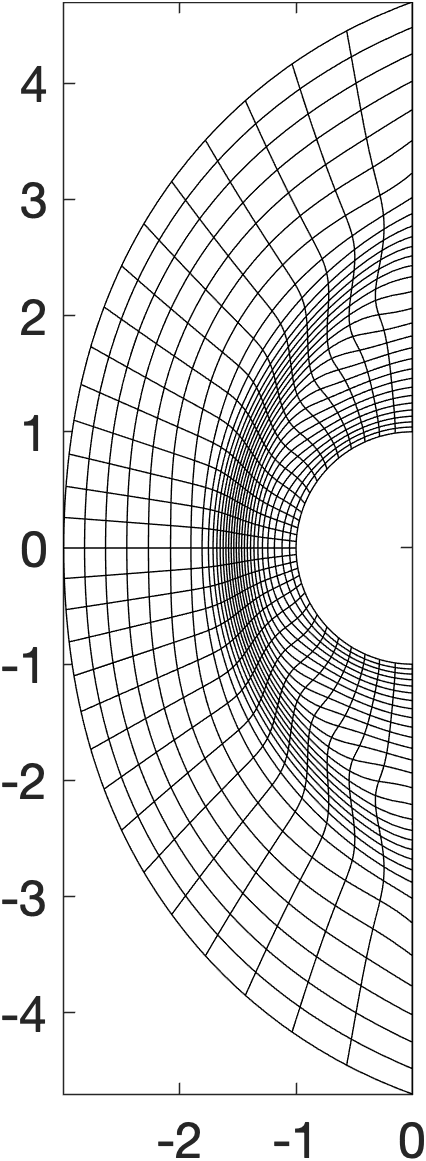}
  \caption{1st adaptive mesh}
	\end{subfigure}
        \hfill
        \begin{subfigure}[b]{0.24\textwidth}
		\centering		\includegraphics[width=0.7\textwidth]{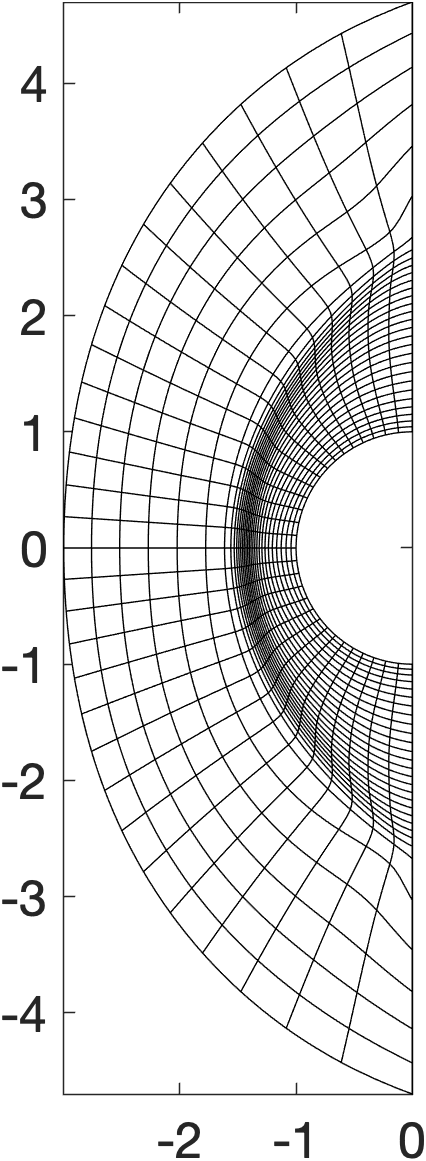}
  \caption{2nd adaptive mesh}
	\end{subfigure}
        \hfill
        \begin{subfigure}[b]{0.24\textwidth}
		\centering		\includegraphics[width=0.7\textwidth]{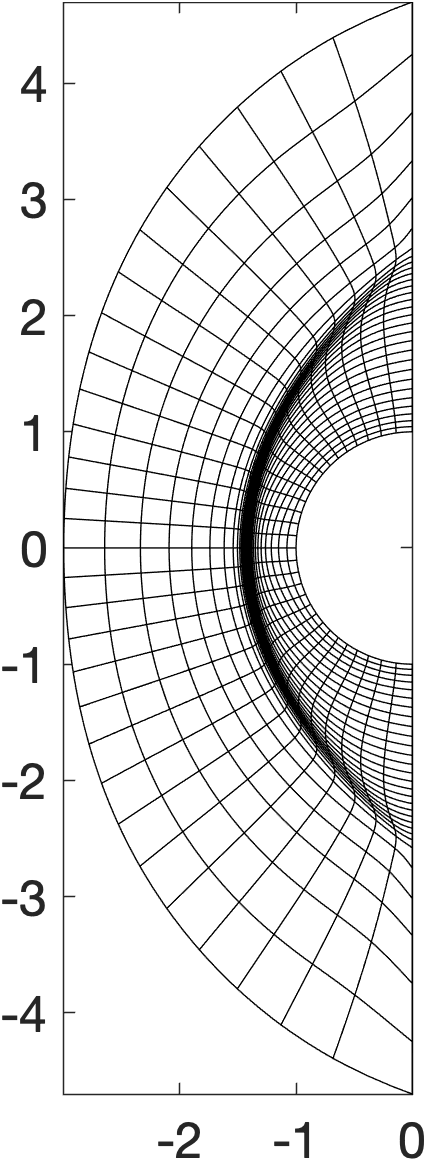}	
  \caption{3rd adaptive mesh}
	\end{subfigure} 
\caption{The initial mesh and three adaptive meshes for inviscid hypersonic flow past the circular cylinder at $M_\infty=7$. These meshes consist of 900 P4 quadrilateral elements.}
\label{cyl70}	
\end{figure}

This test case involves hypersonic flow past a unit circular cylinder at $M_\infty = 7$ and serves to demonstrate the effectiveness of our approach for strong bow shocks in the hypersonic regime . The cylinder wall is modeled with slip wall boundary condition. Supersonic outflow condition is used at the outlet, while supersonic inflow condition is imposed at the inlet. Figure \ref{cyl7a} shows the mesh density functions used in the numerical solution of the Monge-Amp\`ere equation to generate the three adaptive meshes shown in Figure \ref{cyl70}. These mesh density functions are computed from the mesh indicator (\ref{mdf2}) with using the numerical solutions on the initial mesh, the 1st adaptive mesh, and the 2nd adaptive mesh. We notice that the amplitude of the mesh density function increases with the mesh adaptation iteration because the numerical solution becomes sharper due to better resolution of the bow shock. This is because the optimal transport moves the elements toward the shock region and aligns them along the shock curves according to the mesh density function. 

\begin{figure}[htbp]
\centering
        \begin{subfigure}[b]{0.24\textwidth}
		\centering		\includegraphics[width=0.8\textwidth]{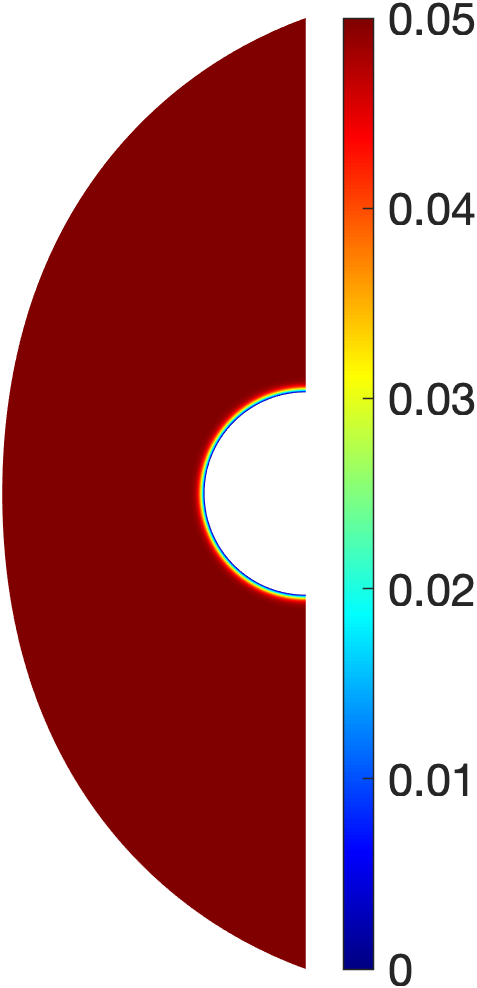}
	\end{subfigure}
	\hfill
	\begin{subfigure}[b]{0.24\textwidth}
		\centering		\includegraphics[width=0.8\textwidth]{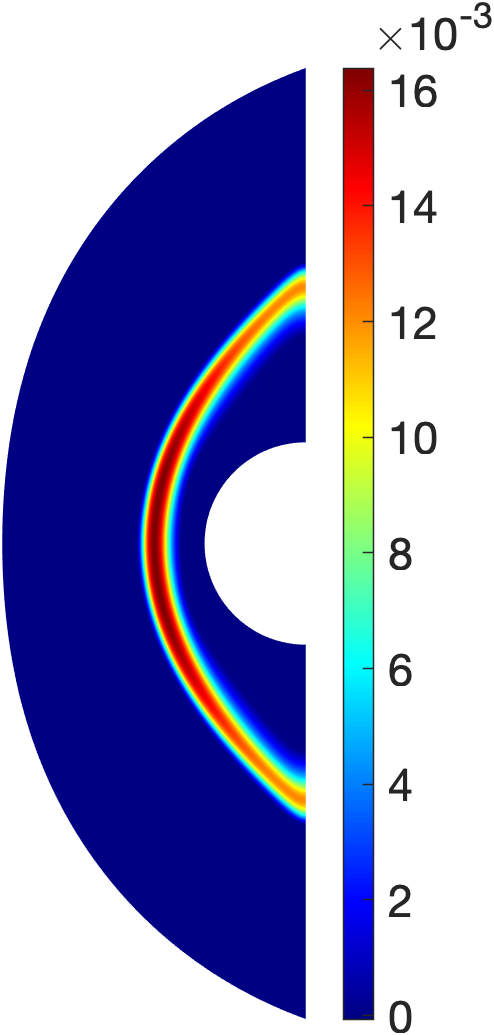}
	\end{subfigure}
        \hfill
        \begin{subfigure}[b]{0.24\textwidth}
		\centering		\includegraphics[width=0.8\textwidth]{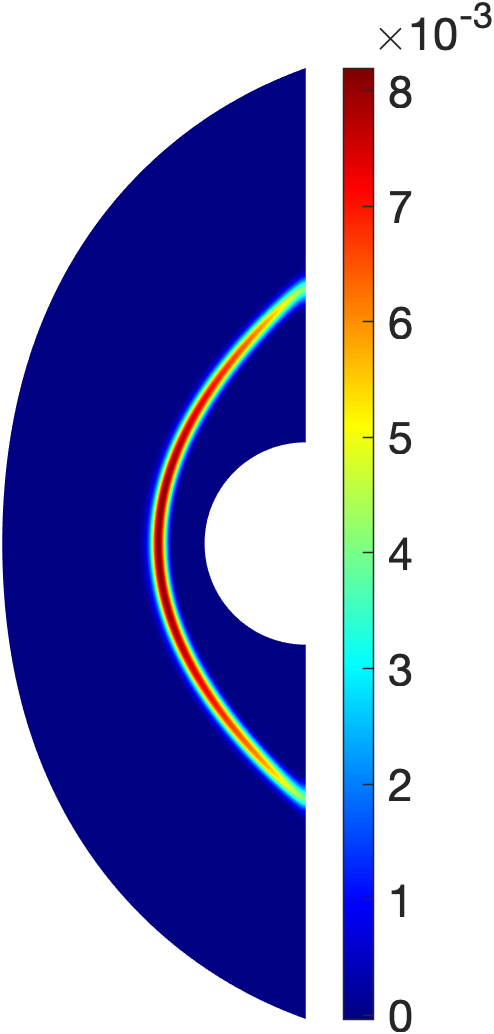}
	\end{subfigure}
        \hfill
        \begin{subfigure}[b]{0.24\textwidth}
		\centering		\includegraphics[width=0.8\textwidth]{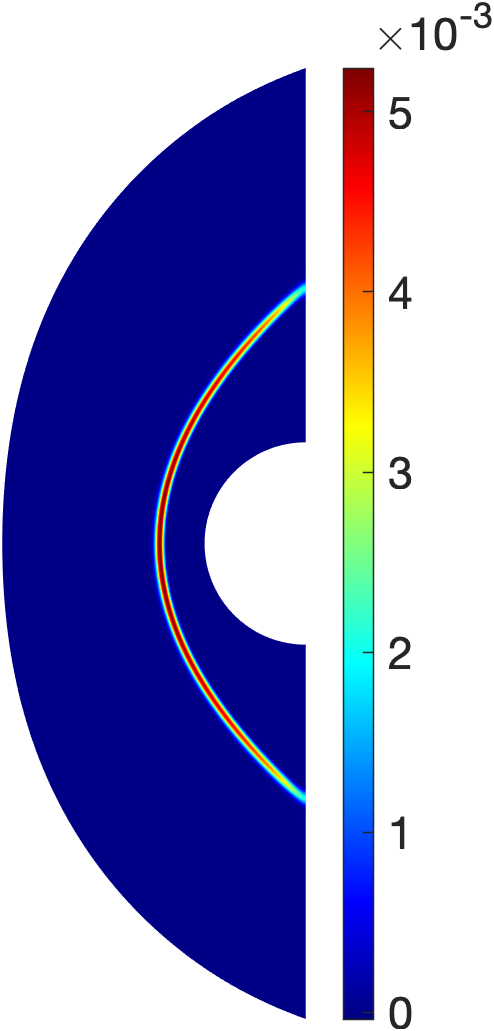}    
	\end{subfigure} \\[2ex]
         \begin{subfigure}[b]{0.24\textwidth}
		\centering		\includegraphics[width=0.7\textwidth]{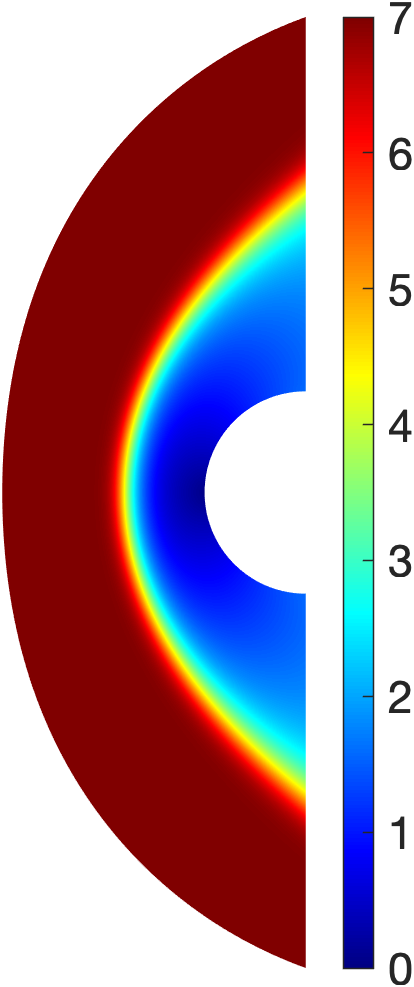}
	\end{subfigure}
	\hfill
	\begin{subfigure}[b]{0.24\textwidth}
		\centering		\includegraphics[width=0.7\textwidth]{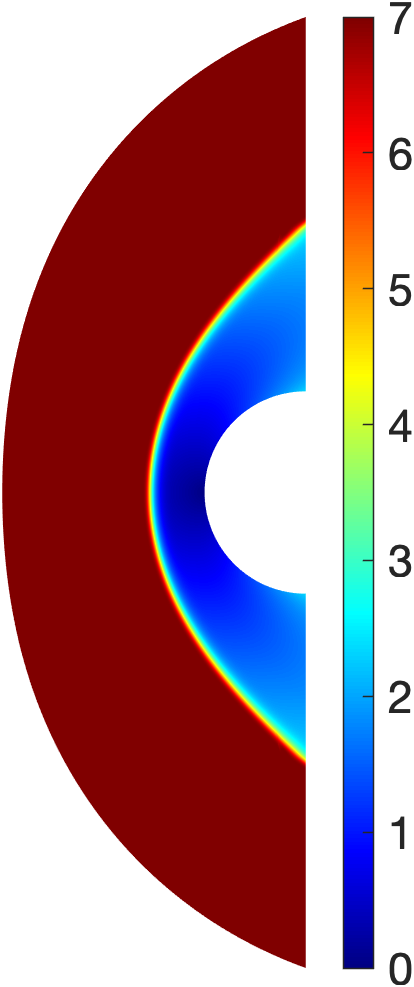}
	\end{subfigure}
        \hfill
        \begin{subfigure}[b]{0.24\textwidth}
		\centering		\includegraphics[width=0.7\textwidth]{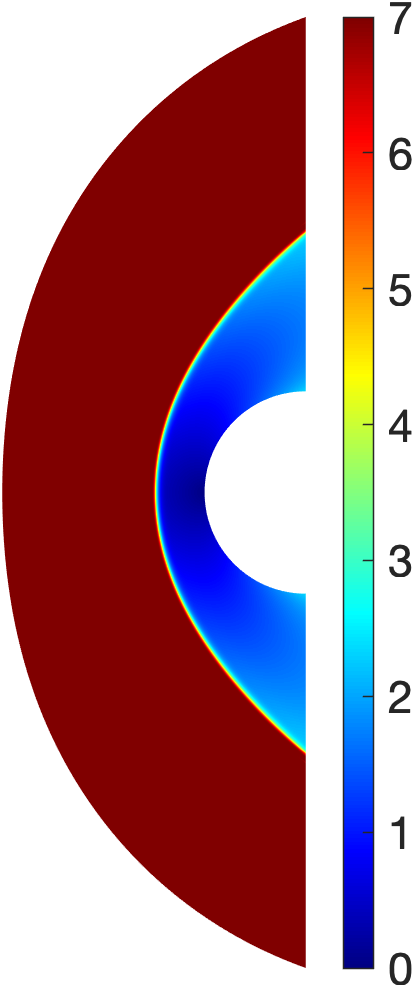}
	\end{subfigure}
        \hfill
        \begin{subfigure}[b]{0.24\textwidth}
		\centering		\includegraphics[width=0.7\textwidth]{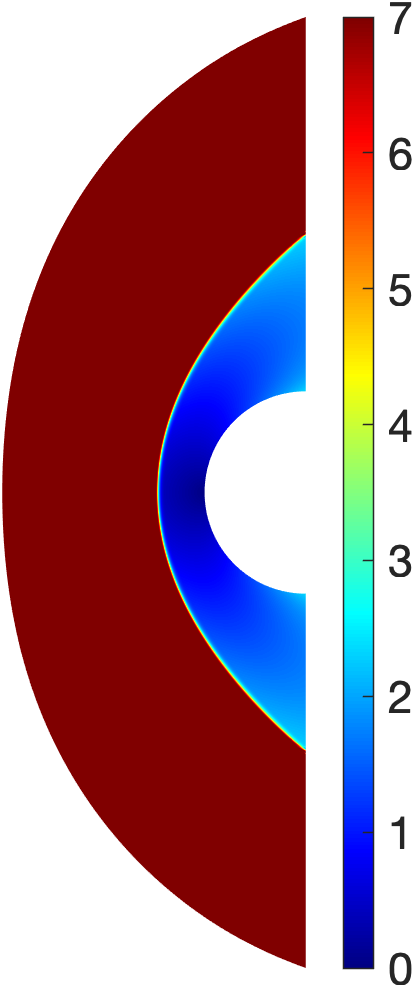}  
	\end{subfigure} \\[2ex]
        \begin{subfigure}[b]{0.24\textwidth}
		\centering		\includegraphics[width=0.77\textwidth]{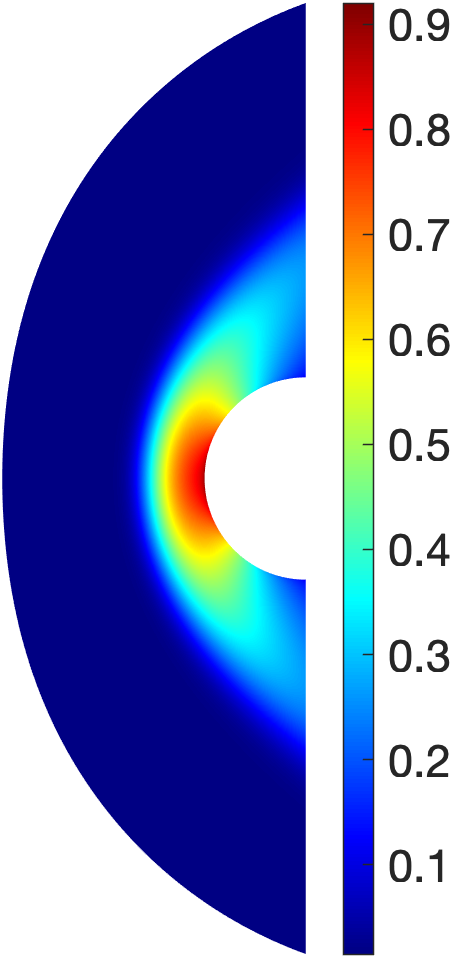}
  \caption{initial mesh}
	\end{subfigure}
	\hfill
	\begin{subfigure}[b]{0.24\textwidth}
		\centering		\includegraphics[width=0.77\textwidth]{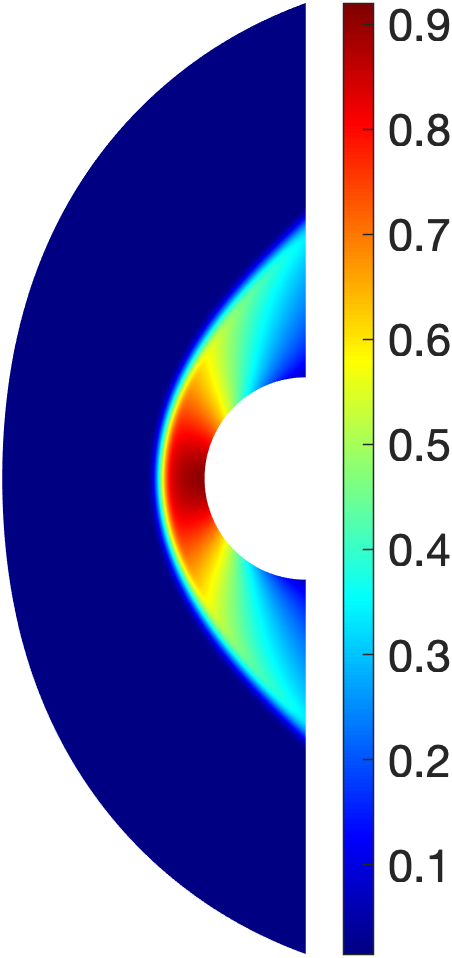}
  \caption{1st adaptive mesh}
	\end{subfigure}
        \hfill
        \begin{subfigure}[b]{0.24\textwidth}
		\centering		\includegraphics[width=0.77\textwidth]{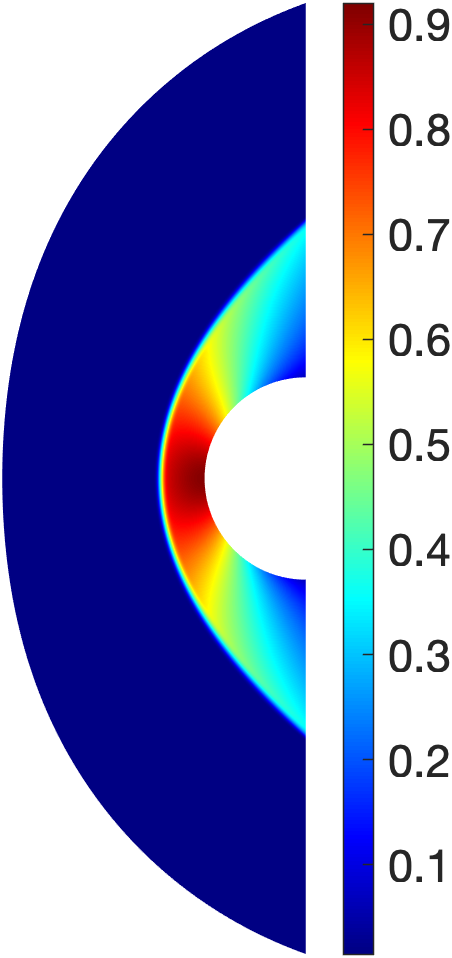}
  \caption{2nd adaptive mesh}
	\end{subfigure}
        \hfill
        \begin{subfigure}[b]{0.24\textwidth}
		\centering		\includegraphics[width=0.77\textwidth]{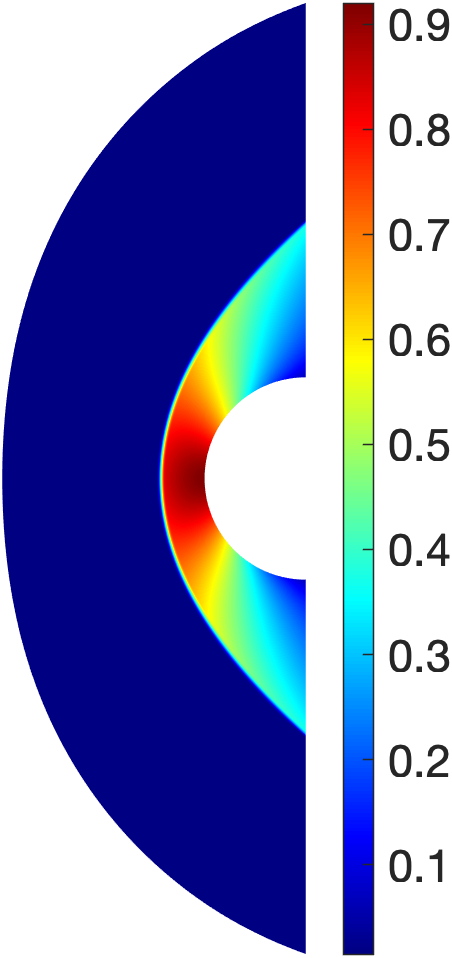}
  \caption{3rd adaptive mesh}
	\end{subfigure} 
\caption{Artificial viscosity (top row), Mach number (middle row), and pressure (bottom row) computed on the initial and adaptive meshes for inviscid hypersonic flow past the circular cylinder at $M_\infty=7$.}
\label{cyl71}	
\end{figure}

\begin{figure}[htbp]
\centering
	\begin{subfigure}[b]{0.49\textwidth}
		\centering		\includegraphics[width=\textwidth]{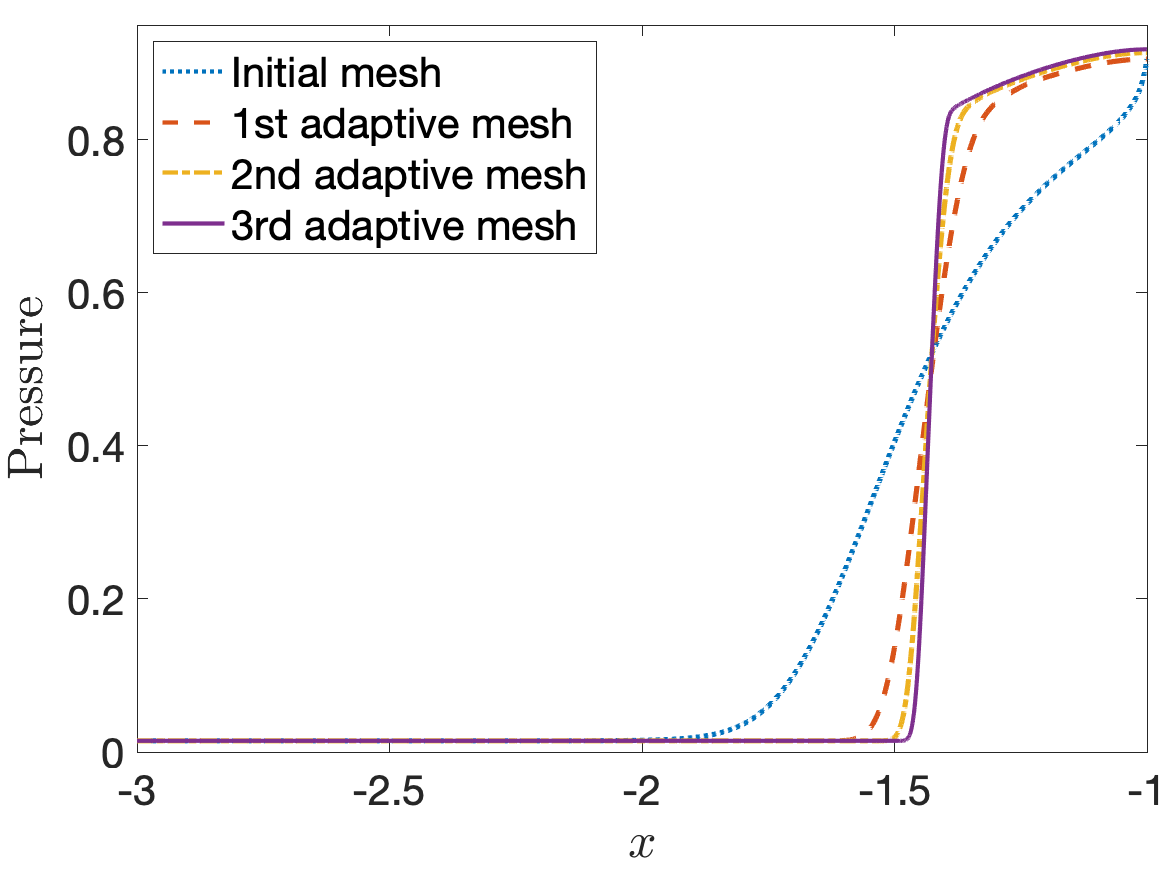}
	\end{subfigure}
	\hfill
	\begin{subfigure}[b]{0.49\textwidth}
		\centering		\includegraphics[width=\textwidth]{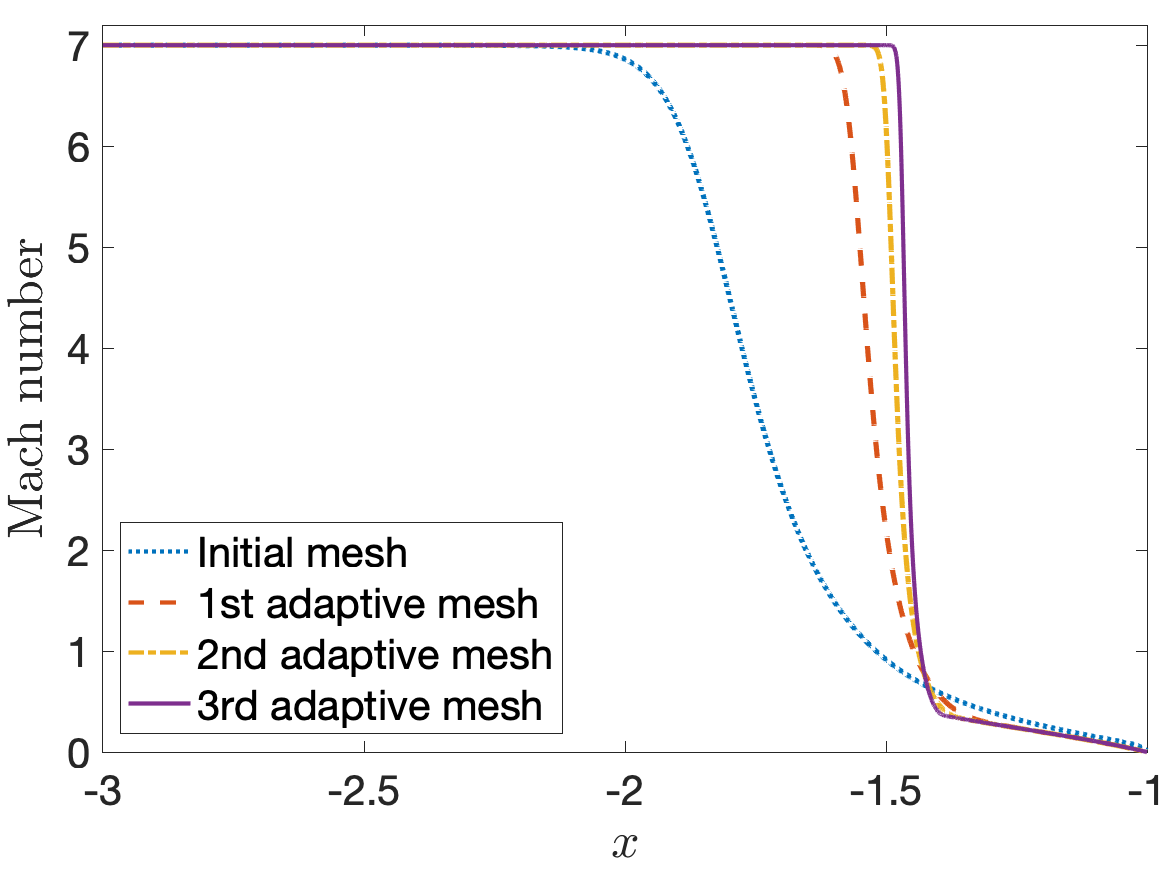}
	\end{subfigure}        
\caption{Profiles of pressure and Mach number along the line $y = 0$ for inviscid hypersonic flow past the circular cylinder at $M_\infty=7$.}
\label{cyl7b}	
\end{figure}

Figure \ref{cyl71} depicts the numerical solution computed on the initial and adaptive meshes. We see how the artificial viscosity fields are reduced in amplitude and width as the mesh adaptation procedure iterates. The numerical solution computed on the final adaptive mesh is clearly more accurate than those computed on the previous adaptive meshes. This can also be seen in Figure \ref{cyl7b} which shows the profiles of pressure and Mach number along the line $y=0$.  We see that these profiles converge rapidly with the adaptation iteration. The profiles computed on the second adaptive mesh are close to those computed on the third adaptive mesh, which are sharp and smooth. There is no oscillation and overshoot in the numerical solution on the final adaptive mesh. These results demonstrate the robustness of the proposed approach for strong bow shocks.

\subsection{Inviscid type IV shock-shock interaction}

Type IV Shock-shock interaction results in a very complex flow field with high pressure and heat flux peak in localized region. It occurs when the incident shock impinges on a bow shock and results in the formation of a supersonic impinging jet, a series of shock waves, expansion waves, and shear layers in a local area of interaction. The supersonic impinging jet, which is bounded by two shear layers separating the jet from the upper and lower subsonic regions, impinges on the body surface, and is terminated by a jet bow shock just ahead of the surface. This impinging jet bow shock wave creates a small stagnation region of high pressure and heating rates. Meanwhile, shear layers are formed to separate the supersonic jet from the lower and upper subsonic regions. 

Type IV hypersonic flows were experimentally studied by Wieting and Holden \cite{WIETING1989}. Over the years, many numerical methods have been used in the study of type IV shock-shock interaction \cite{Hsu1996,Nguyen2020gpu,Thareja1989,Yamamoto1998,Xu2005,Zhong1994a}. In the present work, we consider an inviscid type IV interaction with freestream Mach number $M_\infty = 8.03$. Based on the experimental measurement and the numerical calculations, Thareja et al. \cite{Thareja1989} summarized that the position of incident impinging shock on the cylinder can be approximated by the curve $y = 0.3271x + 0.4147$ for the experiment (Run 21) \cite{WIETING1989}. Boundary conditions are the same as those for the test case presented in Subsection 4.2, where the freemstream state $\bm u_\infty$ is represented by a hyperbolic tangent function to account for the incident impinging shock.

\begin{figure}[htbp]
\centering
	\begin{subfigure}[b]{0.17\textwidth}
		\centering		\includegraphics[width=\textwidth]{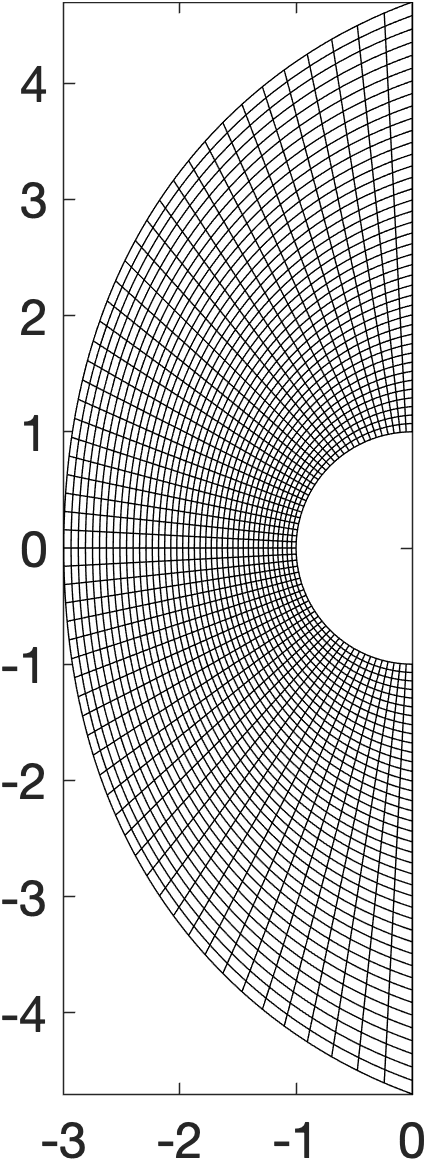}
        \caption{$\mathcal{T}_h^0$}
	\end{subfigure}
	\hfill
	\begin{subfigure}[b]{0.20\textwidth}
		\centering		\includegraphics[width=\textwidth]{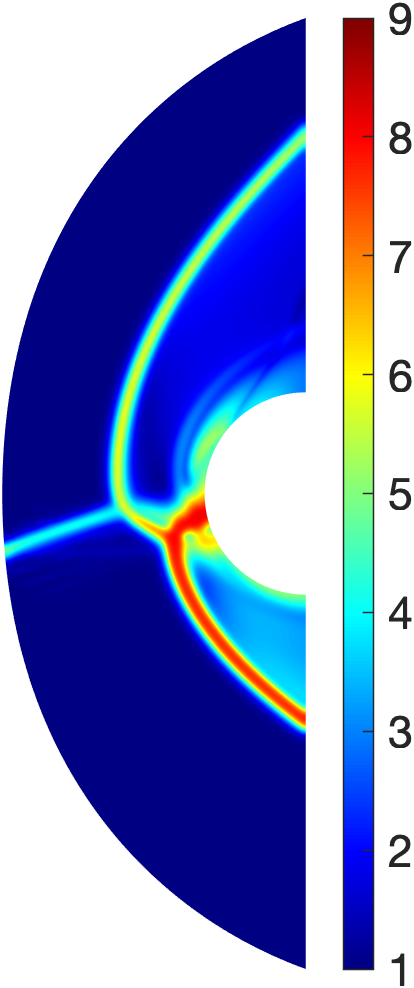}
  \caption{$\varrho'$ on $\mathcal{T}_h^0$}
	\end{subfigure}
        \hfill
        \begin{subfigure}[b]{0.20\textwidth}
		\centering		\includegraphics[width=0.85\textwidth]{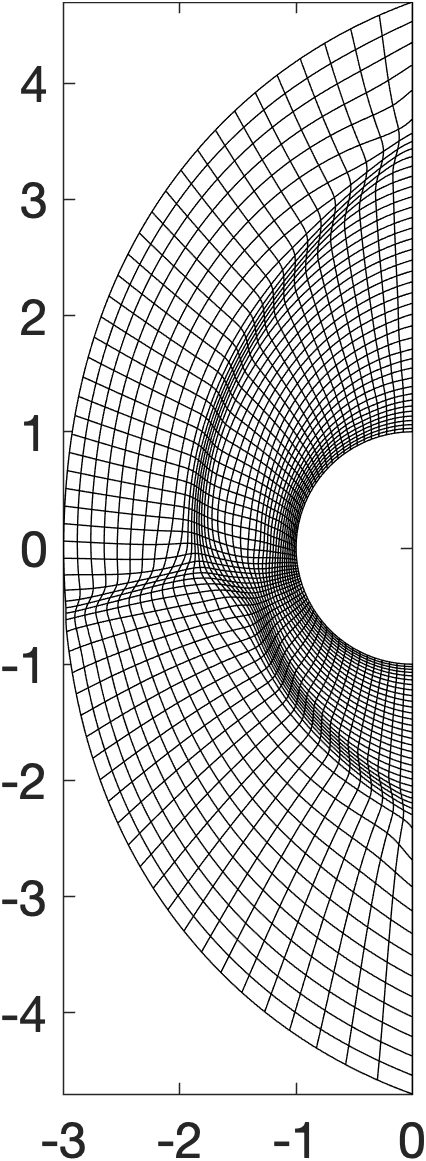}	
  \caption{$\mathcal{T}_h^1$}
	\end{subfigure} 
        \hfill
	\begin{subfigure}[b]{0.20\textwidth}
		\centering		\includegraphics[width=\textwidth]{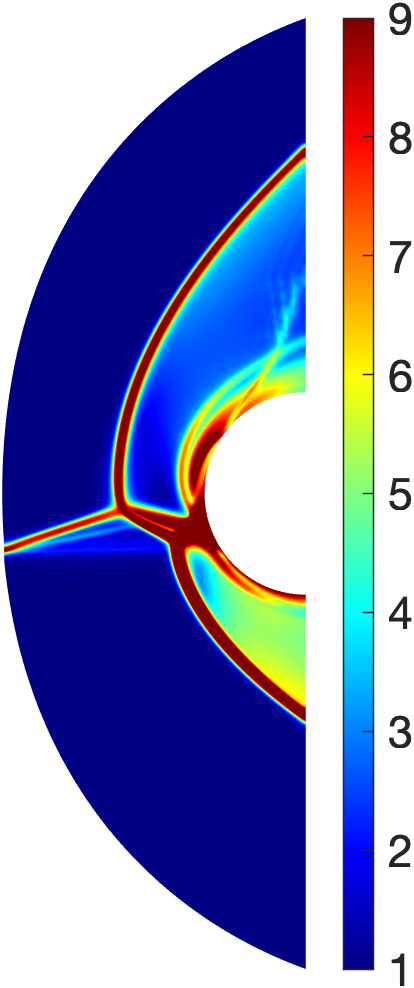}
  \caption{$\varrho'$ on $\mathcal{T}_h^1$}
	\end{subfigure}
        \hfill
        \begin{subfigure}[b]{0.17\textwidth}
		\centering		\includegraphics[width=\textwidth]{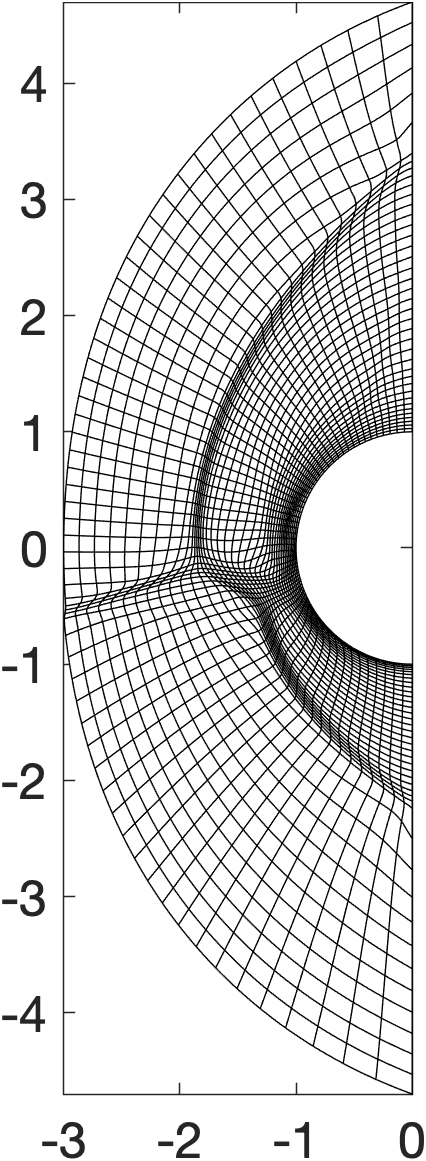}	
  \caption{$\mathcal{T}_h^2$}
	\end{subfigure} 
\caption{The initial mesh $\mathcal{T}_h^0$, the mesh density function $\varrho'$ on $\mathcal{T}_h^0$, the first adaptive mesh $\mathcal{T}_h^1$, the mesh density function $\varrho'$ on $\mathcal{T}_h^1$, the second adaptive mesh $\mathcal{T}_h^2$ for the inviscid type IV shock-shock interaction. These meshes consist of 2400 P4 quadrilaterals.}
\label{typeiva}	
\end{figure}

Figure \ref{typeiva} shows the initial and adaptive meshes as well as the mesh density functions used to obtain the adaptive meshes. The mesh density functions are computed from the mesh indicator (\ref{mdf2}) with using the numerical solutions on the initial mesh and the 1st adaptive mesh. The optimal transport moves the elements toward the shock region and aligns them along the shock curves. Furthermore, it also distributes elements around supersonic impinging jet, jet bow shock, expansion waves, and shear layers according to the mesh density function. As a result, the optimal transport can adapt meshes to capture complicated flow features without increasing the number of elements and modifying data structure.

\begin{figure}[h]
\centering
	\begin{subfigure}[b]{0.23\textwidth}
		\centering		\includegraphics[width=\textwidth]{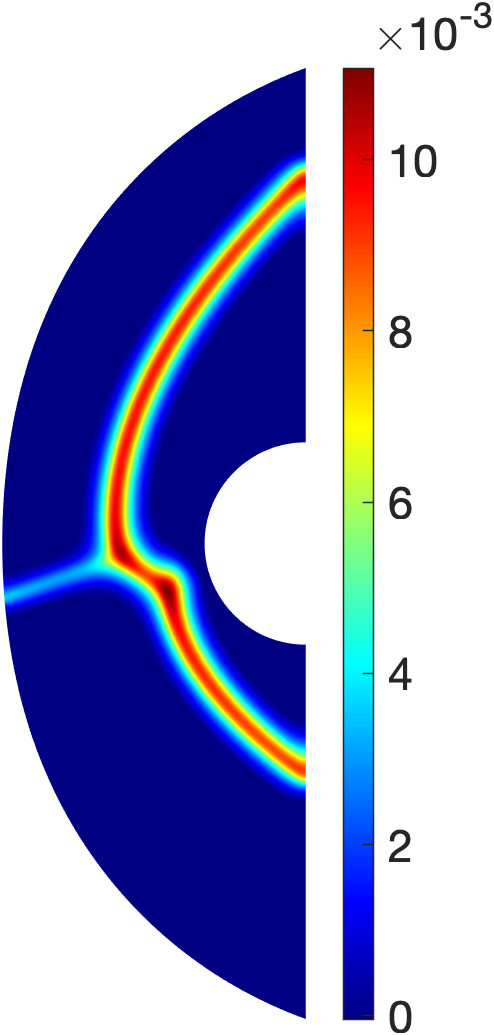}
  \caption{Artificial viscosity}
	\end{subfigure}
        \hfill
        \begin{subfigure}[b]{0.21\textwidth}
		\centering		\includegraphics[width=\textwidth]{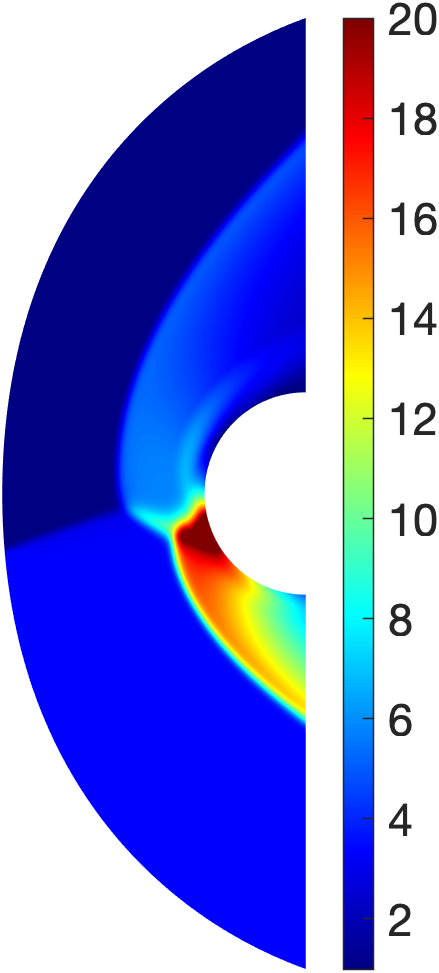}	
  \caption{Density}
	\end{subfigure} 
        \hfill
	\begin{subfigure}[b]{0.22\textwidth}
		\centering		\includegraphics[width=\textwidth]{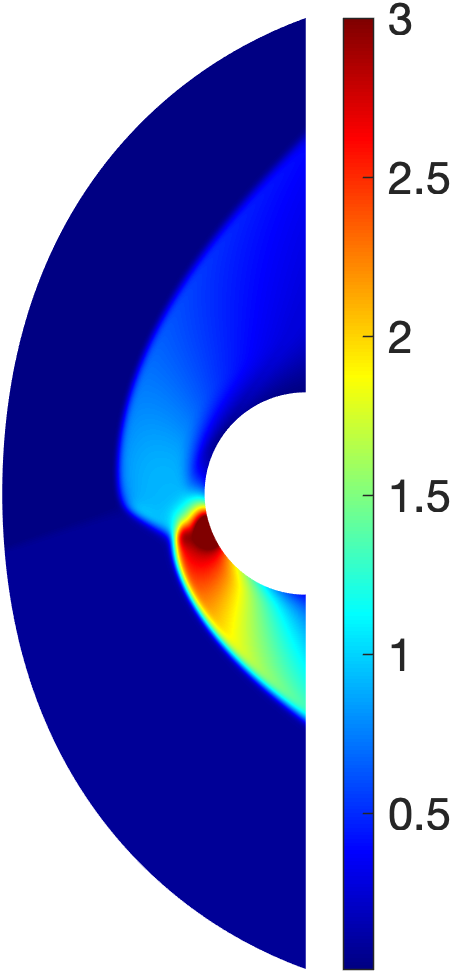}
  \caption{Pressure}
	\end{subfigure}
        \hfill
        \begin{subfigure}[b]{0.195\textwidth}
		\centering		\includegraphics[width=\textwidth]{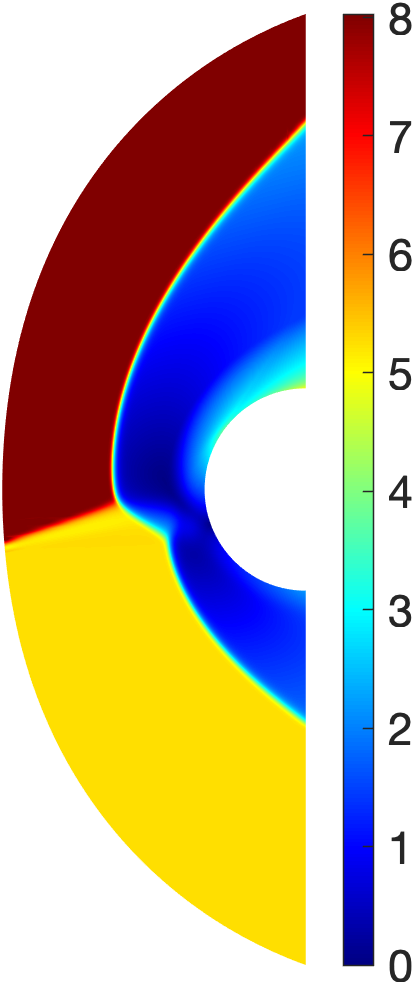}	
  \caption{Mach number}
	\end{subfigure} 
\caption{Numerical solution computed on the initial mesh for the inviscid type IV shock-shock interaction.}
\label{typeivc}	
\end{figure}

\begin{figure}[h!]
\centering
	\begin{subfigure}[b]{0.23\textwidth}
		\centering		\includegraphics[width=\textwidth]{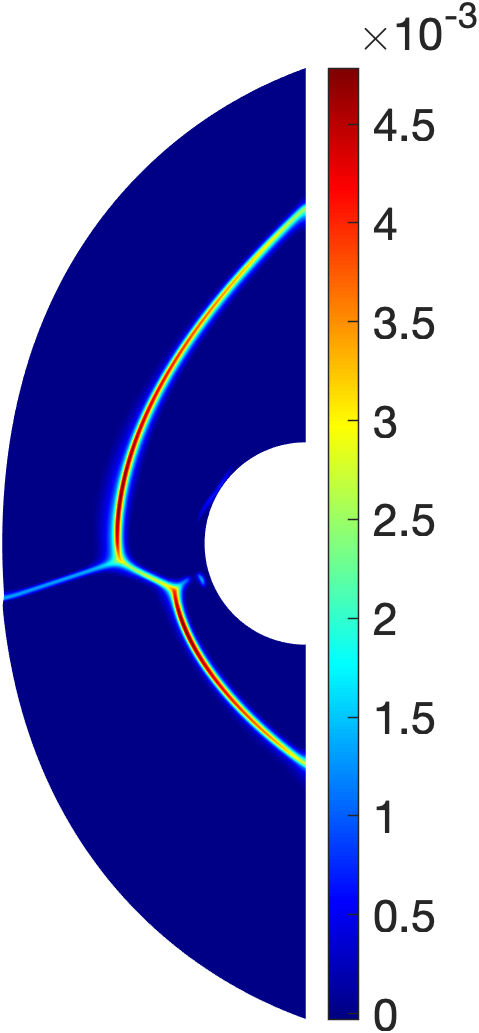}
  \caption{Artificial viscosity}
	\end{subfigure}
        \hfill
        \begin{subfigure}[b]{0.21\textwidth}
		\centering		\includegraphics[width=\textwidth]{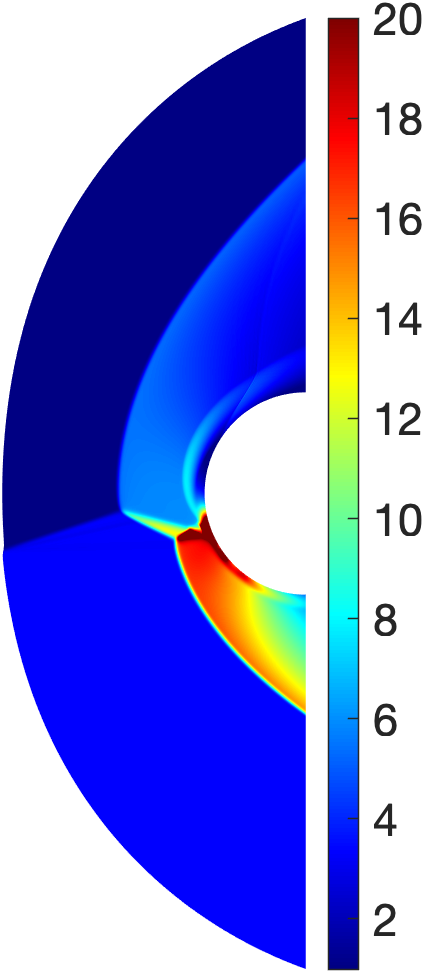}	
  \caption{Density}
	\end{subfigure} 
        \hfill
	\begin{subfigure}[b]{0.22\textwidth}
		\centering		\includegraphics[width=\textwidth]{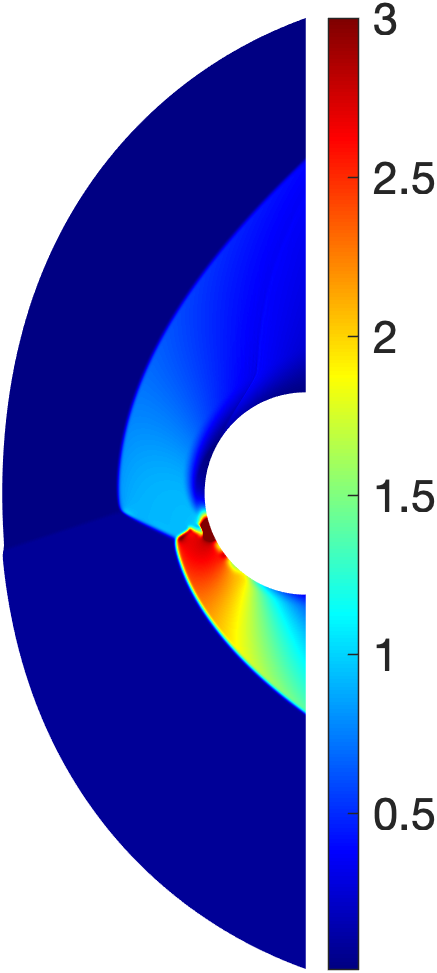}
  \caption{Pressure}
	\end{subfigure}
        \hfill
        \begin{subfigure}[b]{0.195\textwidth}
		\centering		\includegraphics[width=\textwidth]{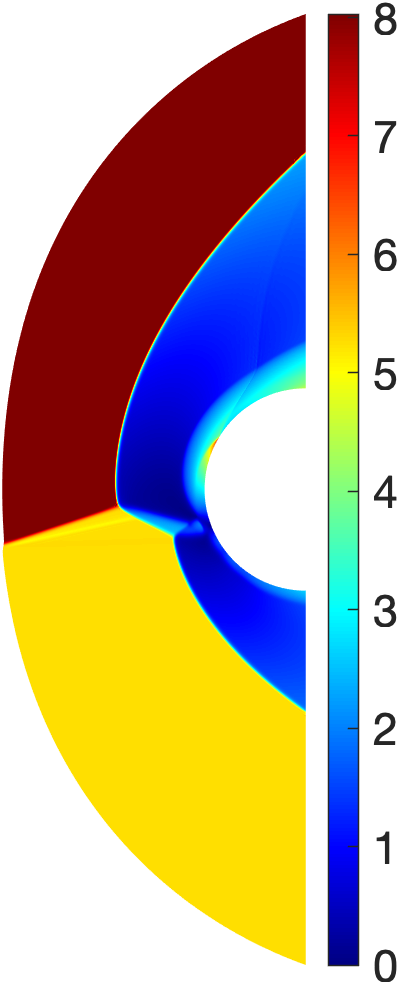}	
  \caption{Mach number}
	\end{subfigure} 
\caption{Numerical solution computed on the second adaptive mesh for the inviscid type IV shock-shock interaction.}
\label{typeivd}	
\end{figure}

We present the numerical solution computed on the initial mesh in Figure \ref{typeivc} and on the second adaptive mesh in Figure \ref{typeivd}. We notice that the numerical solution on the second adaptive mesh reveals supersonic impinging jet, jet bow shock, expansion waves, and shear layers of the flow, whereas the solution on the initial mesh does not possess some of these features. This is because the initial mesh does not have enough grid points to resolve those features even though it has the same number of elements as the second adaptive mesh. By redistributing the elements of the initial mesh to resolve shocks, impinging jet, jet bow shock, expansion waves, and shear layers, the optimal transport  considerably improves the numerical solution. This test case shows the ability of the optimal transport for dealing with complex shock flows.  

\begin{figure}[htbp]
	\centering
 \includegraphics[width=0.7\textwidth]{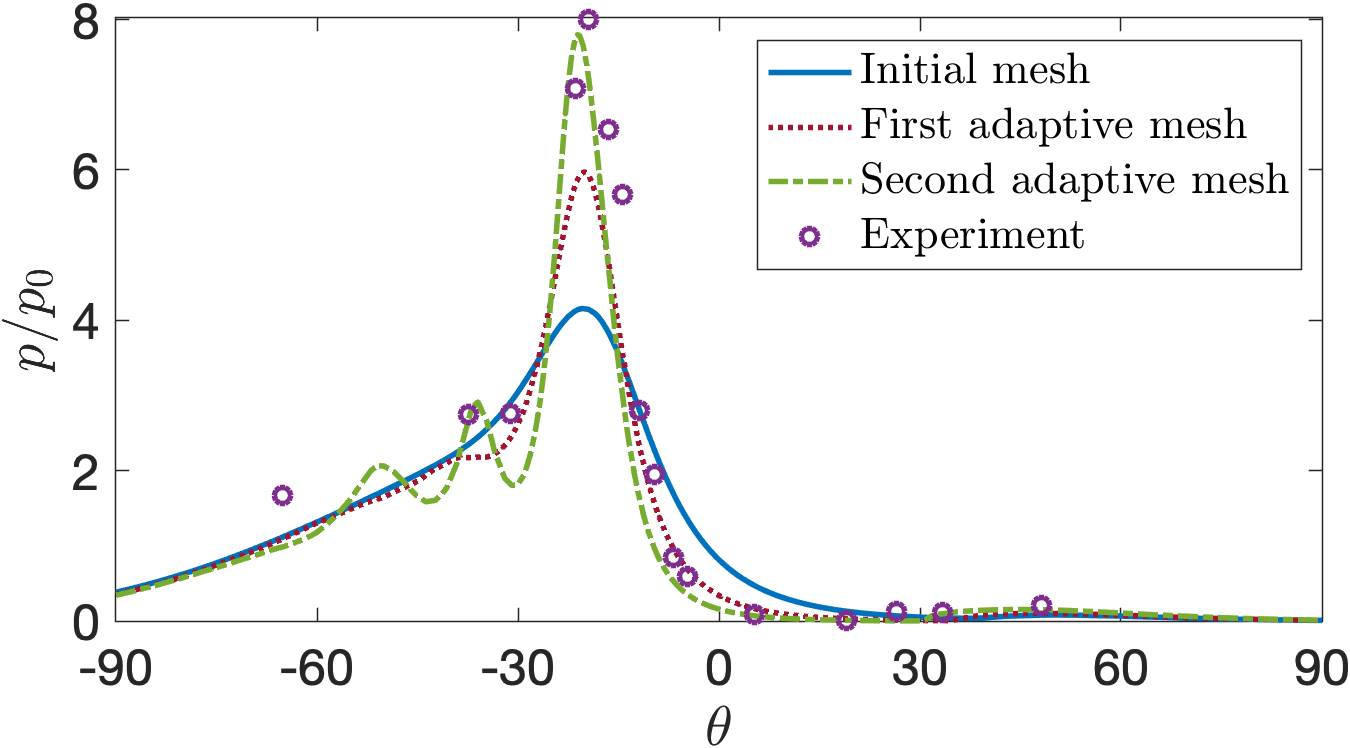}
 \caption{Profiles of the computed pressure ratio $p/p_0$ along the cylindrical surface for the type IV case. Here $p_0$ is the pressure at the stagnation point for inviscid hypersonic flow past the cylinder at $M_\infty = 8.03$. The symbols $\circ$ are the experimental data \cite{WIETING1989}.}
\label{typeive}	
\end{figure}

Finally, we present in Figure \ref{typeive} the profiles of the computed pressure along the cylindrical surface, where the symbols $\circ$ are the experimental data \cite{WIETING1989}. We see that the pressure profile computed on the second adaptive mesh has larger peak than those on the initial mesh and the first adaptive mesh. This is because the second adaptive mesh has a lot more elements in  the supersonic jet region than the initial mesh and the first adaptive mesh. As a result, the computed pressure on the second adaptive mesh agrees with the experimental measurement better than those on the other meshes.

\subsection{Viscous hypersonic flow past unit circular cylinder}

The last test case involves viscous hypersonic flow past unit circular cylinder at $M_\infty = 17.6$ and $Re=376,000$. The freestream temperature is $T_\infty = 200^{\rm o}$ K. The cylinder surface is isothermal with wall temperature $T_{\rm wall} = 500^{\rm o}$ K. Supersonic inflow and outflow boundary conditions are imposed at the inlet and outlet, respectively. This test case serves to demonstrate the ability of the optimal transport approach to deal with very strong bow shocks and extremely thin boundary layers. This problem was studied by Gnoffo and White \cite{Gnoffo2004} comparing the structured code LAURA and the unstructured code FUN3D. The simple geometry and strong shock make it a common benchmark case for assessing the performance of numerical methods and solution algorithms in hypersonic flow predictions \cite{Barter2010,Ching2019,Gnoffo2004a,Kitamura2013,Nguyen2020gpu}. This test case will demonstrate the ability of the optimal transport for dealing with very strong bow shock and extremely thin boundary layer.  

\begin{figure}[htbp]
\centering
	\begin{subfigure}[b]{0.24\textwidth}
		\centering		\includegraphics[width=0.68\textwidth]{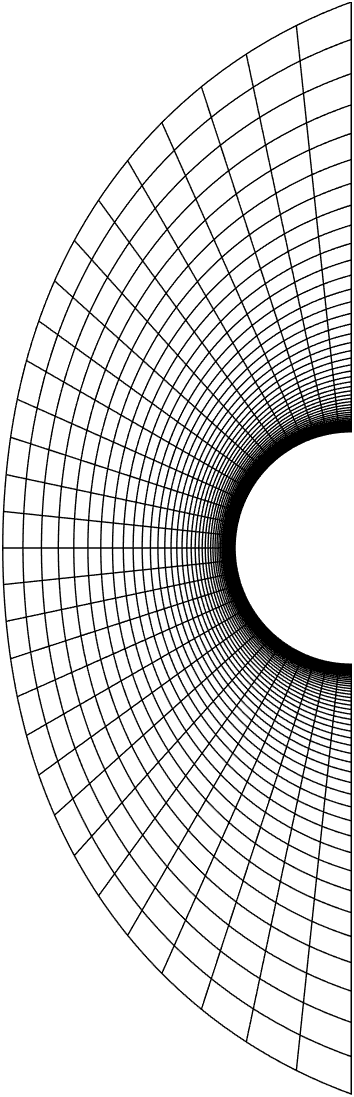}
    \caption{Initial mesh}
	\end{subfigure}
        \hfill
        \begin{subfigure}[b]{0.24\textwidth}
		\centering		\includegraphics[width=0.9\textwidth]{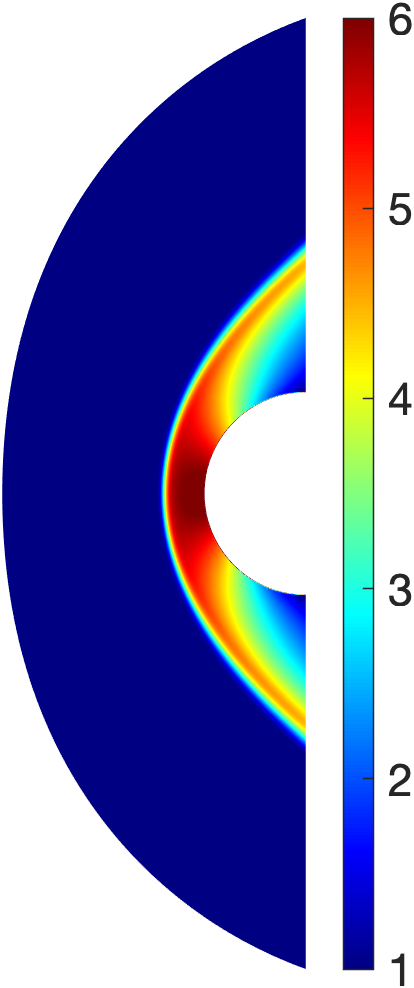}	
  \caption{$\bm u_h$ on the initial mesh}
	\end{subfigure} 
        \hfill
	\begin{subfigure}[b]{0.24\textwidth}
		\centering		\includegraphics[width=0.9\textwidth]{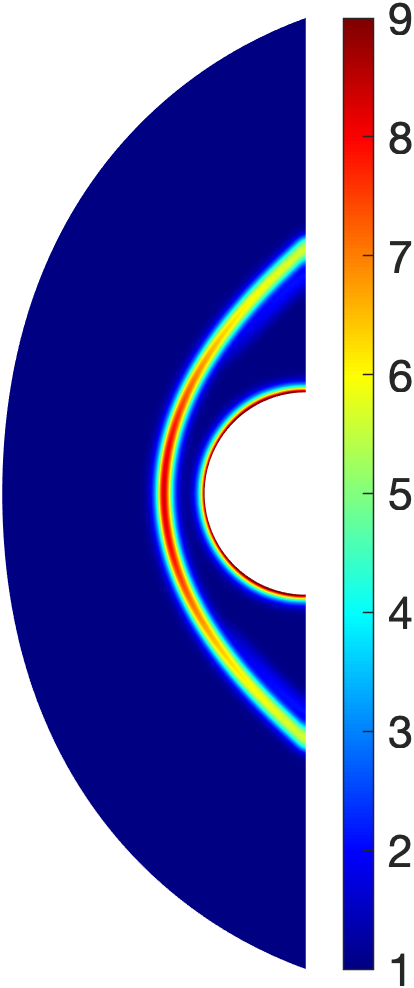}
  \caption{$\varrho'_h$ on the initial mesh}
	\end{subfigure}
        \hfill
        \begin{subfigure}[b]{0.24\textwidth}
		\centering		\includegraphics[width=0.68\textwidth]{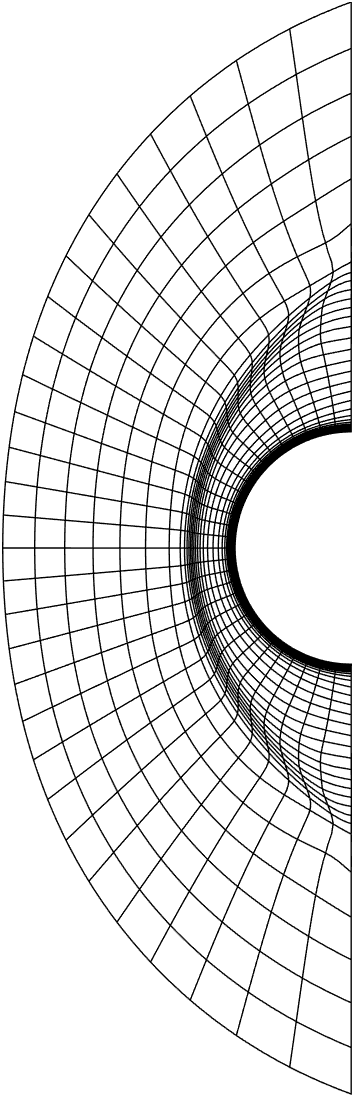}	
  \caption{Adaptive mesh}
	\end{subfigure} 
\caption{The initial mesh $\mathcal{T}_h^0$, the numerical solution $\bm u_h$ and the mesh density function $\varrho'$ on $\mathcal{T}_h^0$, and the  adaptive mesh $\mathcal{T}_h^1$ for the viscous hypersonic flow past a circular cylinder. These meshes consist of 1500 P4 quadrilaterals.}
\label{cyl18a}	
\end{figure}

\begin{figure}[htbp]
\centering
	\begin{subfigure}[b]{0.48\textwidth}
		\centering		\includegraphics[width=\textwidth]{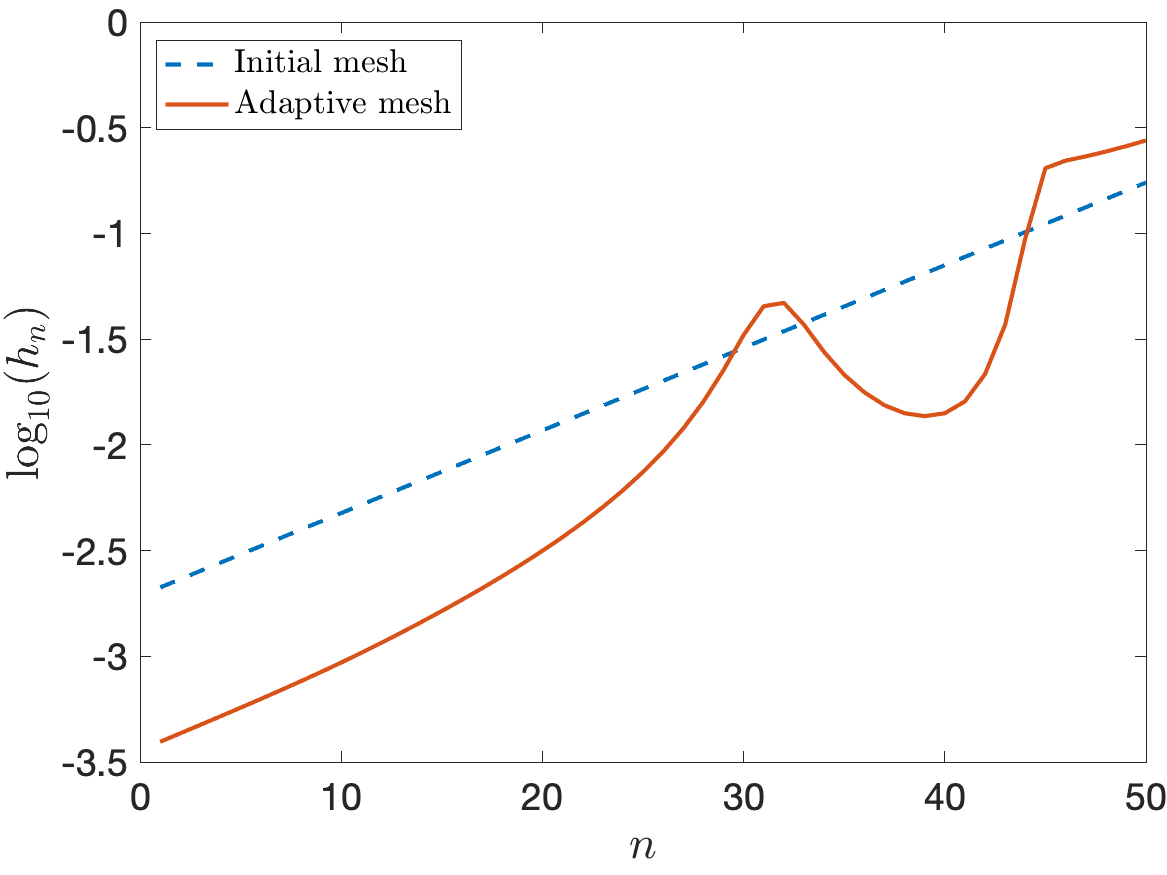}
	\end{subfigure}
	\hfill
	\begin{subfigure}[b]{0.50\textwidth}
		\centering		\includegraphics[width=\textwidth]{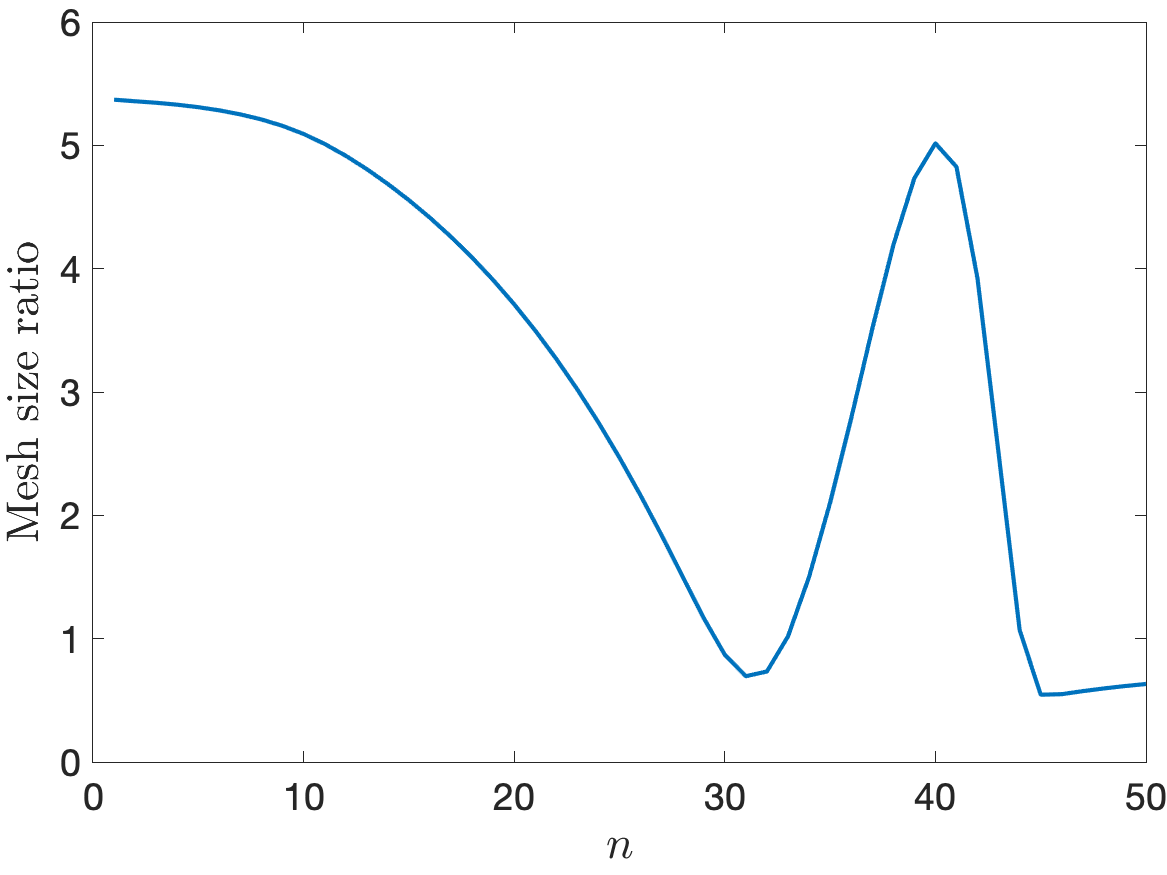}
	\end{subfigure}        
\caption{Logarithm with base 10 of the mesh size $h_n$ along the line $y = 0$, where the subscript $n$ indicate the element number starting from the cylinder wall. Right figure shows the mesh size ratio $h_n^{\rm initial}/h_n^{\rm adaptive}$ between the initial and adaptive mesh.}
\label{cyl18b}	
\end{figure}

Figure \ref{cyl18a} shows the initial and adaptive meshes as well as the mesh density function used to obtain the adaptive mesh.  The mesh density function is computed from the mesh indicator (\ref{mdf2}) with using the numerical solutions on the initial mesh. The optimal transport moves the elements of the initial mesh toward the shock and the boundary layer regions because the mesh density function is high in those regions. As a result, the optimal transport can adapt meshes to capture  shocks and resolve boundary layers. To see this feature more clearly, in Figure \ref{cyl18b}, we plot $\log_{10}(h_n)$ as a function of $n$ for both the initial mesh and the adaptive mesh, where $h_n$ denotes the element size of an $n$th element starting from the cylinder wall along the horizontal line $y=0$. We see that the adaptive mesh has smaller element sizes than the initial mesh near the wall and in the shock region. As a result, the adaptive mesh should be able to resolve the boundary layer and shock better than the initial mesh.

\begin{figure}[htbp]
\centering
	\begin{subfigure}[b]{0.22\textwidth}
		\centering		\includegraphics[width=\textwidth]{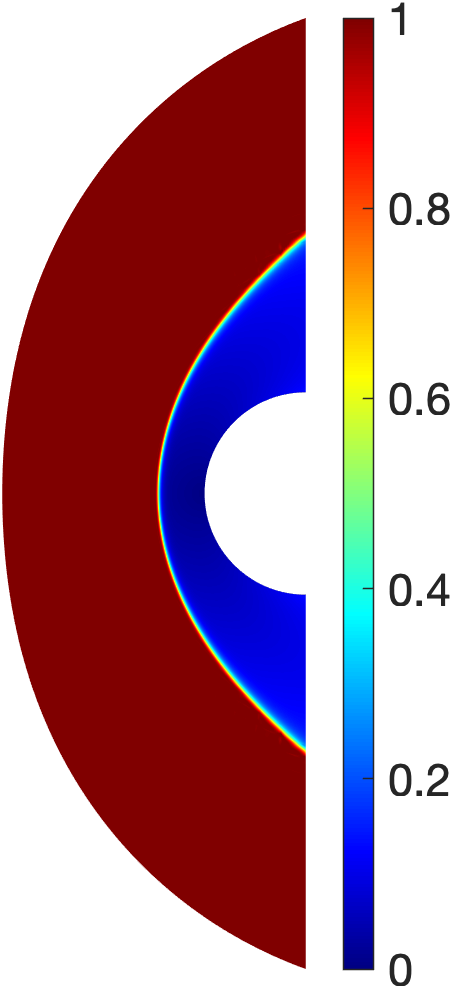}
  \caption{$M/M_\infty$}
	\end{subfigure}
        \hfill
        \begin{subfigure}[b]{0.2\textwidth}
		\centering		\includegraphics[width=\textwidth]{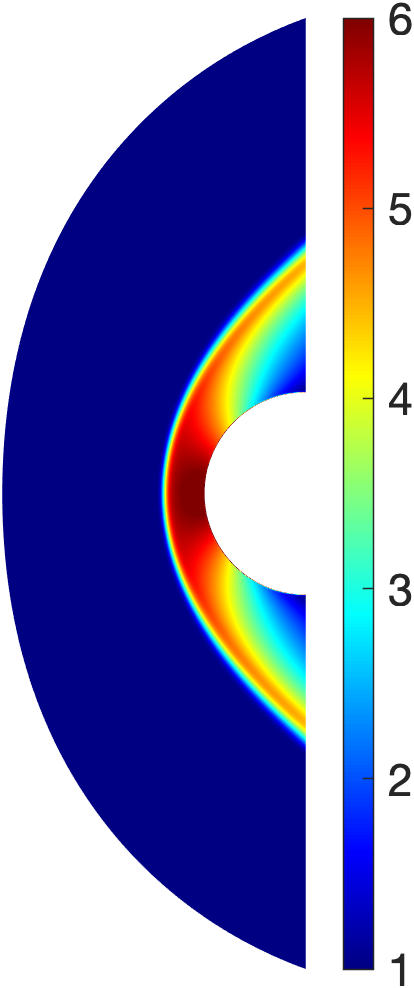}	
  \caption{$\rho/\rho_\infty$}
	\end{subfigure} 
        \hfill
	\begin{subfigure}[b]{0.23\textwidth}
		\centering		\includegraphics[width=\textwidth]{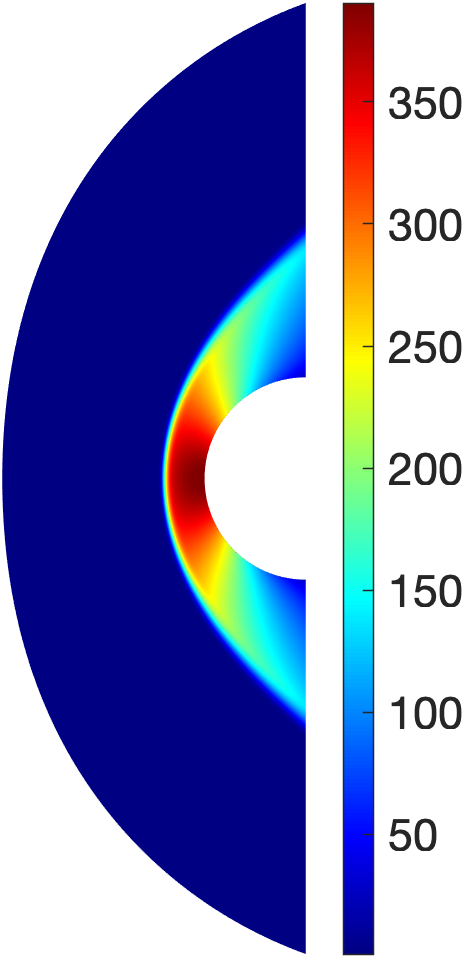}
  \caption{$p/p_\infty$}
	\end{subfigure}
        \hfill
        \begin{subfigure}[b]{0.22\textwidth}
		\centering		\includegraphics[width=\textwidth]{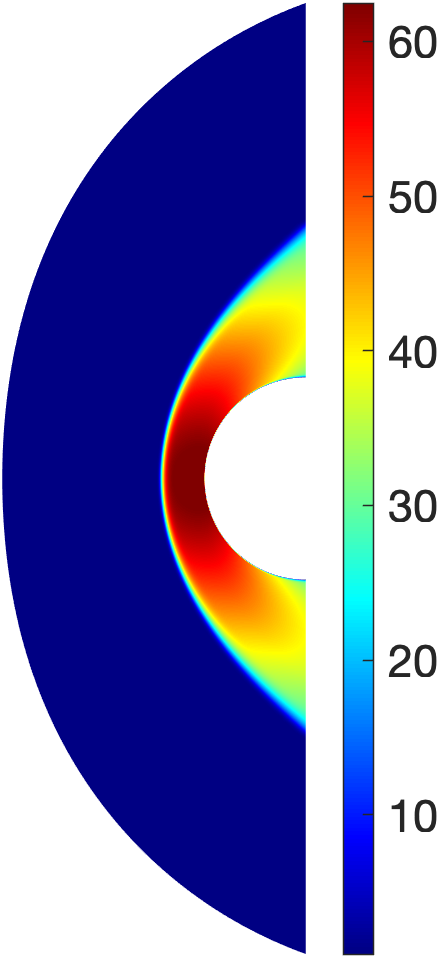}	
  \caption{$T/T_\infty$}
	\end{subfigure} 
\caption{Numerical solution computed on the  initial mesh for the viscous hypersonic flow past a circular cylinder.}
\label{cyl18c}	
\end{figure}

\begin{figure}[htbp]
\centering
	\begin{subfigure}[b]{0.22\textwidth}
		\centering		\includegraphics[width=\textwidth]{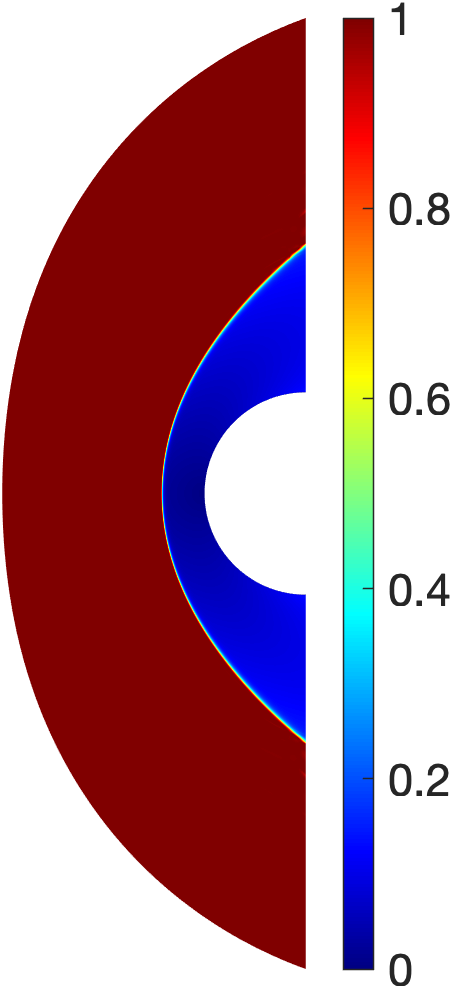}
  \caption{$M/M_\infty$}
	\end{subfigure}
        \hfill
        \begin{subfigure}[b]{0.2\textwidth}
		\centering		\includegraphics[width=\textwidth]{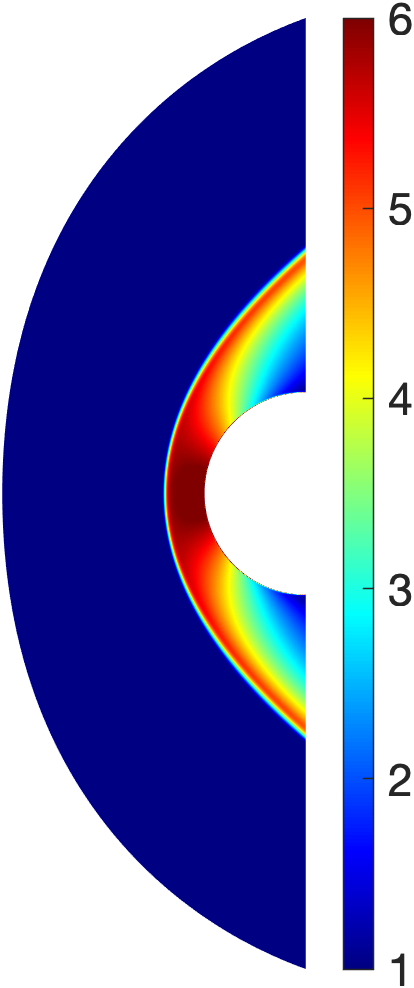}	
  \caption{$\rho/\rho_\infty$}
	\end{subfigure} 
        \hfill
	\begin{subfigure}[b]{0.23\textwidth}
		\centering		\includegraphics[width=\textwidth]{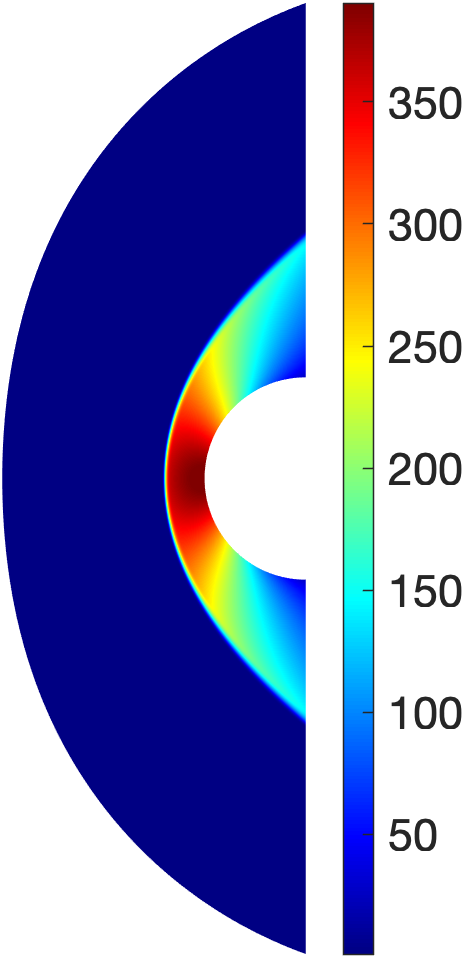}
  \caption{$p/p_\infty$}
	\end{subfigure}
        \hfill
        \begin{subfigure}[b]{0.22\textwidth}
		\centering		\includegraphics[width=\textwidth]{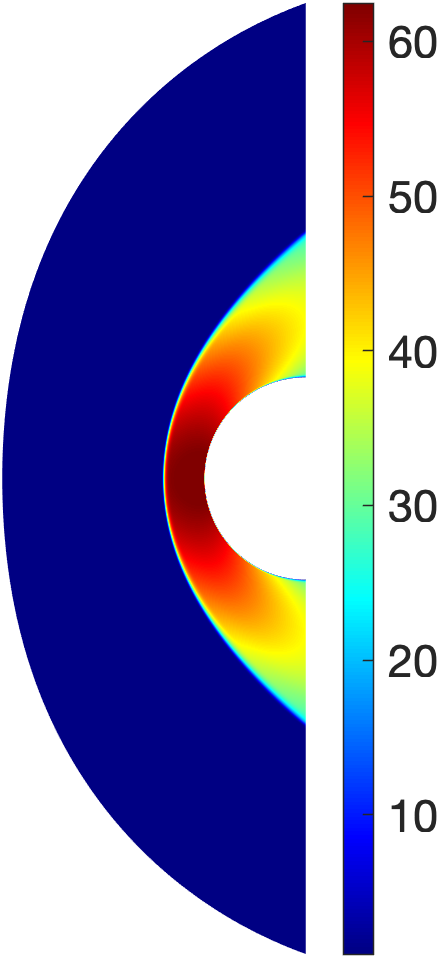}	
  \caption{$T/T_\infty$}
	\end{subfigure} 
\caption{Numerical solution computed on the  adaptive mesh for the viscous hypersonic flow past a circular cylinder.}
\label{cyl18d}	
\end{figure}

\begin{figure}[htbp]
\centering
	\begin{subfigure}[b]{0.49\textwidth}
		\centering		\includegraphics[width=\textwidth]{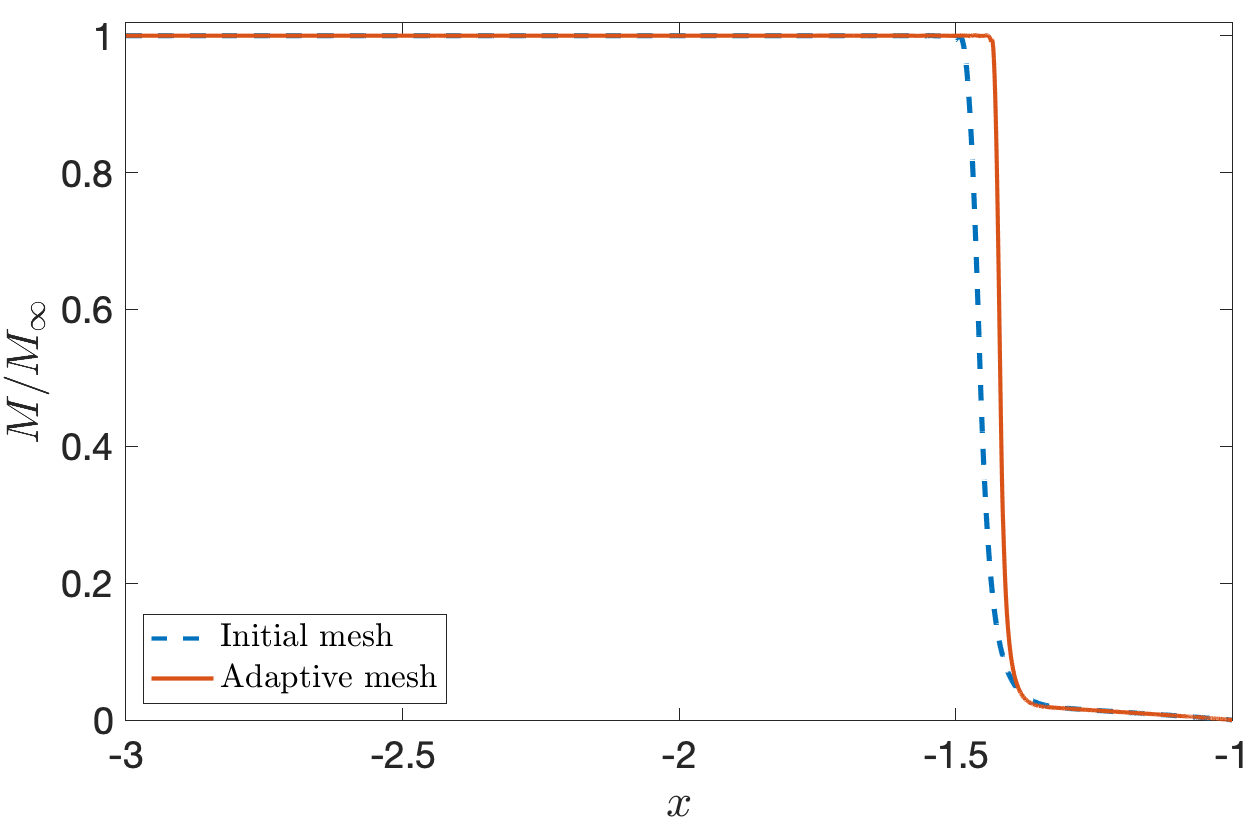}
	\end{subfigure}
	\hfill
	\begin{subfigure}[b]{0.49\textwidth}
		\centering		\includegraphics[width=\textwidth]{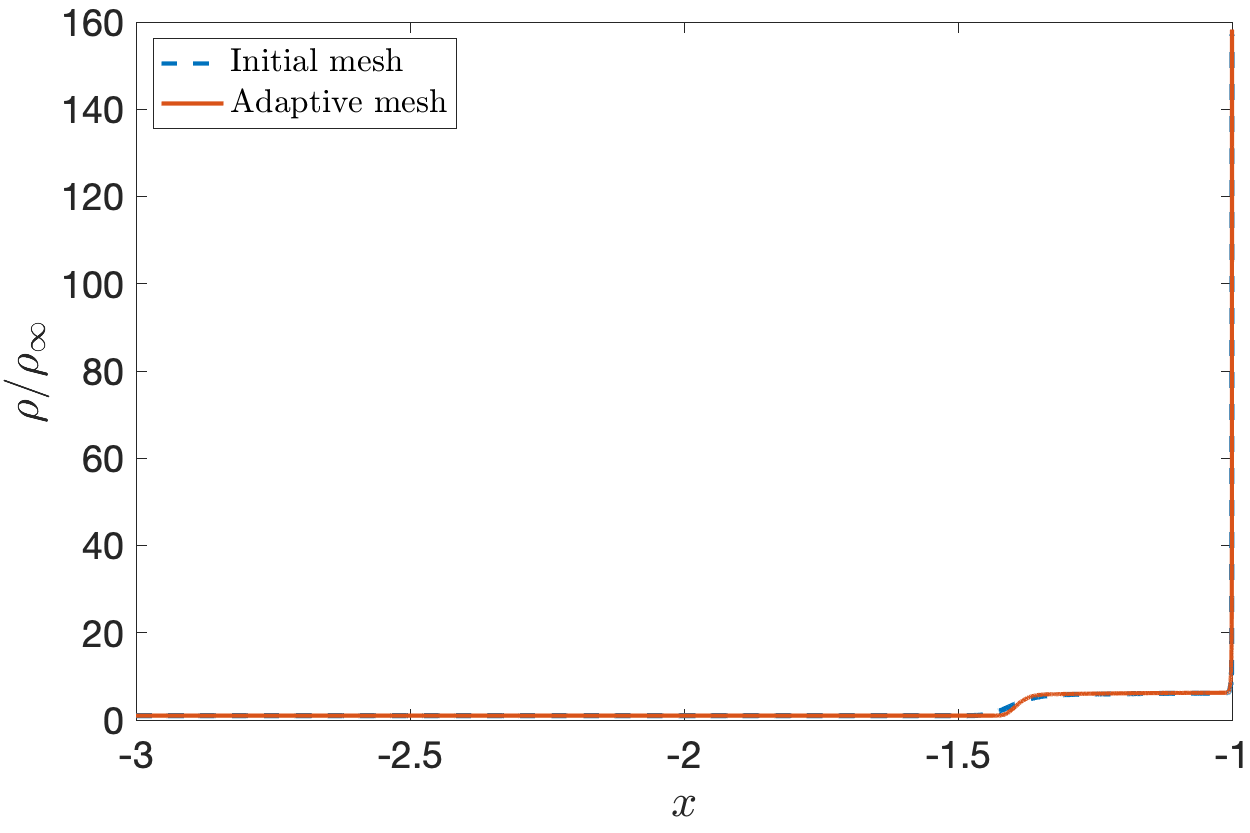}
	\end{subfigure}   \\
 \hfill
	\begin{subfigure}[b]{0.49\textwidth}
		\centering		\includegraphics[width=\textwidth]{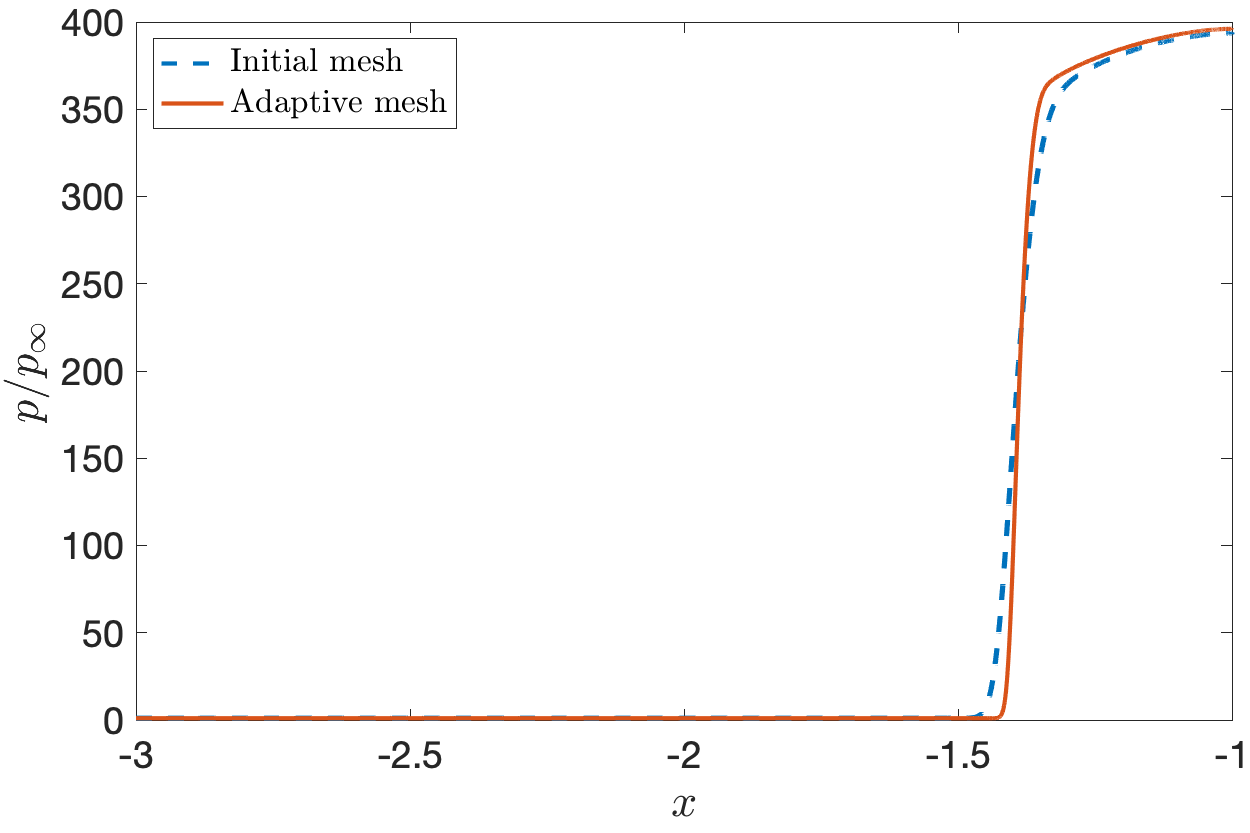}
	\end{subfigure}        
 \hfill
	\begin{subfigure}[b]{0.49\textwidth}
		\centering		\includegraphics[width=\textwidth]{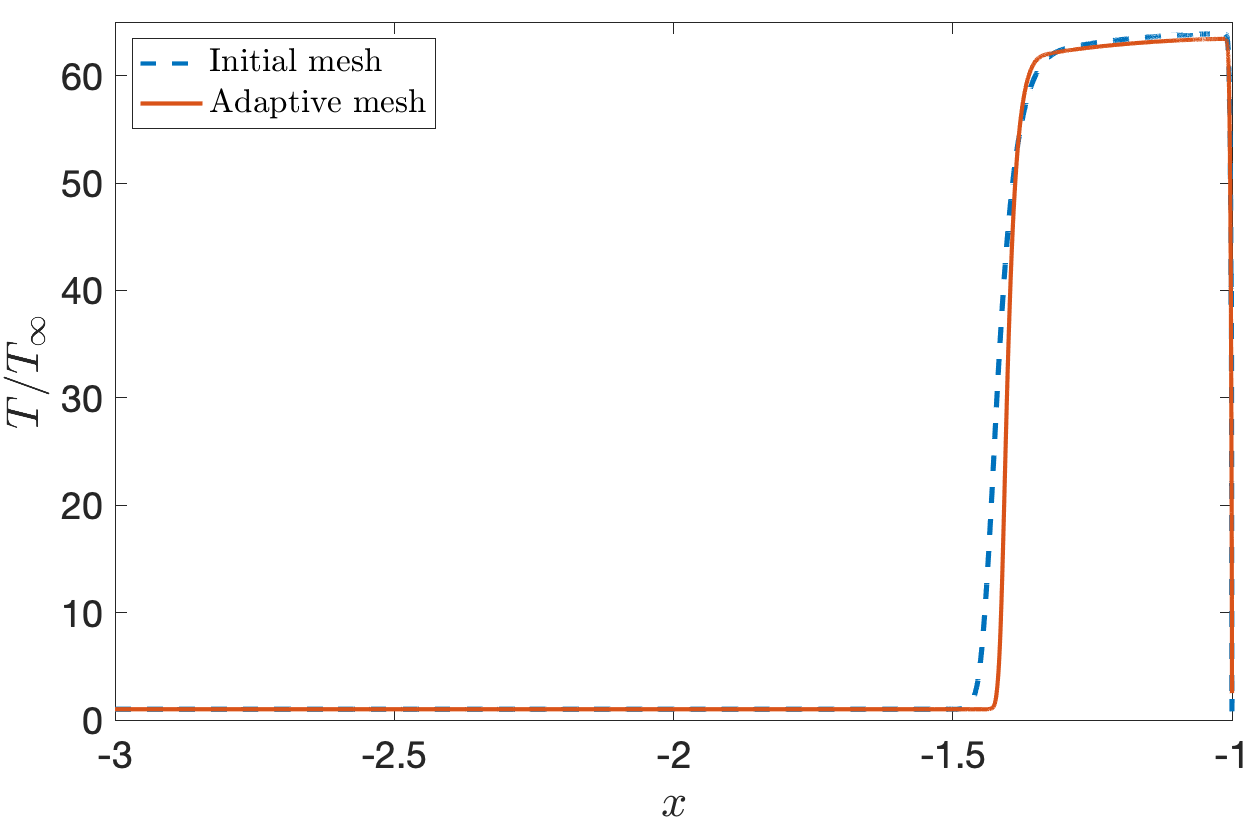}
	\end{subfigure}        
\caption{Profiles of the Mach number, density, pressure, and temperature along the line $y = 0$ for the viscous hypersonic flow past a circular cylinder..}
\label{cyl18e}	
\end{figure}

\begin{figure}[htbp]
\centering
	\begin{subfigure}[b]{0.48\textwidth}
		\centering		\includegraphics[width=\textwidth]{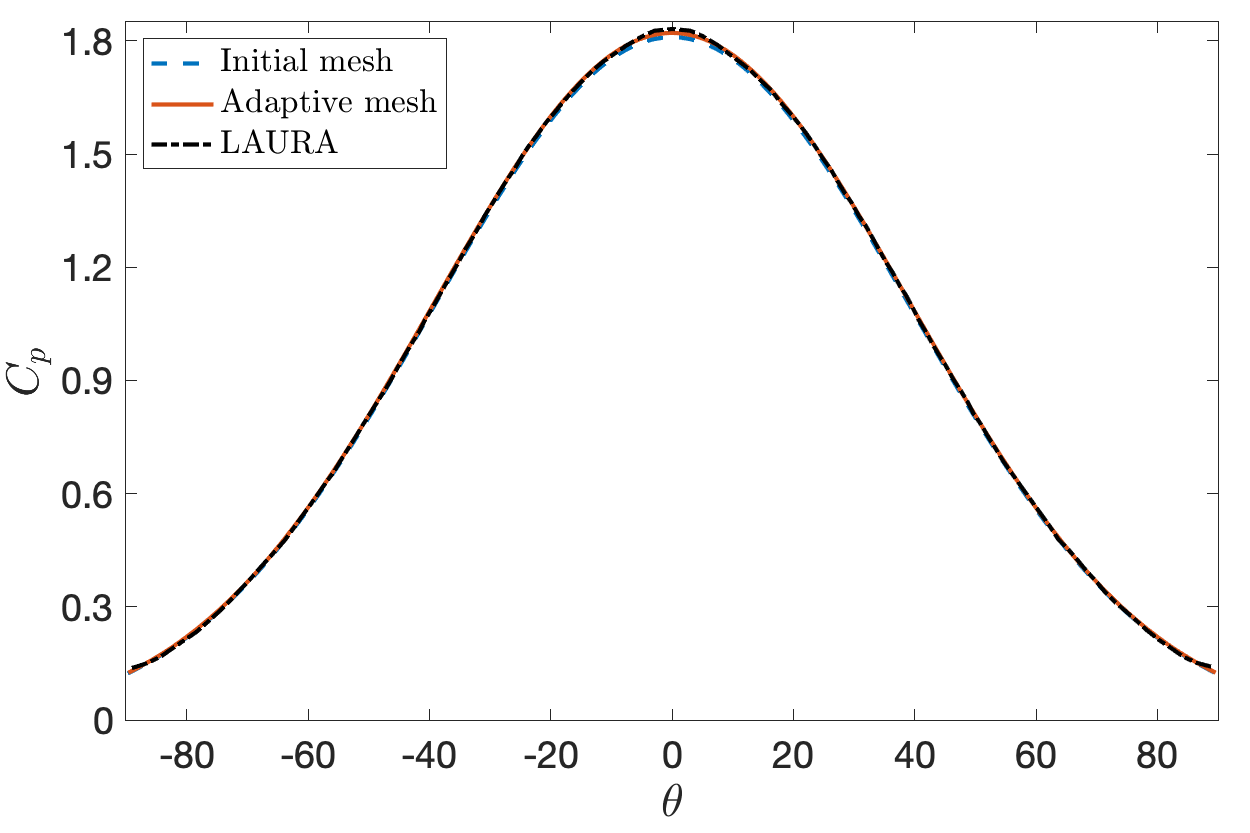}
	\end{subfigure}
	\hfill
	\begin{subfigure}[b]{0.50\textwidth}
		\centering		\includegraphics[width=\textwidth]{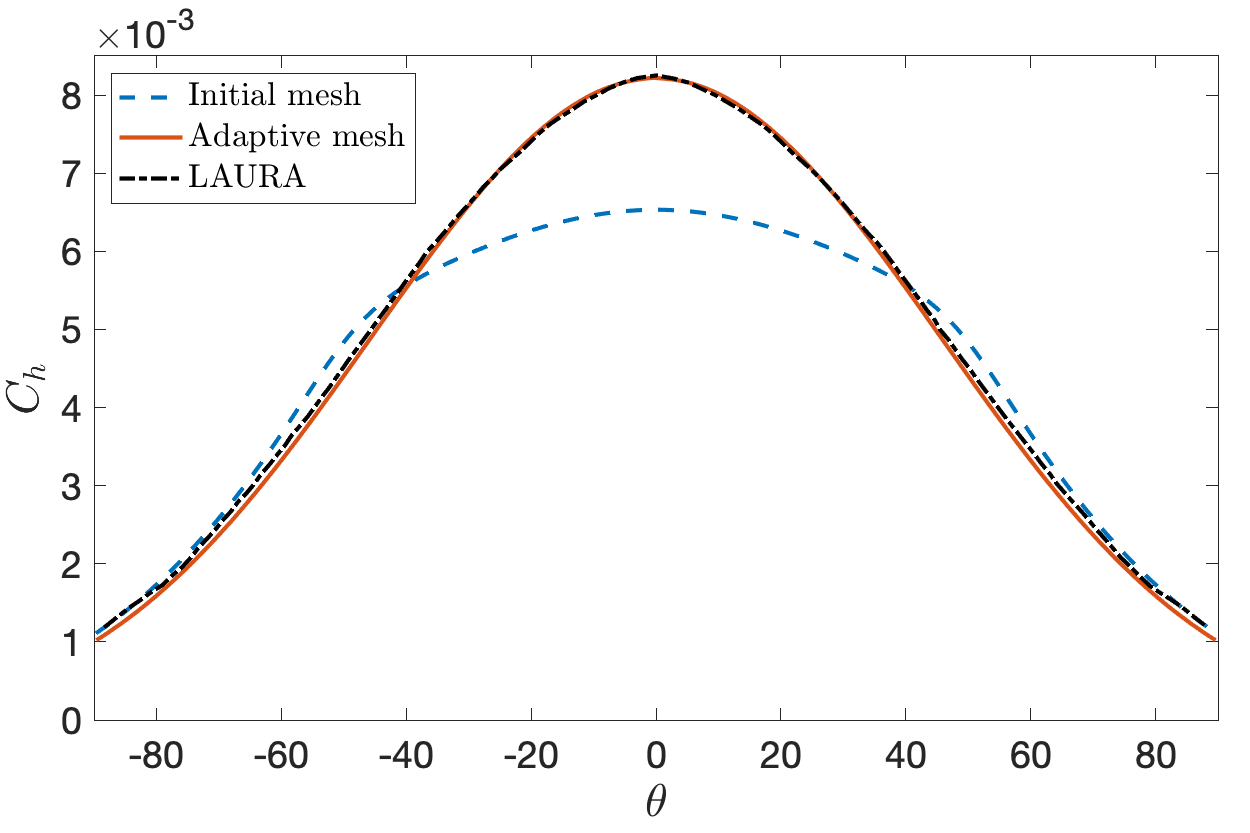}
	\end{subfigure}        
\caption{Pressure coefficient (left) and heat transfer coefficient (right) along the cylinder surface for the viscous hypersonic flow past a circular cylinder. The lines in the third legend correspond to the results obtained using LAURA code by Gnoffo and White \cite{Gnoffo2004}.}
\label{cyl18f}	
\end{figure}

We present the numerical solution computed on the initial mesh in Figure \ref{cyl18c} and on the  adaptive mesh in Figure \ref{cyl18d}. We observe that pressure and temperature rise rapidly behind the bow shock, which create very strong pressure and high temperature environments surrounding the cylinder. In addition, Figure \ref{cyl18e} shows profiles of the numerical solution along the horizontal line $y=0$. We notice that the numerical solution on the adaptive mesh has higher gradient than that on the initial mesh in the shock region and boundary layer. This is because the adaptive mesh has more grid points to resolve those features than the initial mesh. By redistributing the elements of the initial mesh to resolve the bow shock and boundary layer, the optimal transport  can considerably improve the prediction of heating rate as shown in Figure \ref{cyl18f}. We see that while the pressure coefficient on the initial mesh is very similar to that on the adaptive mesh, the heat transfer coefficient on the initial mesh is lower than that on the adaptive mesh. The heat transfer coefficient on the adaptive mesh agrees very well with the prediction by  Gnoffo and White \cite{Gnoffo2004}.

\section{Concluding remarks} \label{sec:conclusions}
We have presented an optimal transport approach for the numerical solution of compressible flows with shock waves. The approach couples an adaptive viscosity regularization method and optimal transport theory in order to capture shocks and adapt meshes. The marriage of optimal transport and viscosity regularization for compressible flows leads to a coupled system of the compressible Euler/Navier-Stokes equations, the Helmholtz equation, and the Monge-Amp\`ere equation. The hybridizable discontinuous Galerkin method is used for the spatial discretization of the governing equations to obtain high-order accurate solutions.  We devise a mesh adaptation  procedure to solve the coupled system in an iterative and sequential fashion.  The approach is found to yield accurate, sharp yet smooth solutions within a few mesh adaptation iterations. We explore two different options to define the mesh indicator function for computing adaptive meshes. The option based on  density gradient is more effective than that based on velocity divergence for dealing with shock flows that have more complex structures such as boundary layers, shear layers, and expansion waves.

We have presented a wide variety of transonic, supersonic and supersonic flows in two dimensions in order to demonstrate the performance of the proposed approach. The approach is capable of moving mesh points to resolve complex shock patterns without creating new mesh points or modifying the  connectivity of the initial mesh. The generated r-adaptive meshes can significantly improve the accuracy of the numerical solution relative to the initial mesh. Accurate prediction of aerodynamic forces and heat transfer rates for viscous shock flows requires meshes to resolve both shocks and boundary layers. Numerical results show that the approach can generate r-adaptive meshes that resolve not only shocks but also boundary layers for viscous shock flows. It yields accurate predictions of pressure and  heat transfer coefficients by adapting the initial mesh to resolve shocks and boundary layers. Moreover, the approach can also adapt the initial mesh to resolve other flow structures such as shear layers and expansion waves.

The approach presented herein can be extended to  chemically reacting hypersonic flows without loss of generality. To this end, different variants of the regularized viscosity can be devised, including physics-based artificial viscosity terms that augment the molecular viscous components. The approach can also be extended to compressible flows in three dimensions. We are going to pursue these extensions in future work.

Another interesting application of the optimal transport approach is model reduction of compressible flows with shock waves. We show in a recent paper \cite{Heyningen2023} that the optimal transport provides an effective treatment of shock waves for model reduction because it can generate snapshots that are aligned well with the shocks. Hence, it results in stable, robust and accurate reduced order models of parametrized compressible flows. In future reserach, we would like to couple the optimal transport theory with the first-order empirical interpolation method \cite{Nguyen2023d} to develop an efficient intrusive reduced order modeling for compressible flows.


\section*{Acknowledgements} \label{}

We gratefully acknowledge the United States  Department of Energy under contract DE-NA0003965, the National Science Foundation for supporting this work (under grant number NSF-PHY-2028125), and the Air Force Office of Scientific Research under Grant No. FA9550-22-1-0356 for supporting this work.


 \bibliographystyle{elsarticle-num} 
\bibliography{library, extrarefs}





\end{document}